\documentclass{amsart}
\usepackage{afterpage}
\usepackage{relsize}

\usepackage{enumerate}
\usepackage{mathdots}
\usepackage{graphicx}
\usepackage{tikz-cd}

\usepackage{booktabs}
\usepackage{tabularx}
\usepackage{array}
\usepackage{microtype}
\usepackage{mathtools}
\usepackage{amssymb}
\usepackage{caption}

\usepackage{tikz}
\usetikzlibrary{arrows.meta,positioning,calc}
\usepackage{pgfplots}
\pgfplotsset{compat=1.18}

\newtheorem{theorem}{Theorem}[section]
\newtheorem{proposition}[theorem]{Proposition}

\newtheorem{corollary}[theorem]{Corollary}

\newcommand{\R}{\mathbb{R}}
\newcommand{\C}{\mathbb{C}}
\newcommand{\tr}{\boldsymbol{\tau}}

\theoremstyle{definition}
\newtheorem{definition}[theorem]{Definition}
\newtheorem{remark}[theorem]{Remark}

\numberwithin{equation}{section} 
  
\begin{document} 

\title{The catastrophes of algebras}

\author{Fred Greensite}

\begin{abstract}
We canonically associate a trace-parametrized family of (pseudo-) Finsler metrics to each real finite-dimensional symmetric Frobenius algebra (and hence to any unital associative algebra via its trivial extension), placing such algebras in the realm of symplectic geometry and Arnol'd's Lagrangian singularity theory. Caustics of the resulting Lagrangian manifolds yield novel computable algebra isomorphism invariants, and any Lagrangian catastrophe of finite Milnor number is synthesized purely from some algebra's internal trace geometry. Equidimensional examples show that these invariants are neither determined by radical-filtration and metric-signature data nor invariant under derived equivalence, distinguishing non-isomorphic algebras with isomorphic Gerstenhaber algebra structures on Hochschild cohomology and isomorphic cyclic (co)homologies. Furthermore, this geometry induces intrinsic \mbox{(pseudo-)Rie}\-mannian metrics on trace-fibered bundles over algebraic deformation spaces. Bypassing coboundary gauge ambiguity and the infinite-dimensional periodic cyclic bar bicomplex, a finite-dimensional variational principle canonically launches trace-parametrized families of maximal analytic Maurer--Cartan paths as geodesic sprays along non-degenerate, trace-compatible unobstructed sheets, determining their trace horizons and making global singularity invariants such as Maslov indices and associated monodromies computationally accessible.
\end{abstract}

\maketitle

\markboth{F. Greensite}{The catastrophes of algebras}

\setcounter{tocdepth}{2}

\tableofcontents

\section{Introduction}
A complete classification of finite-dimensional associative algebras via canonical normal forms or separated moduli spaces is not possible due to wildness of the underlying isomorphism problem \cite{drozd:1979,belitskii:2005}. This indicates that, as a means of furthering understanding of algebra structure, 
 identification of novel classes of algebra invariants independent of prior classes will always be relevant and valuable.  Indeed, this quest has been a central theme of modern algebra, progressing from composition series of the Jordan--H{\"o}lder theorem through the developments of homological algebra and algebraic deformation theory \cite{hochschild:1945,Gerstenhaber:1964,connes:1985,cuntz:1995,Manetti:2004,Dolgushev:2011}. 
In this paper we will fashion a new source of such invariants for real finite-dimensional unital associative algebras by placing such algebras within the framework of symplectic geometry and Lagrangian singularity theory~\cite{arnold:1985}.  This proceeds from our main foundational result that any point in the moduli space canonically carries with it a trace-parametrized family of (pseudo-)Finsler metrics realized via an intrinsic ``algebraic reciprocal transform~(ART)''.    

To illustrate the discriminating power of these symplectic invariants against the backdrop of wildness, we show that any given Lagrangian catastrophe of finite Milnor number can be synthesized purely through some algebra's internal trace geometry, and we exhibit two complementary separation phenomena among equidimensional algebra examples. In the first, we construct real forms of the same complex complete-intersection algebra that share identical radical nilpotency indices, graded Hilbert series, quadratic metric signatures, and Hochschild and cyclic (co)homology dimensions in all degrees, yet via ART realize disjoint umbilic catastrophes ($D_4^+$ versus $D_4^-$). 
In the second, we exhibit equidimensional Brauer tree algebras that are derived equivalent, and therefore possess isomorphic Gerstenhaber algebras $\big(HH^\bullet(A, A), \smile, [\cdot, \cdot]_{\mathrm{G}}\big)$ and isomorphic cyclic (co)homologies across all degrees, yet are separated by the essential cubic ranks of their focal ART phase germs.
Together, these examples demonstrate that the novel invariants introduced in this paper are neither reducible to radical-filtration data nor constrained by derived equivalence (and hence reach beyond the full Gerstenhaber algebra structure on Hochschild cohomology, as well as cyclic (co)homology).

The ART further allows us to devise a finite-dimensional differential-geometric framework to navigate algebraic deformation space that bypasses two classical bottlenecks of homological deformation theory \cite{Gerstenhaber:1964,Manetti:2004,Dolgushev:2011}, these being the non-canonical coboundary gauge choices required at every order when integrating formal Maurer--Cartan series, and the infinite-dimensional $(b, B)$ periodic cyclic bar bicomplex invoked to co-evolve cyclic traces via Getzler's noncommutative Gauss--Manin connection~\cite{getzler:1993} (a bottleneck particularly acute for algebras of infinite global dimension, whose projective bimodule resolutions never terminate).
Specifically, by fibering the admissible trace space over the space of algebra multiplications (and restricting to unobstructed associative deformation sheets), we equip the resulting trace-fibered bundle with an intrinsic (pseudo-)Riemannian metric invariant under unit-preserving changes of basis, allowing families of maximal analytic Maurer--Cartan paths to be launched as geodesic sprays alongside their co-evolving trace calibrations, metric-wall caustics, and Maslov/monodromy invariants via a variational principle.  The respective strengths and limitations of this ART-inspired approach versus classical homological deformation theory are delineated in Remark~\ref{rem:MC_comparison}, showing that the two programs are complementary rather than competitive.

Our primary geometric constructions are formulated for symmetric Frobenius algebras. However, focusing on this class does not restrict the generality of the resulting invariants.  Any algebra $B$ embeds canonically into its symmetric Frobenius trivial extension given by the semidirect product $A = B \ltimes B^*$, so invariants extracted from $A$ serve directly as isomorphism invariants of $B$ (i.e., $B_1 \cong B_2$ implies $B_1 \ltimes B_1^* \cong B_2 \ltimes B_2^*$, although the converse is not always true).

\medskip

The paper proceeds as follows.
\begin{itemize}
\item
Section 2 introduces the ART.  Essentially, we develop an integral transform based on the Maurer-Cartan form that maps a ``normalized admissible trace functional" to the logarithm of a translation-invariant pseudo-Finsler norm. 
\item
Section 3 utilizes the ART to derive a novel class of algebra isomorphism invariants from the invariants of the singularities of Lagrangian manifolds.
\item
Section 4 develops a novel means of exploring algebraic deformation space using ART-inspired methods, in the process adding another new class of algebra invariants from invariants of global singularities of geodesic sprays on pseudo-Riemannian manifolds in the ambient space of algebra multiplications.
 \end{itemize}
 
 \noindent Worked examples are important parts of Section 3 and Section 4.

\subsection*{Notation and conventions} Throughout the paper, we adopt the following conventions: \begin{itemize}
\setlength{\itemsep}{2pt}
\item Unless stated otherwise, the term ``algebra'' refers to a non-zero finite-dimensional unital associative algebra over $\R$ endowed with the standard Euclidean topology.
\item For an algebra $B$, $\mathbf{1}_B$ (or simply $\mathbf{1}$ when clear from context) denotes its multiplicative identity, $\operatorname{rad}(B)$ denotes its Jacobson radical, and $B_0^\times$ denotes the identity connected component of its group of units.
\item $\mathcal{T}(B) \subset B^* = \operatorname{Hom}(B, \R)$ denotes the trace space of $B$, consisting of all linear functionals vanishing on the commutator subspace $[B, B]$.
\item The symbol $A$ is reserved for an $n$-dimensional symmetric Frobenius algebra.
\item The terms ``pseudo-Finsler'' and ``pseudo-Riemannian'' encompass both indefinite and definite Finsler \cite{BaoChernShen2000} and Riemannian signatures.  Following \cite{JavaloyesSanchez2014}, a pseudo-Finsler norm is permitted to have only an open conic domain $\mathcal{U} \subset A \setminus \{0\}$.
\end{itemize}
 
\section{The ART framework}

\begin{definition}\label{def:admissible_traces}
The \emph{parameter space of admissible traces} of $A$ is the set of its Frobenius traces that do not vanish at $\mathbf{1}$,
\begin{equation*}
\mathcal{T}_1(A) := \big\{ \tau \in \mathcal{T}(A) \;\big|\; \tau(\mathbf{1}) \ne 0, \text{ and } \hat{\tau}(x, y) := \tau(xy) \text{ is non-degenerate on } A \big\} . 
\end{equation*}
\end{definition}

Because the left-regular trace satisfies $\operatorname{Tr}_L(\mathbf{1}) = n \neq 0$ (so $\mathbf{1} \notin [A, A]$), the condition $\tau(\mathbf{1}) \neq 0$ defines a non-empty Zariski-open subset of the vector space $\mathcal{T}(A)$. Since non-degeneracy of $\hat{\tau}$ is likewise a non-empty Zariski-open condition on $\mathcal{T}(A)$ (as $A$ is symmetric Frobenius), their intersection $\mathcal{T}_1(A)$ is an open dense (in particular, non-empty) subset of $\mathcal{T}(A)$.

\begin{definition}\label{def:trace_op}
Viewing the tautological inclusion $\tr \in C^\infty\big(\mathcal{T}_1(A), A^*\big)$ as the canonical trace operator parametrized over $\mathcal{T}_1(A)$, its coefficientwise action maps any $A$-valued differential form $z$ to a smooth $\mathcal{T}_1(A)$-parametrized scalar differential form $\tr(z)$ whose evaluation at $\tau \in \mathcal{T}_1(A)$ is given by
\begin{equation}\label{eq:trace_contraction}
\tr(z)\big|_\tau \coloneqq \tau(z).
\end{equation}
\end{definition}

\begin{definition}\label{def:ART}
Let $\mathcal{U} \subset A_0^\times$ be a simply connected open conic neighborhood of units containing the identity ray $\mathbb{R}_{>0}\mathbf{1}$. The \emph{algebraic reciprocal transform (ART)} is the canonical trace-parametrized operator 
\begin{equation}\label{eq:ART_def}
F(y; \tr) \coloneqq \exp\left[ \frac{1}{\tr(\mathbf{1})} \int_{\mathbf{1}}^y \tr\left(t^{-1}dt\right) \right], \qquad y\in \mathcal{U},
\end{equation}
whose evaluation with $\tau\in\mathcal{T}_1(A)$ is the \emph{ART potential} $F(y; \tau)$ resulting from replacement of $\tr$ with $\tau$ in (\ref{eq:ART_def}), so that $F(\cdot\,; \tr) \in C^\infty\big(\mathcal{U} \times \mathcal{T}_1(A)\big)$. 
\end{definition}

Because $\tau$ vanishes on commutators and $d(t^{-1}) = -t^{-1}\,dt\,t^{-1}$, the scalar $1$-form $\omega_\tau \coloneqq \tau(t^{-1}dt)$ satisfies $d\omega_\tau(u, v) = -\tau\big([t^{-1}u, t^{-1}v]\big) = 0$ for all $t \in A_0^\times$ and $u, v \in A$. On the simply connected domain $\mathcal{U}$, $\omega_\tau$ is therefore exact, rendering the path integral in (\ref{eq:ART_def}) independent of the chosen path after evaluation with $\tau \in \mathcal{T}_1(A)$.

Since (\ref{eq:ART_def}) is invariant under non-zero scalar rescalings of $\tau$, it descends to the Zariski-open \emph{projective admissible trace space} $\mathbb{P}(\mathcal{T}_1(A)) \coloneqq \mathcal{T}_1(A)/\R^\times \subset \mathbb{P}(\mathcal{T}(A))$, which the canonical section $[\tau] \mapsto \bar{\tau} \coloneqq \frac{\tau}{\tau(\mathbf{1})}$ identifies diffeomorphically with the affine-hyperplane locus of \emph{normalized admissible traces} (i.e., normalized Frobenius traces) $\{\bar{\tau} \in \mathcal{T}_1(A) : \bar{\tau}(\mathbf{1}) = 1\}$. Theorem~\ref{thm:pseudo_finsler} will thus allow us to interpret (\ref{eq:ART_def}) as mapping normalized admissible traces to conic pseudo-Finsler norms. Use of the word ``transform'' in the nomenclature ART reflects the structure of the exponent in (\ref{eq:ART_def}), which is reminiscent of a classical integral transform.

Formulating geometric and symplectic objects in terms of $\tr$, including the ART and the subsequent ray Hamiltonian, caustic generating family, and deformation space metric, exhibits each as a smoothly parametrized family over the base manifold $\mathcal{T}_1(A)$, evaluated fiberwise at a chosen trace parameter $\tau \in \mathcal{T}_1(A)$.
In the context of singularity and deformation theory, $\mathbb{P}\left(\mathcal{T}_1(A)\right)$ serves as the space of external control parameters that drive caustic bifurcations, unfold local catastrophes, and stratify moduli-space horizons.

\begin{remark}\label{rem:CQ}
Because any $\tau \in \mathcal{T}(A)$ vanishes identically on the commutator subspace $[A, A]$, one might suspect that evaluating the ART potential purely through traces discards essential non-commutative information. This is completely remedied by lifting $A$ to its formal infinitesimal neighborhoods using the Cuntz--Quillen framework~\cite{cuntz:1995}. Specifically, for any presentation $0 \to I \to R \to A \to 0$ of $A$ by a finitely generated free unital algebra $R$ with ideal $I = \ker(R \to A)$, each $k$-th order infinitesimal thickening $A_{(k)} \coloneqq R/I^{k+1}$ ($k \ge 1$) remains finite-dimensional, and evaluating traces across these thickenings (or their symmetric Frobenius trivial extensions) restores sensitivity to $[A, A]$. In fact, because $\dim A < \infty$, this can be performed canonically via the universal presentation $R = T(A)$ on the underlying vector space of $A$ (or $T(A/\R\mathbf{1})$ in the normalized unital setting), which requires no choice of generators and is functorial under algebra isomorphisms. Thus, Definition~\ref{def:ART} naturally generalizes across the tower of thickenings $A_{(k)}$ (via their trivial extensions as necessary), yielding a canonical higher-order ART collection $\{F_{(k)}\}_{k \ge 0}$. In the present work we restrict our attention to the base algebra ($k = 0$), demonstrating that even at this baseline level, the interplay between trace parameters and the nilpotent radical already captures rich families of Lagrangian catastrophes.
\end{remark}

\begin{proposition}\label{prop:WM_factorization}
Given a Wedderburn--Malcev decomposition $A = S \ltimes \operatorname{rad}(A)$, let $\mathcal{U}_S \subset S_0^\times$ be a simply connected open cone containing $\mathbf{1}$. On the simply connected open conic neighborhood $\mathcal{U} = \mathcal{U}_S(\mathbf{1} + \operatorname{rad}(A)) \subset A_0^\times$ of $\R_{>0}\mathbf{1}$, every $y \in \mathcal{U}$ decomposes uniquely as $y = s + r$ with $s \in \mathcal{U}_S$ and $r \in \operatorname{rad}(A)$, and the ART factors canonically as
\begin{equation}\label{eq:WM_factor} 
F(y; \tr) = R_{\boldsymbol{\tau}}(s)e^{P_{\boldsymbol{\tau}}\left(s^{-1}r\right)}. 
\end{equation}
Here, $R_{\boldsymbol{\tau}}(s)$ is a strictly positive real-analytic function on $\mathcal{U}_S$ given by a product of the simple-block determinant polynomial norms of $S$ raised to real exponents linear in $\frac{\tr}{\tr(\mathbf{1})}$, multiplied across any complex simple factors of $S$ by exponential determinant-phase factors whose exponents are linear in $\frac{\tr}{\tr(\mathbf{1})}$.  Furthermore, $P_{\boldsymbol{\tau}}(\rho)$ is an intrinsic polynomial function on $\operatorname{rad}(A)$ without constant term whose homogeneous components are linear in $\frac{\tr}{\tr(\mathbf{1})}$. Finally, the values of both factors $R_{\boldsymbol{\tau}}(s)$ and $P_{\boldsymbol{\tau}}(s^{-1}r)$ at $y$ are independent of the choice of Wedderburn--Malcev complement $S$.
\end{proposition}

\begin{proof}
Fixing a Wedderburn--Malcev decomposition $A = S \ltimes \operatorname{rad}(A)$ and a simply connected open cone $\mathcal{U}_S \subset S_0^\times$ containing $\mathbf{1}$ (for instance, $\mathcal{U}_S = \R_{>0}B_\epsilon(\mathbf{1})$ for a small open convex ball $B_\epsilon(\mathbf{1}) \subset S_0^\times$), every $y \in \mathcal{U} = \mathcal{U}_S(\mathbf{1} + \operatorname{rad}(A))$ factors uniquely as $y = s(\mathbf{1} + \rho)$ with $s \in \mathcal{U}_S$ and $\rho = s^{-1}r \in \operatorname{rad}(A)$, so $\mathcal{U} \cong \mathcal{U}_S \times \operatorname{rad}(A)$ is likewise a simply connected open cone. Following evaluation with $\tau\in\mathcal{T}_1(A)$, we perform the integration indicated in (\ref{eq:ART_def}) along the concatenated path $\gamma_1 * \gamma_2 \subset \mathcal{U}$, where $\gamma_1 \subset \mathcal{U}_S$ connects $\mathbf{1}$ to $s$, and $\gamma_2(v) = s(\mathbf{1} + v\rho)$ for $v \in [0, 1]$ connects $s$ to $y$. Because $\mathcal{U}$ is simply connected, the integral along $\gamma_1 * \gamma_2$ equals the integral along any path in $\mathcal{U}$ from $\mathbf{1}$ to $y$.

By the Wedderburn--Artin theorem, $S \cong \bigoplus_{k=1}^m S_k$, with each simple factor $S_k \cong M_{n_k}(D_k)$ for $D_k \in \{\R, \C, \mathbb{H}\}$. For $D_k \in \{\R, \mathbb{H}\}$, the commutator subspace $[S_k, S_k]$ has real codimension $1$ in $S_k$, so the space of $\R$-linear traces on $S_k$ is one-dimensional and spanned by $\operatorname{Tr}_k(t_k) = \alpha_k \operatorname{Re}\sum_{i=1}^{n_k}(t_k)_{ii}$ ($\alpha_k = 1$ for $D_k = \R$ and $\alpha_k = 2$ for $D_k = \mathbb{H}$), whose logarithmic derivative $\operatorname{Tr}_k(t_k^{-1}dt_k) = d(\log\det_k(t_k))$ generates the positive polynomial determinant norm $\det_k(t_k) > 0$ on $(S_k)_0^\times$ (namely, $\det_\R(t_k)$ for $D_k = \R$ and the Study determinant $\operatorname{Nrd}_{M_{n_k}(\mathbb{H})}(t_k)$ for $D_k = \mathbb{H}$). For $D_k = \C$, $[S_k, S_k] = \mathfrak{sl}_{n_k}(\C)$ has real codimension $2$ in $S_k$, so the space of $\R$-linear traces on $S_k$ is two-dimensional, spanned by $\operatorname{Tr}_k(t_k) \coloneqq 2\operatorname{Re}\operatorname{Tr}_\C(t_k)$ and $\operatorname{Im}\operatorname{Tr}_\C(t_k)$. Setting $\det_k(t_k) \coloneqq |\det_\C(t_k)|^2 > 0$ for $D_k = \C$, the complex Jacobi formula $\operatorname{Tr}_\C(t_k^{-1}dt_k) = d(\log\det_\C(t_k))$ yields both $\operatorname{Tr}_k(t_k^{-1}dt_k) = d(\log\det_k(t_k))$ and $\operatorname{Im}\operatorname{Tr}_\C(t_k^{-1}dt_k) = d(\arg\det_\C(t_k))$. Consequently, on $S = \bigoplus_{k=1}^m S_k$, any $\tau \in \mathcal{T}_1(A)$ restricts to
\begin{equation*}
\tau|_S(t) = \sum_{k=1}^m c_k \operatorname{Tr}_k(t_k) + \sum_{k:\,D_k = \C} d_k \operatorname{Im}\operatorname{Tr}_\C(t_k)
\end{equation*}
for constants $c_k, d_k \in \R$. Integrating $\tau(t^{-1}dt)$ from $\mathbf{1}$ to $s = (s_1, \dots, s_m) \in \mathcal{U}_S$ along $\gamma_1 \subset \mathcal{U}_S$ and dividing by $\tau(\mathbf{1})$ gives
\begin{equation}\label{eq:log_R_tau}
\log R_\tau(s) \coloneqq \frac{1}{\tau(\mathbf{1})}\left[ \sum_{k=1}^m c_k \log\operatorname{det}_k(s_k) + \sum_{k:\,D_k = \C} d_k \operatorname{Arg}_{\mathcal{U}_S}\!\big(\det_\C(s_k)\big) \right], 
\end{equation}
where $\operatorname{Arg}_{\mathcal{U}_S}\!\big(\det_\C(s_k)\big)$ denotes the unique continuous branch of $\arg\det_\C(s_k)$ along $\mathcal{U}_S$ vanishing at $\mathbf{1}$.

Along $\gamma_2$ we have $t^{-1}dt = (\mathbf{1} + v\rho)^{-1}\rho\,dv$. Denoting by $N$ the nilpotency index of $\operatorname{rad}(A)$, we have the terminating geometric series $(\mathbf{1} + v\rho)^{-1} = \sum_{j=0}^{N-1}(-1)^j v^j \rho^j$. Integrating term by term and noting that $\rho^N = 0$ gives
\begin{equation*}
P_{\tau}\left(s^{-1}r\right) \coloneqq
\frac{1}{\tau(\mathbf{1})}\int_{\gamma_2}\tau \left(t^{-1}dt\right) = \frac{1}{\tau(\mathbf{1})}\sum_{j=1}^{N-1}\frac{(-1)^{j-1}}{j}\tau\big((s^{-1}r)^j\big) . 
\end{equation*}
According to Definition~\ref{def:ART}, adding the above to (\ref{eq:log_R_tau}) and exponentiating the result gives $F(y;\tau)$. Finally, by Malcev's uniqueness theorem, any other semisimple complement $\tilde{S}$ takes the form $\tilde{S} = u S u^{-1}$ for some $u \in \mathbf{1} + \operatorname{rad}(A)$. Under $y = \tilde{s} + \tilde{r}$, congruence modulo $\operatorname{rad}(A)$ forces $\tilde{s} = u s u^{-1} \in \mathcal{U}_{\tilde{S}} \coloneqq u\mathcal{U}_S u^{-1} \subset \tilde{S}_0^\times$, and because $\mathbf{1} + \operatorname{rad}(A)$ is normal in $A_0^\times$, we have $(u\mathcal{U}_S u^{-1})(\mathbf{1} + \operatorname{rad}(A)) = \mathcal{U}_S(\mathbf{1} + \operatorname{rad}(A)) = \mathcal{U}$. Both the conjugated path $\tilde{\gamma}_1 = u \gamma_1 u^{-1}$ and the linear segment $\tilde{\gamma}_2(v) = \tilde{s}(\mathbf{1} + v\tilde{s}^{-1}\tilde{r})$ therefore remain in $\mathcal{U}$. Because $\tilde{\gamma}_1^{-1}d\tilde{\gamma}_1 = u(\gamma_1^{-1}d\gamma_1)u^{-1}$ and $\tau$ vanishes on $[A, A]$, trace cyclicity along $u \gamma_1 u^{-1}$ yields $R_\tau(\tilde{s}) = R_\tau(s)$, whence $P_\tau(\tilde{s}^{-1}\tilde{r}) = P_\tau(s^{-1}r)$. Since $\tau\in\mathcal{T}_1(A)$ is arbitrary, the proposition follows.
\end{proof}

\begin{theorem}\label{thm:pseudo_finsler}
For $\tau \in \mathcal{T}_1(A)$, consider the ART potential $F(y; \tau)$ as specified in Definition~\ref{def:ART}. 
Then there exists a simply connected conic open neighborhood $\mathcal{V} \subset \mathcal{U}$ containing the identity ray $\R_{>0}\mathbf{1}$ on which $F(y; \tau)$ is a conic pseudo-Finsler norm, thereby defining a translation-invariant, and hence flat, conic pseudo-Finsler structure on the open subbundle $A \times \mathcal{V} \subset TA \setminus \{0\}$.
\end{theorem}

\begin{proof}
By (\ref{eq:ART_def}), $F(y; \tau) > 0$ is smooth throughout $\mathcal{U}$. Because $\mathcal{U}$ is conic and contains $\R_{>0}\mathbf{1}$, for any $\lambda > 0$ and $y \in \mathcal{U}$, scaling a path in $\mathcal{U}$ from $\mathbf{1}$ to $y$ by $\lambda$ yields a path in $\mathcal{U}$ from $\lambda\mathbf{1}$ to $\lambda y$. Since $t^{-1}dt = (\lambda t)^{-1}d(\lambda t)$, we may evaluate the contracted integral in (\ref{eq:ART_def}) at $\lambda y$ along the concatenation of the linear segment from $\mathbf{1}$ to $\lambda\mathbf{1}$ (parametrized by $u$) with this scaled path from $\lambda\mathbf{1}$ to $\lambda y$. Performing a change of variable on the second path, we obtain
\begin{align*}
\frac{1}{\tau(\mathbf{1})}\int_{\mathbf{1}}^{\lambda y}\tau(t^{-1}dt) 
&= \frac{1}{\tau(\mathbf{1})}\int_1^\lambda \tau\left(u^{-1}\mathbf{1}\,du\right) + \frac{1}{\tau(\mathbf{1})}\int_{\mathbf{1}}^y \tau(t^{-1}dt) \\
&= \log \lambda + \frac{1}{\tau(\mathbf{1})}\int_{\mathbf{1}}^y \tau(t^{-1}dt),
\end{align*}
whence $F(\lambda y; \tau) = \lambda F(y; \tau)$ upon exponentiation. Thus we have degree-1 positive homogeneity on $\mathcal{U}$ (and hence on any conic subset $\mathcal{V} \subset \mathcal{U}$), a first requirement for a pseudo-Finsler norm. 

We next show that $g_{\mu\nu} := \frac{1}{2} \partial_\nu \partial_\mu (F^2)$ is non-degenerate on a simply connected conic open neighborhood $\mathcal{V} \subset \mathcal{U}$. 
Let 
\begin{equation}\label{eq:Psi_def}
\Psi(y; \tau) := \log F(y; \tau) = \frac{1}{\tau(\mathbf{1})}\int_{\mathbf{1}}^y \tau\left(t^{-1}dt\right).
\end{equation}
Because $\tau$ is linear and the inversion map $t \mapsto t^{-1}$ is smooth on $A_0^\times$, the differential $1$-form $d\Psi|_y \in T_y^* A$ evaluated on a tangent vector $v \in A$ is given by the integrand,
\begin{equation*}
d\Psi|_y(v) = \frac{1}{\tau(\mathbf{1})} \tau\left(y^{-1}v\right).
\end{equation*}
Expressing $y = y^\mu e_\mu$ in a linear basis $\{e_\mu\}$ (summing over repeated indices), the first partial derivatives of $\Psi$ are obtained from (\ref{eq:Psi_def}) as
\begin{equation}\label{eq:dPsi}
\partial_\mu \Psi(y) = \frac{1}{\tau(\mathbf{1})}\tau\left(y^{-1}e_\mu\right).
\end{equation}
Evaluating at $y = \mathbf{1}$ yields
\begin{equation}\label{eq:dPsi_one}
\partial_\mu \Psi(\mathbf{1}) = \frac{\tau(e_\mu)}{\tau(\mathbf{1})}.
\end{equation}

To determine the second partial derivatives, we exploit the identity $d(y y^{-1}) = 0$. Along the coordinate direction $e_\nu$ we then have $e_\nu y^{-1} + y \, \partial_\nu(y^{-1}) = 0$, so that
\begin{equation}\label{eq:d_inv}
 \partial_\nu(y^{-1}) = -y^{-1} e_\nu y^{-1}.
\end{equation}
Equations (\ref{eq:dPsi}) and (\ref{eq:d_inv}) then imply
\begin{equation*}
\partial_\nu \partial_\mu \Psi(y) = \frac{1}{\tau(\mathbf{1})}\tau\left(\partial_\nu(y^{-1}) e_\mu\right) = -\frac{1}{\tau(\mathbf{1})}\tau\left(y^{-1}e_\nu y^{-1}e_\mu\right).
\end{equation*}
Evaluating at $y = \mathbf{1}$ we obtain
\begin{equation}\label{eq:ddPsi_one}
\partial_\nu \partial_\mu \Psi(\mathbf{1}) = -\frac{1}{\tau(\mathbf{1})}\tau(e_\nu e_\mu) = -\frac{1}{\tau(\mathbf{1})}\tau(e_\mu e_\nu).
\end{equation}

Now we use the above to evaluate $\frac{1}{2} \partial_\nu \partial_\mu (F^2)$. Expressing $F^2 = \exp(2\Psi)$, the first derivative with respect to $y^\mu$ is
\begin{equation*}
\partial_\mu (F^2) = 2 \exp(2\Psi) \, \partial_\mu \Psi = 2 F^2 \partial_\mu \Psi.
\end{equation*}
Differentiating again, we obtain
\begin{equation}\label{eq:g_munu}
g_{\mu\nu}(y; \tau) = 2 F(y; \tau)^2 (\partial_\mu \Psi)(\partial_\nu \Psi) + F(y; \tau)^2 (\partial_\nu \partial_\mu \Psi).
\end{equation}
Substituting \eqref{eq:dPsi_one}, \eqref{eq:ddPsi_one}, and $F(\mathbf{1}; \tau) = 1$ (from (\ref{eq:ART_def})) into~\eqref{eq:g_munu} gives
\begin{equation*}
g_{\mu\nu}(\mathbf{1}; \tau) = \frac{2}{\tau(\mathbf{1})^2}\tau(e_\mu)\tau(e_\nu) - \frac{1}{\tau(\mathbf{1})}\tau(e_\mu e_\nu).
\end{equation*}
Evaluating with arbitrary tangent vectors $u = u^\mu e_\mu$ and $v = v^\nu e_\nu$ produces the coordinate-free expression for $g_{\mu\nu}$ at the identity $\mathbf{1}\in A$, notated as \begin{equation}\label{eq:g_one}
g_{\mathbf{1}}(u, v; \tau) = \frac{2}{\tau(\mathbf{1})^2}\tau(u)\tau(v) - \frac{1}{\tau(\mathbf{1})}\tau(uv).
\end{equation}

We next establish that the symmetric bilinear form $g_{\mathbf{1}}(u, v; \tau)$ is non-degenerate. Indeed, the latter is degenerate if and only if we can find $u\ne 0$ such that $g_{\mathbf{1}}(u, v; \tau) = 0$ for all $v \in A$. In that case, for all $v\in A$, (\ref{eq:g_one}) implies
\begin{equation}\label{eq:g_one_deg1}
\tau(uv) = \frac{2\tau(u)}{\tau(\mathbf{1})}\tau(v) = \tau\left( \frac{2\tau(u)}{\tau(\mathbf{1})}\mathbf{1}  v \right).
\end{equation}
Equivalently, $\tau\left( \left(u - \frac{2\tau(u)}{\tau(\mathbf{1})}\mathbf{1}\right) v \right) = 0$ for all $v \in A$. Since we have $\tau\in \mathcal{T}_1(A)$, the associated bilinear form $\hat{\tau}(x,y)\coloneqq \tau(xy)$ is non-degenerate and $\tau(\mathbf{1})\ne 0$, so that the equation in the prior sentence implies 
\begin{equation}
\label{eq:g_one_deg2}
u = \frac{2\tau(u)}{\tau(\mathbf{1})}\mathbf{1}.
\end{equation}
Applying $\tau$ to both sides gives $\tau(u) = 2\tau(u)$, which implies $\tau(u) = 0$. Substituting this into (\ref{eq:g_one_deg2}) yields $u = 0$. Thus, $g_{\mathbf{1}}(u, v; \tau)$ is non-degenerate.

Smoothness and degree-$1$ positive homogeneity of $F(y;\tau)$ on $\mathcal{U}$ imply that $g_{\mu\nu}(y;\tau) = \frac{1}{2}\partial_\nu\partial_\mu\big(F(y;\tau)^2\big)$ is smooth and positively homogeneous of degree zero in $y$ throughout $\mathcal{U}$. Since $g_{\mu\nu}(y; \tau)$ has thus far been shown non-degenerate only at $y = \mathbf{1}$ (and may degenerate elsewhere in $\mathcal{U}$), we invoke continuity at $\mathbf{1}$ to choose an open convex ball $B_\epsilon(\mathbf{1}) \subset \mathcal{U}$ (with $0 \notin B_\epsilon(\mathbf{1})$) on which $g_{\mu\nu}(y; \tau)$ remains non-degenerate. Because $\mathcal{U}$ is conic, setting $\mathcal{V} \coloneqq \mathbb{R}_{>0} B_\epsilon(\mathbf{1})$ yields a convex (and hence simply connected) conic open neighborhood $\mathcal{V} \subset \mathcal{U}$ containing $\mathbb{R}_{>0}\mathbf{1}$ on which $g_{\mu\nu}(y;\tau)$ is everywhere non-degenerate by degree-zero homogeneity, endowing $A \times \mathcal{V} \subset TA \setminus \{0\}$ with a flat conic pseudo-Finsler structure.
\end{proof}

Because the ART is constructed without reference to coordinates or choices of basis, any unital algebra isomorphism $\varphi \colon A_1 \to A_2$ induces a linear isomorphism $\varphi^* \colon \mathcal{T}(A_2) \to \mathcal{T}(A_1)$ restricting to a bijection $\varphi^* \colon \mathcal{T}_1(A_2) \to \mathcal{T}_1(A_1)$ that satisfies $F_{A_2}(\varphi(y); \tau) = F_{A_1}(y; \varphi^*\tau)$ and $g_{\mathbf{1}}^{A_2}(\varphi(u), \varphi(v); \tau) = g_{\mathbf{1}}^{A_1}(u, v; \varphi^*\tau)$. Simultaneously, its cotangent lift $T^*\varphi^{-1} \colon T^*A_1 \to T^*A_2$ is a fiber-preserving symplectomorphism and hence a Lagrangian equivalence with respect to the cotangent bundle projections $T^*A_i \to A_i$, that intertwines the induced pseudo-Finsler structures across corresponding admissible traces.

\begin{corollary}\label{cor:iso_invariance}
Any trace-parametrized family of Lagrangian submanifolds of the cotangent bundle constructed canonically from the ART pseudo-Finsler geometry is intertwined via Lagrangian equivalence under unital algebra isomorphisms. Consequently, the resulting caustic singularity classes and their trace-parametrized unfolding invariants are algebra isomorphism invariants.
\end{corollary}

\section{Catastrophe synthesis}
\label{sec:catastrophes}

By Theorem~\ref{thm:pseudo_finsler}, each $\tau \in \mathcal{T}_1(A)$ endows a conic neighborhood of $\R_{>0}\mathbf{1} \subset A$ with a translation-invariant pseudo-Finsler norm $F(y;\tau)$. Via the fiberwise Legendre transform of the kinetic Lagrangian $\mathcal{L}(y;\tau) = \frac{1}{2}F(y;\tau)^2$ near $y = \mathbf{1}$, which is a local diffeomorphism onto a neighborhood of $\mathbf{1}^*_\tau \coloneqq d\mathcal{L}|_{\mathbf{1}}(\tau) = \frac{1}{\tau(\mathbf{1})}\tau \in A^*$ with Hessian $g_{\mathbf{1}}(\tau)$, unit-speed wavefront propagation on $T^*A$ is governed by the dual co-indicatrix
\begin{equation}\label{eq:co_indicatrix}
F^*(\mathbf{p}; \tr) = 1. 
\end{equation}

To probe local Lagrangian singularities along a $\kappa$-dimensional subspace $V \subset \ker\tau^*$  on which $g_{\mathbf{1},\perp}(\tau^*) \coloneqq g_{\mathbf{1}}(\cdot,\cdot;\tau^*)|_{V \times V}$ is non-degenerate for a reference trace $\tau^* \in \mathcal{T}_1(A)$, which is called an \emph{active transverse channel}, let $\mathcal{O}_V \subset \{\tau \in \mathcal{T}_1(A) : \tau|_V = 0\}$ be an open neighborhood of $\tau^*$ on which $g_{\mathbf{1},\perp}(\tau)$ remains non-degenerate.  Because $g_{\mathbf{1}}(\mathbf{1},\mathbf{1};\tau) = 1$ and $\mathbf{1} \perp_{g_{\mathbf{1}}(\tau)} \ker\tau$, each $\tau \in \mathcal{O}_V$ induces the $g_{\mathbf{1}}(\tau)$-orthogonal splitting
\begin{equation}\label{eq:orth_split}
A = \R\mathbf{1} \oplus V \oplus W(\tau), \qquad W(\tau) \coloneqq (\R\mathbf{1} \oplus V)^{\perp_{g_{\mathbf{1}}(\tau)}}. 
\end{equation}
Restricting momenta near $\mathbf{1}^*_\tau$ to the annihilator $W(\tau)^0 \subset A^*$ by writing $\mathbf{p} = (1 - h)\mathbf{1}^*_\tau + \boldsymbol{\sigma}$ with $h \in \R$ near $0$ and $\boldsymbol{\sigma} \in V^*$ near $\mathbf{0}$ (extended by zero on $\R\mathbf{1} \oplus W(\tau)$), the squared co-indicatrix condition $F^*\big((1 - h)\mathbf{1}^*_\tau + \boldsymbol{\sigma}; \tau\big)^2 = 1$ expands across (\ref{eq:orth_split}) as
\begin{equation*}
(1 - h)^2 + g_{\mathbf{1},\perp}(\tau)^{-1}(\boldsymbol{\sigma}, \boldsymbol{\sigma}) + O\big(|h|^3 + |\boldsymbol{\sigma}|^3\big) = 1. 
\end{equation*}
Canceling $1$ from both sides leaves $2h = g_{\mathbf{1},\perp}(\tau)^{-1}(\boldsymbol{\sigma}, \boldsymbol{\sigma}) + O(|h|^2 + |\boldsymbol{\sigma}|^3)$, so by the Implicit Function Theorem there is a unique smooth function $H_\perp(\boldsymbol{\sigma}; \tau)$ near $\boldsymbol{\sigma} = \mathbf{0}$ with $H_\perp(\mathbf{0}; \tau) = 0$ such that $\mathbf{p} = \big(1 - H_\perp(\boldsymbol{\sigma}; \tau)\big)\mathbf{1}^*_\tau + \boldsymbol{\sigma}$ satisfies (\ref{eq:co_indicatrix}), given to quadratic order by
\begin{equation}\label{eq:H_perp}
H_\perp(\boldsymbol{\sigma}; \tr) = \frac{1}{2}g_{\mathbf{1},\perp}(\tr)^{-1}(\boldsymbol{\sigma}, \boldsymbol{\sigma}) + O(|\boldsymbol{\sigma}|^3). 
\end{equation}
Via symplectic reduction along $\mathbf{p}|_{W(\tau)} = 0$ onto $T^*(\R\mathbf{1} \oplus V)$, evaluating $\mathbf{p}(z\mathbf{1} + x_\perp) = z + \langle x_\perp, \boldsymbol{\sigma} \rangle - z H_\perp(\boldsymbol{\sigma}; \tau)$ for $(z, x_\perp) \in \R \times V$, adding an initial wavefront phase $S_0(\boldsymbol{\sigma})$ at $z = 0$, and dropping the leading $\boldsymbol{\sigma}$-independent offset summand $z$ on the right-hand side, yields the Hamilton--Jacobi generating family~\cite{arnold:1985,BerryUpstill}
\begin{equation}\label{eq:paraxial_phase}
\Phi(x_\perp, z, \tr; \boldsymbol{\sigma}) = \langle x_\perp, \boldsymbol{\sigma} \rangle + S_0(\boldsymbol{\sigma}) - z H_\perp(\boldsymbol{\sigma}; \tr), 
\end{equation}
whose fiber-critical locus $\partial_{\boldsymbol{\sigma}}\Phi = 0$ generates the evolved Lagrangian wavefront with transverse Hamiltonian $H_\perp$. Fixing the canonical quadratic probe $S_0(\boldsymbol{\sigma}) = \frac{1}{2}k_0 Q_{\tau^*}(\boldsymbol{\sigma})$ ($k_0 \neq 0$), where $Q_{\tau^*}(\boldsymbol{\sigma}) \coloneqq g_{\mathbf{1},\perp}(\tau^*)^{-1}(\boldsymbol{\sigma}, \boldsymbol{\sigma})$, cancels the quadratic Hessian of $\Phi$ at the focal point $\mathbf{c}^* = (\mathbf{0}, z^*, \tau^*)$ with $z^* = k_0$. By Corollary~\ref{cor:iso_invariance}, the right-equivalence class of the focal germ $\Phi(\mathbf{c}^*; \cdot)$ and the Lagrangian ($R^+$-)equivalence class of the unfolding $\Phi$ across admissible pairs $(\tau^*, V)$ (for fixed $k_0 \neq 0$) are algebra isomorphism invariants.

\subsection{Universal catastrophe realizability}
We now prove that for any given local catastrophe class there is an algebra such that the catastrophe is accessible purely through the algebra's internal trace geometry under a canonical quadratic incident probe (i.e., without introducing an engineered higher-order incident phase profile).

The essence of the proof is straightforward. We define an algebra (which happens not to be symmetric Frobenius) with as many nilpotent generators of specified nilpotency index as required, and pass to its symmetric Frobenius trivial extension (so that the ambient Legendre transform near $\mathbf{1}_A$ is non-singular). Because of the nilpotent polynomial exponent in (\ref{eq:WM_factor}) of Proposition~\ref{prop:WM_factorization}, arbitrary multinomial cross-terms appear on the dual co-indicatrix in momentum space and can be independently tuned via the trace controls. The desired catastrophes then follow from standard catastrophe theory principles.

\begin{theorem}\label{thm:universal_catastrophe}
Let $\mathcal{K}$ be an Arnol'd Lagrangian singularity class of corank $\kappa \ge 1$ with finite Milnor number $\mu \ge 2$, and fix $k_0 \neq 0$. Then there exist a symmetric Frobenius algebra $A$, a reference trace $\tau^* \in \mathcal{T}_1(A)$, and a $\kappa$-dimensional $g_{\mathbf{1}}(\tau^*)$-non-degenerate active transverse channel $V \subset \ker\tau^*$ such that, under the canonical quadratic probe $S_0(\boldsymbol{\sigma}) = \frac{1}{2}k_0 Q_{\tau^*}(\boldsymbol{\sigma})$, the ART phase family~(\ref{eq:paraxial_phase}) at $\mathbf{c}^* = (\mathbf{0}, k_0, \tau^*)$ realizes $\mathcal{K}$ and its versal unfolding over $V \times \R \times \mathcal{O}_V$.
\end{theorem}

\begin{proof}
Let $\mathfrak{m}_{\boldsymbol{\sigma}}$ denote the maximal ideal of smooth function germs on $V^* \cong \mathbb{R}^\kappa$ vanishing at $\boldsymbol{\sigma} = \mathbf{0}$, and represent $\mathcal{K}$ by a polynomial normal form $P_0(\boldsymbol{\sigma}) \in \mathfrak{m}_{\boldsymbol{\sigma}}^3$ with Milnor number $\mu = \dim Q_{P_0}$, where $Q_{P_0} = C^\infty_{\mathbf{0}}(V^*) / \langle \nabla P_0 \rangle$ is the local Milnor algebra~\cite{arnold:1985}; by Mather's determinacy theorem~\cite{Mather1968}, $P_0$ is $k$-determined for any fixed $k \ge \mu + 1 \ge 3$. Set
\begin{equation*}
B = \mathbb{R}[\mathbf{e}_1, \dots, \mathbf{e}_\kappa] / \mathfrak{m}_{\mathbf{e}}^{k+1} \cong \mathbb{R}\mathbf{1}_B \oplus \bigoplus_{m=1}^k S^m(V), 
\end{equation*}
where $V = \operatorname{span}_{\mathbb{R}}\{\mathbf{e}_1, \dots, \mathbf{e}_\kappa\}$, so that $\operatorname{rad}(B)^{k+1} = \{0\}$ while $S^2(V) \neq \{0\}$ since $k \ge 3$. Because $\dim\operatorname{soc}(B) = \dim S^k(V) > 1$ renders every trace pairing degenerate on $B$ whenever $\kappa \ge 2$, we pass uniformly for all $\kappa \ge 1$ to its symmetric Frobenius trivial extension $A = B \ltimes B^*$ with multiplication $(b_1, \phi_1)(b_2, \phi_2) = (b_1 b_2, b_1 \cdot \phi_2 + \phi_1 \cdot b_2)$, where $(b \cdot \phi)(x) = \phi(xb)$ and $(\phi \cdot b)(x) = \phi(bx)$. Every $\tau_B \in B^*$ with $\tau_B(\mathbf{1}_B) \neq 0$ lifts canonically to $\tau \in \mathcal{T}_1(A)$ via
\begin{equation}\label{eq:triv_ext_trace}
\tau(b, \phi) \coloneqq \phi(\mathbf{1}_B) + \tau_B(b), 
\end{equation}
whose bilinear pairing
\begin{equation}\label{eq:triv_ext_pairing}
\hat{\tau}\big((b_1, \phi_1), (b_2, \phi_2)\big) = \phi_1(b_2) + \phi_2(b_1) + \tau_B(b_1 b_2)
\end{equation}
is globally non-degenerate on $A$ via its off-diagonal evaluation block. Embedding $V \cong V \times \{0\} \subset \operatorname{rad}(A)$, we compute along the smooth submanifold $\mathcal{T}_{1,V}(A) \subset \{\tau \in \mathcal{T}_1(A) : \tau|_V = 0\}$ of canonical lifts (\ref{eq:triv_ext_trace}) satisfying $\tau_B|_V = 0$ and $\det g_{\mathbf{1},\perp}(\tau) \neq 0$.

Because $B$ is commutative, for any unit $t = (b, \phi) \in A_0^\times$, $t^{-1} = (b^{-1}, -b^{-2}\cdot\phi)$ gives $t^{-1}dt = \big(b^{-1}db, \, d(b^{-1}\cdot\phi)\big)$. Integrating from $\mathbf{1}_A = (\mathbf{1}_B, 0)$ to $(b, \phi)$ via (\ref{eq:Psi_def}) and (\ref{eq:triv_ext_trace}) yields
\begin{equation}\label{eq:Psi_triv_ext}
\Psi(b, \phi; \tau) = \frac{1}{\tau_B(\mathbf{1}_B)}\left[ \int_{\mathbf{1}_B}^b \tau_B(s^{-1}ds) + \phi(b^{-1}) \right]. 
\end{equation}
By Euler's theorem (since $d\mathcal{L} = F^2 d\Psi$ and $F = e^\Psi = 1$), $\mathbf{p} \in A^*$ near $\mathbf{1}^*_\tau$ satisfies $F^*(\mathbf{p}; \tau) = 1$ if and only if $\mathbf{p} = d\Psi|_{(b,\phi)}$ at a point $(b, \phi)$ with $\Psi(b, \phi; \tau) = 0$. For $\tau \in \mathcal{T}_{1,V}(A)$ and $\mathbf{p} = (1 - h)\mathbf{1}^*_\tau + \boldsymbol{\sigma} \in W(\tau)^0$ ($h \in \R$ near $0$, $\boldsymbol{\sigma} \in V^*$ near $\mathbf{0}$), setting $u_{\boldsymbol{\sigma}} \coloneqq g_{\mathbf{1},\perp}(\tau)^{-1}\boldsymbol{\sigma} \in V$ gives $\mathbf{p} = g_{\mathbf{1}}\big((1 - h)\mathbf{1}_A + (u_{\boldsymbol{\sigma}}, 0), \cdot; \tau\big)$ across (\ref{eq:orth_split}), so using $\tau_B(u_{\boldsymbol{\sigma}}) = 0$ and (\ref{eq:g_one}) yields
\begin{multline*}
\mathbf{p}(\delta b, \delta\phi) = g_{\mathbf{1}}\big((1 - h)\mathbf{1}_A + (u_{\boldsymbol{\sigma}}, 0), (\delta b, \delta\phi); \tau\big) \\ = \frac{1}{\tau_B(\mathbf{1}_B)}\Big[ \tau_B\big(((1 - h)\mathbf{1}_B - u_{\boldsymbol{\sigma}})\delta b\big) + \delta\phi\big((1 - h)\mathbf{1}_B - u_{\boldsymbol{\sigma}}\big) \Big]. 
\end{multline*}
Matching this against $d\Psi|_{(b,\phi)}(\delta b, \delta\phi) = \frac{1}{\tau_B(\mathbf{1}_B)}\big[ \tau_B(b^{-1}\delta b) - \phi(b^{-2}\delta b) + \delta\phi(b^{-1}) \big]$ in $\delta\phi \in B^*$ forces $b^{-1} = (1 - h)\mathbf{1}_B - u_{\boldsymbol{\sigma}}$, whereupon the $\delta b \in B$ variation yields $\phi(b^{-2}\delta b) = 0$ and hence $\phi = 0$. Substituting $b = ((1 - h)\mathbf{1}_B - u_{\boldsymbol{\sigma}})^{-1}$ and $\phi = 0$ into $\Psi(b, \phi; \tau) = 0$, expanding the terminating logarithm on $\operatorname{rad}(B)$ ($\operatorname{rad}(B)^{k+1} = \{0\}$), and using $g_{\mathbf{1},\perp}(u_1, u_2; \tau) = -\frac{1}{\tau_B(\mathbf{1}_B)}\tau_B(u_1 u_2)$ on $V$, the dispersion condition $F^*\big((1 - h)\mathbf{1}^*_\tau + \boldsymbol{\sigma}; \tau\big) = 1$ reduces in closed form to
\begin{equation}\label{eq:triv_ext_dispersion}
\log(1 - h) + \frac{1}{2(1 - h)^2} g_{\mathbf{1},\perp}(\tau)^{-1}(\boldsymbol{\sigma}, \boldsymbol{\sigma}) - \sum_{m=3}^k \frac{1}{m(1 - h)^m \, \tau(\mathbf{1}_A)} \tau\Big( \big(g_{\mathbf{1},\perp}(\tau)^{-1} \boldsymbol{\sigma}\big)^m \Big) = 0. 
\end{equation}

Solving (\ref{eq:triv_ext_dispersion}) near $(h, \boldsymbol{\sigma}) = (0, \mathbf{0})$ for $h = H_\perp(\boldsymbol{\sigma}; \tau)$ by recursive series reversion through degree $k$ gives
\begin{multline}\label{eq:H_perp_reversion}
H_\perp(\boldsymbol{\sigma}; \tau) = \frac{1}{2} g_{\mathbf{1},\perp}(\tau)^{-1}(\boldsymbol{\sigma}, \boldsymbol{\sigma}) \\ - \sum_{m=3}^k \left[ \frac{1}{m \, \tau(\mathbf{1}_A)} \tau\Big( \big(g_{\mathbf{1},\perp}(\tau)^{-1}\boldsymbol{\sigma}\big)^m \Big) + \widetilde{\mathcal{R}}_m(\boldsymbol{\sigma}; \tau) \right] + O(|\boldsymbol{\sigma}|^{k+1}), 
\end{multline}
where each degree-$m$ homogeneous polynomial $\widetilde{\mathcal{R}}_m(\boldsymbol{\sigma}; \tau)$ depends polynomially only on $g_{\mathbf{1},\perp}(\tau)^{-1}$ and lower-degree trace moments on $\bigoplus_{j=2}^{m-1} S^j(V)$. Because the symmetric trace moments $\{\tau_B(\mathbf{e}_1^{a_1}\cdots\mathbf{e}_\kappa^{a_\kappa})\}_{a_1+\dots+a_\kappa = m}$ provide independent linear coordinates on $S^m(V)^*$ and $g_{\mathbf{1},\perp}(\tau)^{-1} \colon V^* \to V$ is a linear isomorphism, (\ref{eq:H_perp_reversion}) implies that on the slice $\{\tau \in \mathcal{T}_{1,V}(A) : \tau_B(\mathbf{1}_B) = 1\}$, the truncated $k$-jet map $\tau \mapsto j^k_{\mathbf{0}} H_\perp(\cdot; \tau) \in \mathfrak{m}_{\boldsymbol{\sigma}}^2 / \mathfrak{m}_{\boldsymbol{\sigma}}^{k+1}$ is a triangular smooth diffeomorphism onto the space of quadratic-to-degree-$k$ polynomials with non-degenerate quadratic form.

Now set $z^* = k_0 \neq 0$, fix $\tau_B^*(\mathbf{1}_B) = 1$, and choose $\tau_B^*(\mathbf{e}_i \mathbf{e}_j) = -\varepsilon_i \delta_{ij}$ ($\varepsilon_i \in \{\pm 1\}$), identifying $V^* \cong \R^\kappa$ via the dual basis $\{\mathbf{e}_i^*\}$, so that $g_{\mathbf{1},\perp}(\mathbf{e}_i, \mathbf{e}_j; \tau^*) = \varepsilon_i \delta_{ij}$ and the quadratic part of $z^* H_\perp(\boldsymbol{\sigma}; \tau^*)$ cancels $S_0(\boldsymbol{\sigma}) = \frac{1}{2}k_0 Q_{\tau^*}(\boldsymbol{\sigma})$ in $\Phi(\mathbf{0}, z^*, \tau^*; \boldsymbol{\sigma})$. By triangular surjectivity of the $k$-jet map in (\ref{eq:H_perp_reversion}), sequential back-substitution across $m = 3, \dots, k$ uniquely specifies $\tau_B^*|_{S^m(V)}$ so that the central phase germ $f(\boldsymbol{\sigma}) \coloneqq \Phi(\mathbf{0}, z^*, \tau^*; \boldsymbol{\sigma}) = P_0(\boldsymbol{\sigma}) + O(|\boldsymbol{\sigma}|^{k+1})$ matches $P_0$ through degree $k$; since $P_0$ is $k$-determined, $f$ is smoothly right-equivalent to $P_0$. Finally, on an open neighborhood $\mathcal{O}_V \subset \{\tau \in \mathcal{T}_1(A) : \tau|_V = 0\}$ of $\tau^*$ where $g_{\mathbf{1},\perp}(\tau)$ remains non-degenerate, varying $x_\perp \in V$ at $\mathbf{c}^* = (\mathbf{0}, z^*, \tau^*)$ yields $\partial_{x_{\perp,i}}\Phi|_{\mathbf{c}^*} = \sigma_i$ (spanning $\mathfrak{m}_{\boldsymbol{\sigma}} / \mathfrak{m}_{\boldsymbol{\sigma}}^2$), while holding $S_0$ fixed at $\tau^*$ and varying $\tau$ along $\mathcal{T}_{1,V}(A) \cap \mathcal{O}_V$ maps $T_{\tau^*}\mathcal{T}_{1,V}(A)$ surjectively onto $\mathfrak{m}_{\boldsymbol{\sigma}}^2 / \mathfrak{m}_{\boldsymbol{\sigma}}^{k+1}$ because $\tau \mapsto j^k_{\mathbf{0}} H_\perp(\cdot; \tau)$ is a local diffeomorphism at $\tau^*$ and $z^* \neq 0$. Thus the parameter derivatives of $\Phi$ over $V \times \R \times \mathcal{O}_V$ span $\mathfrak{m}_{\boldsymbol{\sigma}} / \mathfrak{m}_{\boldsymbol{\sigma}}^{k+1}$ and hence surject onto the reduced Milnor algebra $\widetilde{Q}_f = \mathfrak{m}_{\boldsymbol{\sigma}} / \langle \nabla f \rangle$ (as $\mathfrak{m}_{\boldsymbol{\sigma}}^{k+1} \subset \langle \nabla f \rangle$ by $k$-determinacy), proving $R^+$-versality (and hence Lagrangian versality~\cite{arnold:1985}) by the Mather--Martinet theorem~\cite{Mather1968, Martinet1982}.
\end{proof}

While the commutative local algebra $B$ provides a convenient vehicle for proving Theorem~\ref{thm:universal_catastrophe}, neither commutativity nor locality is required once $\operatorname{rad}(A) \neq \{0\}$. In Section~\ref{sec:quiver_example} below, we explicitly illustrate this by realizing and versally unfolding both umbilic catastrophes within the trivial extension of an indecomposable, noncommutative, non-local bound quiver algebra.

\begin{corollary}\label{cor:simple_obstruction}
Let $A$ be a simple finite-dimensional real associative algebra, so that $A \cong M_n(D)$ with $n \ge 1$ and $D \in \{\R, \C, \mathbb{H}\}$ by the Wedderburn and Frobenius theorems. Then $\operatorname{rad}(A) = \{0\}$ and the projective admissible trace space $\mathbb{P}(\mathcal{T}_1(A))$ has dimension $0$ when $D \in \{\R, \mathbb{H}\}$ and $1$ when $D = \C$, independent of $n$. Consequently, under the canonical quadratic probe $S_0(\boldsymbol{\sigma}) = \frac{1}{2}k_0 Q_{\tau^*}(\boldsymbol{\sigma})$ ($k_0 \neq 0$) at any $\tau^* \in \mathcal{T}_1(A)$ and along any $\kappa$-dimensional ($\kappa \ge 1$) $g_{\mathbf{1}}(\tau^*)$-non-degenerate channel $V \subset \ker\tau^*$, regardless of $\dim A$, the ART phase family~(\ref{eq:paraxial_phase}) over $V \times \R \times \mathcal{O}_V$ at the focal point $\mathbf{c}^* = (\mathbf{0}, k_0, \tau^*)$ can versally unfold only Lagrangian singularities of corank $\kappa \le 2$, restricted in Arnol'd's classification to at most $\{A_2, A_3^\pm, D_4^\pm\}$ ($\mu \le \kappa + 2 \le 4$) when $D \in \{\R, \mathbb{H}\}$ and at most $\{A_2, A_3^\pm, A_4, D_4^\pm, D_5^\pm\}$ ($\mu \le \kappa + 3 \le 5$) when $D = \C$.
\end{corollary}

\begin{proof}
Since $A \cong M_n(D)$ is simple, $\operatorname{rad}(A) = \{0\}$. Because $[M_n(D), M_n(D)]$ consists of matrices whose diagonal sum lies in $[D, D]$ (with $[\R, \R] = [\C, \C] = \{0\}$ and $[\mathbb{H}, \mathbb{H}] = \operatorname{Im}\mathbb{H}$), the trace space $\mathcal{T}(A) = (A/[A, A])^*$ is spanned over $\R$ by $\operatorname{Tr}_\R$ for $D = \R$, by $\operatorname{Tr}_{\mathbb{H}}(X) = \operatorname{Re}\sum_{i=1}^n X_{ii}$ for $D = \mathbb{H}$, and by $\{\operatorname{Re}\operatorname{Tr}_\C, \operatorname{Im}\operatorname{Tr}_\C\}$ for $D = \C$. Thus $d_{\mathcal{T}} \coloneqq \dim_\R \mathcal{T}(A)$ equals $1$ for $D \in \{\R, \mathbb{H}\}$ and $2$ for $D = \C$, for all $n \ge 1$. Because $\mathcal{T}_1(A)$ is a non-empty conic open subset of $\mathcal{T}(A)$ (namely $\R^\times \operatorname{Tr}_D$ for $D \in \{\R, \mathbb{H}\}$ and $\{c\operatorname{Re}\operatorname{Tr}_\C + d\operatorname{Im}\operatorname{Tr}_\C : c \neq 0\}$ for $D = \C$), $\mathbb{P}(\mathcal{T}_1(A))$ is a smooth manifold of dimension $d_{\mathcal{T}} - 1 \in \{0, 1\}$, with $P_{\boldsymbol{\tau}} \equiv 0$ in Proposition~\ref{prop:WM_factorization}.

At the focal point $\mathbf{c}^* = (\mathbf{0}, k_0, \tau^*)$, $S_0(\boldsymbol{\sigma})$ cancels the quadratic part of $k_0 H_\perp(\boldsymbol{\sigma}; \tau^*)$ in (\ref{eq:H_perp})--(\ref{eq:paraxial_phase}), so the central focal germ $f(\boldsymbol{\sigma}) \coloneqq \Phi(\mathbf{c}^*; \boldsymbol{\sigma})$ lies in $\mathfrak{m}_{\boldsymbol{\sigma}}^3$ on $V^* \cong \R^\kappa$ and has corank $\kappa$. Because $\mathbf{1}^*_{\lambda\tau} = \mathbf{1}^*_\tau$ and $F(y; \lambda\tau) = F(y; \tau)$ for all $\lambda \in \R^\times$, both $H_\perp(\boldsymbol{\sigma}; \tau)$ and $\Phi(x_\perp, z, \tau; \boldsymbol{\sigma})$ factor through $V \times \R \times \mathbb{P}(\mathcal{O}_V)$ (where $\mathbb{P}(\mathcal{O}_V) \subset \mathbb{P}(\mathcal{T}_1(A))$ denotes the projective image of $\mathcal{O}_V$), a manifold of dimension $\kappa + 1 + \dim\mathbb{P}(\mathcal{O}_V) \le \kappa + d_{\mathcal{T}}$. By the Mather--Martinet theorem~\cite{Mather1968, Martinet1982}, versal unfolding of $f$ (of finite Milnor number $\mu$) requires the parameter derivatives of $\Phi$ at $\mathbf{c}^*$ to surject onto $\widetilde{Q}_f = \mathfrak{m}_{\boldsymbol{\sigma}} / \langle \nabla f \rangle$, forcing $\mu - 1 = \dim_\R \widetilde{Q}_f \le \kappa + d_{\mathcal{T}}$. On the other hand, since $f \in \mathfrak{m}_{\boldsymbol{\sigma}}^3$ implies $\partial_{\sigma_i} f \in \mathfrak{m}_{\boldsymbol{\sigma}}^2$ ($1 \le i \le \kappa$) and $\mathfrak{m}_{\boldsymbol{\sigma}} \cdot \mathfrak{m}_{\boldsymbol{\sigma}}^2 = \mathfrak{m}_{\boldsymbol{\sigma}}^3$, the image of $\langle \nabla f \rangle$ in $\mathfrak{m}_{\boldsymbol{\sigma}} / \mathfrak{m}_{\boldsymbol{\sigma}}^3$ is spanned over $\R$ by the $\kappa$ cosets $\partial_{\sigma_1} f, \dots, \partial_{\sigma_\kappa} f \bmod \mathfrak{m}_{\boldsymbol{\sigma}}^3$ and has dimension at most $\kappa$. Since $\dim_\R(\mathfrak{m}_{\boldsymbol{\sigma}} / \mathfrak{m}_{\boldsymbol{\sigma}}^3) = \kappa + \frac{1}{2}\kappa(\kappa + 1)$, projecting $\widetilde{Q}_f \twoheadrightarrow \mathfrak{m}_{\boldsymbol{\sigma}} / (\mathfrak{m}_{\boldsymbol{\sigma}}^3 + \langle \nabla f \rangle)$ gives
\begin{equation*}
\frac{\kappa(\kappa + 1)}{2} \le \mu - 1 \le \kappa + d_{\mathcal{T}}. 
\end{equation*}
For every $\kappa \ge 3$, $\frac{1}{2}\kappa(\kappa + 1) - \kappa = \frac{1}{2}\kappa(\kappa - 1) \ge 3 > 2 \ge d_{\mathcal{T}}$, which is impossible; hence $\kappa \le 2$ and $\mu \le \kappa + d_{\mathcal{T}} + 1$. Comparing with Arnol'd's classification of corank-$\le 2$ singularities~\cite{arnold:1985} yields $\{A_2, A_3^\pm\}$ ($\kappa = 1, \mu \le 3$) and $\{D_4^\pm\}$ ($\kappa = 2, \mu \le 4$) when $d_{\mathcal{T}} = 1$, and $\{A_2, A_3^\pm, A_4\}$ ($\kappa = 1, \mu \le 4$) and $\{D_4^\pm, D_5^\pm\}$ ($\kappa = 2, \mu \le 5$) when $d_{\mathcal{T}} = 2$.
\end{proof}

Together, Theorem~\ref{thm:universal_catastrophe} and Corollary~\ref{cor:simple_obstruction} show that within an indecomposable algebra (where semisimplicity reduces to simplicity) versally unfolding catastrophes of corank $\kappa \ge 3$ or Milnor number $\mu > 5$ requires a non-trivial nilpotent radical filtration $\operatorname{rad}(A)$. 

\subsection{Umbilic synthesis in a noncommutative bound quiver algebra}
\label{sec:quiver_example}

Let $Q$ be the quiver~\cite{gabriel:1992} with vertices $\{1, 2\}$, loop $x \colon 1 \to 1$, and transit arrows $a \colon 1 \to 2$, $b \colon 2 \to 1$,
\begin{equation*}
Q: \quad \begin{tikzcd}
1 \arrow[out=150,in=210,loop,swap,"x"] \arrow[r,shift left=1ex,"a"] & 2 \arrow[l,shift left=1ex,"b"]
\end{tikzcd}, 
\end{equation*}
and define $B = \mathbb{R}Q / \mathcal{I}$ (using left-to-right path concatenation), where $\mathcal{I}$ is generated by all paths of length $\ge 4$ together with $aba = 0$ and $bab = 0$. Grading $B = \bigoplus_{m=0}^3 B_m$ by path length gives bases
\begin{align*}
B_0 &= \operatorname{span}\{e_1, e_2\}, & B_1 &= \operatorname{span}\{x, a, b\}, \\
B_2 &= \operatorname{span}\{x^2, ab, ba, xa, bx\}, & B_3 &= \operatorname{span}\{x^3, xab, bxa, abx, x^2a, bx^2\}
\end{align*}
($B_{\ge 4} = \{0\}$), so $\dim B = 16$, $\mathbf{1}_B = e_1 + e_2$, and $\operatorname{rad}(B) = \bigoplus_{m=1}^3 B_m$ has nilpotency index $4$.

The algebra $B$ is indecomposable (any central idempotent projects into \\ $\{0, e_1, e_2, \mathbf{1}_B\} \subset B/\operatorname{rad}(B) \cong \R e_1 \oplus \R e_2$, and $(e_1 + r)a \equiv a \not\equiv 0 \equiv a(e_1 + r) \pmod{B_{\ge 2}}$ for all $r \in \operatorname{rad}(B)$), non-local, and noncommutative ($ab \neq ba$). Moreover, $B$ is not symmetric Frobenius: $x^2a = [e_1, x^2a] \in [B, B]$ satisfies $x^2a B = \R x^2a \subset [B, B]$, forcing $x^2a \in \ker\hat{\tau}_B$ for every $\tau_B \in \mathcal{T}(B)$. Passing to its $32$-dimensional trivial extension $A = B \ltimes B^*$, every trace $\tau_B \in \mathcal{T}(B) = (B/[B, B])^*$ with $\tau_B(\mathbf{1}_B) \neq 0$ lifts via (\ref{eq:triv_ext_trace})--(\ref{eq:triv_ext_pairing}) to $\tau \in \mathcal{T}_1(A)$ (since $(b_1\cdot\phi_2 - \phi_2\cdot b_1)(\mathbf{1}_B) = 0$ for all $b_1 \in B$, $\phi_2 \in B^*$).

Embed the two-dimensional channel $V = \operatorname{span}\{v_1, v_2\} \subset \operatorname{rad}(B) \times \{0\} \subset \operatorname{rad}(A)$ via $v_1 = (x, 0)$ and $v_2 = (a + b, 0)$, identifying $V$ with its base projection in $B$ since products in $B \times \{0\}$ have zero $B^*$ component. Direct commutator inspection yields $[B, B] \cap B_1 = \operatorname{span}\{a, b\}$ (from $a = [e_1, a]$, $b = [b, e_1]$), $[B, B] \cap B_2 = \operatorname{span}\{xa, bx, ab - ba\}$, and $[B, B] \cap B_3 = \operatorname{span}\{x^2a, bx^2, xab - abx, abx - bxa\}$. Hence $\tau|_V = 0$ reduces to $\tau_B(x) = 0$. On $B_2$, $\tau(v_1 v_2 + v_2 v_1) = \tau_B(xa + bx) = 0$ vanishes identically while $\tau_B(x^2)$ and $2\tau_B(ab) = \tau_B(ab + ba)$ form independent coordinates on $(B_2 / ([B, B] \cap B_2))^*$; setting $\tau_B^*(x^2) = 2\tau_B^*(ab) = -\tau^*(\mathbf{1}_A)$ calibrates $g_{\mathbf{1},\perp}(\tau^*) = I_2$ ($Q_{\tau^*}(\boldsymbol{\sigma}) = \sigma_1^2 + \sigma_2^2$) by (\ref{eq:g_one}).

On $B_3$, expanding $w^3$ for $w = \sigma_1 v_1 + \sigma_2 v_2$ via $xb = ax = aba = bab = 0$ gives $w^3 = \sigma_1^3 x^3 + \sigma_1^2\sigma_2(x^2a + bx^2) + \sigma_1\sigma_2^2(xab + bxa + abx)$. Since $\tau$ vanishes on $x^2a, bx^2 \in [B, B]$ and satisfies $\tau(xab) = \tau(bxa) = \tau(abx)$, we obtain
\begin{equation}\label{eq:quiver_cubic_trace}
\frac{1}{3\tau(\mathbf{1}_A)}\tau(w^3) = \alpha \sigma_1^3 + \beta \sigma_1 \sigma_2^2, \qquad \alpha \coloneqq \frac{\tau_B(x^3)}{3\tau(\mathbf{1}_A)}, \quad \beta \coloneqq \frac{\tau_B(xab)}{\tau(\mathbf{1}_A)}, 
\end{equation}
where $\alpha, \beta \in \R$ are independently tunable since $\{[x^3], [xab]\}$ is a basis of $B_3 / ([B, B] \cap B_3)$.

Writing $y = y_0\mathbf{1}_A + r + \tilde{w}$ across (\ref{eq:orth_split}) ($r \in V$, $\tilde{w} \in W(\tau)$), $g_{\mathbf{1}}(\tau)$-orthogonality restricts $(r, \tilde{w})$ couplings in $\mathcal{L}(y; \tau) = \frac{1}{2}F(y; \tau)^2$ to $O(|r|^2|\tilde{w}|)$, so $\mathbf{p}|_{W(\tau)} = \partial_{\tilde{w}}\mathcal{L} = 0$ forces $\tilde{w} = O(|r|^2)$, affecting $\mathcal{L}$ only at $O(|r|^4)$. On $\R\mathbf{1}_A \oplus V$ ($\tau(r) = 0$), expanding $\mathcal{L}(y_0\mathbf{1}_A + r; \tau) = \frac{1}{2}y_0^2 + \frac{1}{2}g_{\mathbf{1},\perp}(r, r; \tau) + \frac{1}{3\tau(\mathbf{1}_A)}\tau(r^3) + O(|r|^4)$ for $y_0 = 1 + O(|r|^2)$ and Legendre-transforming to $\mathbf{p} = (1 - h)\mathbf{1}^*_\tau + \boldsymbol{\sigma} \in W(\tau)^0$ ($h = O(|\boldsymbol{\sigma}|^2)$) flips the cubic sign in $F^*\big((1 - h)\mathbf{1}^*_\tau + \boldsymbol{\sigma}; \tau\big)^2 = 1 - 2h + g_{\mathbf{1},\perp}(\tau)^{-1}(\boldsymbol{\sigma}, \boldsymbol{\sigma}) - \frac{2}{3\tau(\mathbf{1}_A)}\tau(w_{\boldsymbol{\sigma}}^3) + O(|\boldsymbol{\sigma}|^4) = 1$ with $w_{\boldsymbol{\sigma}} \coloneqq g_{\mathbf{1},\perp}(\tau)^{-1}\boldsymbol{\sigma} \in V$, yielding
\begin{equation}\label{eq:H_perp_cubic}
H_\perp(\boldsymbol{\sigma}; \tau) = \frac{1}{2}g_{\mathbf{1},\perp}(\tau)^{-1}(\boldsymbol{\sigma}, \boldsymbol{\sigma}) - \frac{1}{3\tau(\mathbf{1}_A)}\tau(w_{\boldsymbol{\sigma}}^3) + O(|\boldsymbol{\sigma}|^4). 
\end{equation}

At $\tau^*$ ($g_{\mathbf{1},\perp}(\tau^*) = I_2$, $w_{\boldsymbol{\sigma}} = \sigma_1 v_1 + \sigma_2 v_2$), (\ref{eq:quiver_cubic_trace})--(\ref{eq:H_perp_cubic}) give $H_\perp(\boldsymbol{\sigma}; \tau^*) = \frac{1}{2}(\sigma_1^2 + \sigma_2^2) - (\alpha\sigma_1^3 + \beta\sigma_1\sigma_2^2) + O(|\boldsymbol{\sigma}|^4)$. Under $S_0(\boldsymbol{\sigma}) = \frac{1}{2}k_0(\sigma_1^2 + \sigma_2^2)$ ($k_0 \neq 0$), the focal germ of (\ref{eq:paraxial_phase}) at $\mathbf{c}^* = (\mathbf{0}, k_0, \tau^*)$ is
\begin{equation*}
f_0(\boldsymbol{\sigma}) = k_0\mathcal{C}(\boldsymbol{\sigma}) + O(|\boldsymbol{\sigma}|^4), \qquad \mathcal{C}(\boldsymbol{\sigma}) = \sigma_1(\alpha\sigma_1^2 + \beta\sigma_2^2).
\end{equation*}
Whenever $\alpha\beta \neq 0$, $\langle \nabla\mathcal{C} \rangle = \langle 3\alpha\sigma_1^2 + \beta\sigma_2^2, \, 2\beta\sigma_1\sigma_2 \rangle$ contains $\mathfrak{m}_{\boldsymbol{\sigma}}^3$, so $f_0$ is $3$-determined~\cite{Mather1968} with Milnor number $\mu = \dim(C^\infty_{\mathbf{0}}(V^*)/\langle \nabla\mathcal{C} \rangle) = 4$, realizing the elliptic umbilic $D_4^-$ when $\beta/\alpha = 3\tau_B^*(xab)/\tau_B^*(x^3) < 0$ and the hyperbolic umbilic $D_4^+$ when $\beta/\alpha > 0$. Finally, varying $x_\perp \in V$ spans $\operatorname{span}\{\sigma_1, \sigma_2\}$ while varying the independent quadratic trace ratios $-\tau_B(x^2)/\tau(\mathbf{1}_A)$ and $-2\tau_B(ab)/\tau(\mathbf{1}_A)$ in $\mathcal{O}_V$ independently varies $\sigma_1^2$ and $\sigma_2^2$ modulo $\mathfrak{m}_{\boldsymbol{\sigma}}^3 \subset \langle \nabla f_0 \rangle$, spanning the $3$-dimensional reduced Milnor algebra $\widetilde{Q}_{f_0}$ and versally unfolding $D_4^\pm$ (with both classes persisting across open neighborhoods of unmarked non-degenerate $2$-planes in $\operatorname{rad}(A)$ by stability of $3$-determinacy).

\subsection{An example of equidimensional algebras with identical homological dimensions that are distinguished via ART catastrophes}
\label{sec:homological_separation}

Whereas extracting caustic invariants in generic pseudo-Finsler geometry requires integrating nonlinear geodesic equations, translation invariance and Proposition~\ref{prop:WM_factorization} reduce local ART singularity classification to finite-dimensional multilinear algebra. To see how ART catastrophe classes detect real algebraic twists invisible to linear homological dimensions ($\dim HH^k$, $\dim HC_k$) and quadratic metric signatures, consider the six-dimensional commutative $\R$-algebras
\begin{equation*}
A_{\mathrm{hyp}} \coloneqq \frac{\mathbb{R}[x, y]}{\big(xy, \; x^3 - y^3\big)},  \qquad
A_{\mathrm{ell}} \coloneqq \frac{\mathbb{R}[x, y]}{\big(x^2 + y^2, \; y^3\big)}. 
\end{equation*}
In $A_{\mathrm{hyp}}$, $xy = 0$ and $x^3 = y^3$ give $x^2y = xy^2 = x^4 = 0$ and graded basis $\{\mathbf{1}\} \cup \{x, y\} \cup \{x^2, y^2\} \cup \{x^3\}$; in $A_{\mathrm{ell}}$, $x^2 + y^2 = y^3 = 0$ gives $x^2y = x^4 = 0$, $xy^2 = -x^3$, and graded basis $\{\mathbf{1}\} \cup \{x, y\} \cup \{x^2, xy\} \cup \{x^3\}$. As local complete intersections, both are $6$-dimensional symmetric Frobenius algebras with Hilbert series $H(t) = 1 + 2t + 2t^2 + t^3$, radical nilpotency index $N = 4$, and socle $\R x^3$.

Over $\C$, the substitution $x \mapsto u + v$, $y \mapsto i(u - v)$ sends $x^2 + y^2 \mapsto 4uv$ and $y^3 \mapsto -i(u^3 - v^3) \pmod{uv}$, exhibiting $A_{\mathrm{hyp}}$ and $A_{\mathrm{ell}}$ as real forms of the same complex algebra, \begin{equation*}
A_{\mathrm{ell}} \otimes_{\mathbb{R}} \mathbb{C} \xrightarrow{\;\sim\;} \frac{\mathbb{C}[u, v]}{\big(uv, \; u^3 - v^3\big)} \cong A_{\mathrm{hyp}} \otimes_{\mathbb{R}} \mathbb{C}. 
\end{equation*}
Since Hochschild cohomology and cyclic homology commute with flat field extensions, $\dim HH^k(A_{\mathrm{hyp}}, A_{\mathrm{hyp}}) = \dim HH^k(A_{\mathrm{ell}}, A_{\mathrm{ell}})$ and $\dim HC_k(A_{\mathrm{hyp}}) = \dim HC_k(A_{\mathrm{ell}})$ for all $k \ge 0$. Quadratic metric signatures likewise fail to separate them: for every $\tau \in \mathcal{T}_1$ (characterized on both algebras by $\tau(\mathbf{1}) \neq 0$ and $s \coloneqq \tau(x^3) \neq 0$), $\operatorname{rad}(A)^2 \cap \ker\tau$ is totally $g_{\mathbf{1}}(\tau)$-isotropic and pairs non-degenerately with $\operatorname{rad}(A)/\operatorname{rad}(A)^2$, forcing $g_{\mathbf{1}}(\tau)$ to have signature $(2, 2)$ on $\operatorname{rad}(A) \cap \ker\tau$ and---with the positive-definite $(\mathbf{1}, x^3)$-block---signature $(4, 2)$ on $A$ (while $\hat{\tau}$ has signature $(3, 3)$).

Although in this small commutative example one can separate the two real forms directly via the $GL_2(\R)$-orbit of the cubing form $\operatorname{Sym}^3(\operatorname{rad}(A)/\operatorname{rad}(A)^2) \to \operatorname{soc}(A)$ (or cochain-level cup products and mixed Gerstenhaber brackets on $HH^\bullet(A, A)$, even though the real Lie algebras $HH^1(A, A)$ remain isomorphic), the ART framework converts this distinction directly into disjoint symplectic caustic classes on $T^*A$ via trace evaluations alone. Let $\tau \in \mathcal{T}_1$ (with $s = \tau(x^3) \neq 0$) and let $V \subset \operatorname{rad}(A) \cap \ker\tau$ be any two-dimensional $g_{\mathbf{1}}(\tau)$-non-degenerate channel. In particular, positive-definite Euclidean channels ($g_{\mathbf{1},\perp}(\tau) = I_2$) exist on both algebras: at $\tau(\mathbf{1}) = \tau(x^3) = 1$ with $\tau|_{A_1 \oplus A_2} = 0$, one checks directly via (\ref{eq:g_one}) that $v_1 = \frac{1}{\sqrt{2}}(x - x^2), \, v_2 = \frac{1}{\sqrt{2}}(y - y^2)$ in $A_{\mathrm{hyp}}$ and $v_1 = \frac{1}{\sqrt{2}}(x - x^2), \, v_2 = \frac{1}{\sqrt{2}}(y + xy)$ in $A_{\mathrm{ell}}$ are $g_{\mathbf{1}}(\tau)$-orthonormal.

By (\ref{eq:H_perp_cubic}), $H_\perp(\boldsymbol{\sigma}; \tau) = \frac{1}{2}g_{\mathbf{1},\perp}(\tau)^{-1}(\boldsymbol{\sigma}, \boldsymbol{\sigma}) - \frac{1}{3\tau(\mathbf{1})}\tau(w_{\boldsymbol{\sigma}}^3) + O(|\boldsymbol{\sigma}|^4)$ with $w_{\boldsymbol{\sigma}} = g_{\mathbf{1},\perp}(\tau)^{-1}\boldsymbol{\sigma} \in V$. Decomposing $w_{\boldsymbol{\sigma}} = r_1 + r_2 + r_3$ ($r_m \in A_m$), $\operatorname{rad}(A)^4 = \{0\}$ implies $w_{\boldsymbol{\sigma}}^3 = r_1^3$. 
Writing $r_1(\boldsymbol{\sigma}) = a(\boldsymbol{\sigma})x + b(\boldsymbol{\sigma})y$ on $V^* \cong \R^2$, the cubic trace form $\frac{1}{3\tau(\mathbf{1})}\tau(w_{\boldsymbol{\sigma}}^3)$ governing the focal $3$-jet $j^3_{\mathbf{0}}\Phi(\mathbf{0}, k_0, \tau; \boldsymbol{\sigma}) = \frac{k_0}{3\tau(\mathbf{1})}\tau(w_{\boldsymbol{\sigma}}^3)$ under $S_0(\boldsymbol{\sigma}) = \frac{1}{2}k_0 Q_\tau(\boldsymbol{\sigma})$ ($k_0 \neq 0$, $z^* = k_0$) evaluates across all $\tau \in \mathcal{T}_1$ and non-degenerate channels $V \subset \operatorname{rad}(A) \cap \ker\tau$ as follows.
 \begin{itemize}
    \item In $A_{\mathrm{hyp}}$ ($xy = 0$, $y^3 = x^3$), $\frac{1}{3\tau(\mathbf{1})}\tau(w_{\boldsymbol{\sigma}}^3) = \frac{s}{3\tau(\mathbf{1})}\big(a(\boldsymbol{\sigma})^3 + b(\boldsymbol{\sigma})^3\big)$. Whenever $w \mapsto r_1$ has rank $2$ on $V$, $(a, b)$ is a linear isomorphism on $\R^2$, so this binary cubic has one real root line ($a = -b$) and a complex-conjugate pair. Thus $A_{\mathrm{hyp}}$ exclusively realizes the \emph{hyperbolic umbilic} $D_4^+$---giving $j^3_{\mathbf{0}} H_\perp = \frac{1}{2}(\sigma_1^2 + \sigma_2^2) - \frac{1}{6\sqrt{2}}(\sigma_1^3 + \sigma_2^3)$ on the orthonormal channel above---and never $D_4^-$ (while rank $\le 1$ projections yield non-isolated singularities of infinite Milnor number).
    \item In $A_{\mathrm{ell}}$ ($x^2y = y^3 = 0$, $xy^2 = -x^3$), $\frac{1}{3\tau(\mathbf{1})}\tau(w_{\boldsymbol{\sigma}}^3) = \frac{s}{3\tau(\mathbf{1})}\big(a(\boldsymbol{\sigma})^3 - 3a(\boldsymbol{\sigma})b(\boldsymbol{\sigma})^2\big) = \frac{s}{3\tau(\mathbf{1})}a(\boldsymbol{\sigma})\big(a(\boldsymbol{\sigma}) - \sqrt{3}\,b(\boldsymbol{\sigma})\big)\big(a(\boldsymbol{\sigma}) + \sqrt{3}\,b(\boldsymbol{\sigma})\big)$. Whenever $(a, b)$ has rank $2$, this cubic factors into three distinct real root lines, so $A_{\mathrm{ell}}$ exclusively realizes the \emph{elliptic umbilic} $D_4^-$---giving $j^3_{\mathbf{0}} H_\perp = \frac{1}{2}(\sigma_1^2 + \sigma_2^2) - \frac{1}{6\sqrt{2}}(\sigma_1^3 - 3\sigma_1\sigma_2^2)$ on the orthonormal channel above---and never $D_4^+$.
\end{itemize}
By Corollary~\ref{cor:iso_invariance}, this disjointness of realizable isolated umbilic classes ($D_4^+$ vs.\ $D_4^-$) confirms $A_{\mathrm{hyp}} \not\cong_{\R} A_{\mathrm{ell}}$ (Table~\ref{tab:homological_separation}).

\begin{table}[htbp]
\centering
\footnotesize
\renewcommand{\arraystretch}{1.2}
\setlength{\tabcolsep}{5pt}
\caption{Comparison of classical, homological, quadratic-metric, and ART catastrophe invariants for $A_{\mathrm{hyp}} = \R[x,y]/(xy,\,x^3 - y^3)$ and $A_{\mathrm{ell}} = \R[x,y]/(x^2 + y^2,\,y^3)$.}
\label{tab:homological_separation}
\begin{tabularx}{\textwidth}{>{\raggedright\arraybackslash}X c c c}
\toprule
\textbf{Algebraic or Geometric Invariant} & $A_{\mathrm{hyp}}$ & $A_{\mathrm{ell}}$ & \textbf{Separates?} \\
\midrule
Dimension $\dim A$ & $6$ & $6$ & No \\
Radical nilpotency index $N$ & $4$ & $4$ & No \\
Commutative symmetric Frobenius (complete intersection) & Yes & Yes & No \\
Graded Hilbert series $H_A(t)$ & $1 + 2t + 2t^2 + t^3$ & $1 + 2t + 2t^2 + t^3$ & No \\
Complexification $A \otimes_{\R} \C$ & $\dfrac{\C[u,v]}{(uv,\,u^3-v^3)}$ & $\dfrac{\C[u,v]}{(uv,\,u^3-v^3)}$ & No \\
Hochschild cohomology dimensions $\dim HH^k(A, A)$ ($\forall k \ge 0$) & Identical & Identical & No \\
Cyclic homology dimensions $\dim HC_k(A)$ ($\forall k \ge 0$) & Identical & Identical & No \\
Global signature of $g_{\mathbf{1}}(\tau)$ on $A$ ($\forall \tau \in \mathcal{T}_1$) & $(4, 2)$ & $(4, 2)$ & No \\
Signature of $g_{\mathbf{1}}(\tau)$ on $\operatorname{rad}(A) \cap \ker\tau$ ($\forall \tau \in \mathcal{T}_1$) & $(2, 2)$ & $(2, 2)$ & No \\
\midrule
\textbf{Isolated ART radical catastrophe ($\kappa = 2$, $z^* = k_0 \neq 0$)} & \textbf{$D_4^+$ (Hyperbolic)} & \textbf{$D_4^-$ (Elliptic)} & \textbf{Yes} \\
\bottomrule
\end{tabularx}
\end{table}

\subsection{Separating derived-equivalent algebras via ART cubic invariants}
\label{sec:derived_separation}

By theorems of Rickard and Keller~\cite{Rickard1989,Rickard1991,Keller1994}, any two finite-dimensional algebras $A_1$ and $A_2$ with equivalent bounded derived module categories $D^b(\operatorname{mod}\text{-}A_1) \simeq D^b(\operatorname{mod}\text{-}A_2)$ share canonically isomorphic centers $Z(A_1) \cong Z(A_2)$, isomorphic full Gerstenhaber algebras $\big(HH^\bullet(A_1, A_1), \smile, [\cdot, \cdot]_{\mathrm{G}}\big) \cong \big(HH^\bullet(A_2, A_2), \smile, [\cdot, \cdot]_{\mathrm{G}}\big)$, and isomorphic cyclic homologies $HC_\bullet(A_1) \cong HC_\bullet(A_2)$ in all degrees. However, because a derived equivalence is mediated by a chain complex of bimodules rather than an algebra homomorphism, it need not preserve the nonlinear radical power maps governing Proposition~\ref{prop:WM_factorization}.

For any smooth germ $f \in \mathfrak{m}_{\boldsymbol{\sigma}}^3$ on $\R^\kappa$, the \emph{essential cubic rank} $\operatorname{rank}_{\mathrm{ess}}(j^3_{\mathbf{0}} f)$ of its cubic $3$-jet $C = j^3_{\mathbf{0}} f$, which is defined as the codimension of $\ker_{\mathrm{lin}}(C) \coloneqq \{ \mathbf{v} \in \R^\kappa : D_{\mathbf{v}} C \equiv 0 \}$, or equivalently the minimal number of linear coordinates needed to express $C$, is invariant under all local diffeomorphisms fixing $\mathbf{0}$ and is thus an invariant of the smooth right-equivalence class of $f$.

\begin{theorem}\label{thm:derived_separation}
There exist non-isomorphic, equidimensional symmetric Frobenius algebras over $\R$ that are derived equivalent, and hence indistinguishable by cyclic (co)homology or any Gerstenhaber-algebraic operation on Hochschild cohomology, but whose ART focal phase germs across non-degenerate radical channels have distinct maximal essential cubic ranks.
\end{theorem}

\begin{proof}
Consider two Brauer tree algebras over $\R$ with $e = 4$ edges and exceptional multiplicity $m = 3$~\cite{Rickard1989} (using left-to-right path concatenation), \begin{enumerate}
    \item Let $A_{\mathrm{star}}$ correspond to the star tree of degree sequence $(4, 1, 1, 1, 1)$ with multiplicity $m = 3$ at the leaf of edge $1$: vertices $\{1, 2, 3, 4\}$, a leaf loop $x \colon 1 \to 1$, and a $4$-cycle $1 \xrightarrow{a_1} 2 \xrightarrow{a_2} 3 \xrightarrow{a_3} 4 \xrightarrow{a_4} 1$, modulo $x a_1 = a_4 x = 0$, $x^3 = \omega_1$, $x^4 = 0$, and $\omega_i a_i = 0$ where $\omega_i \coloneqq a_i a_{i+1} a_{i+2} a_{i+3}$ ($i = 1, 2, 3, 4 \bmod 4$).
    \item Let $A_{\mathrm{line}}$ correspond to the line tree of degree sequence $(2, 2, 2, 1, 1)$ with multiplicity $m = 3$ at the internal vertex joining edges $1$ and $2$: vertices $\{1, 2, 3, 4\}$ and three $2$-cycles $1 \underset{b}{\overset{a}{\rightleftarrows}} 2 \underset{d}{\overset{c}{\rightleftarrows}} 3 \underset{g}{\overset{f}{\rightleftarrows}} 4$, modulo $ac = db = cf = gd = 0$, $(ab)^3 \eqqcolon \omega_1$, $(ba)^3 = cd \eqqcolon \omega_2$, $dc = fg \eqqcolon \omega_3$, $gf \eqqcolon \omega_4$, and $(ab)^3 a = (ba)^3 b = 0$.
\end{enumerate}
Both $A_{\mathrm{star}}$ and $A_{\mathrm{line}}$ are symmetric Frobenius algebras of dimension $\sum_v m_v(\deg v)^2 = 22$ sharing the Brauer tree invariants $(e, m) = (4, 3)$, so $D^b(\operatorname{mod}\text{-}A_{\mathrm{star}}) \simeq D^b(\operatorname{mod}\text{-}A_{\mathrm{line}})$ by Rickard's theorem~\cite{Rickard1989}.

For any $\kappa$-dimensional $g_{\mathbf{1}}(\tau^*)$-non-degenerate channel $V \subset \operatorname{rad}(A) \cap \ker\tau^*$ at the focal point $z^* = k_0 \neq 0$ under $S_0(\boldsymbol{\sigma}) = \frac{1}{2}k_0 Q_{\tau^*}(\boldsymbol{\sigma})$, (\ref{eq:paraxial_phase}) and (\ref{eq:H_perp_cubic}) give $j^3_{\mathbf{0}} f_0(\boldsymbol{\sigma}) = \frac{k_0}{3\tau^*(\mathbf{1})}\tau^*(w_{\boldsymbol{\sigma}}^3)$ with $w_{\boldsymbol{\sigma}} = g_{\mathbf{1},\perp}(\tau^*)^{-1}\boldsymbol{\sigma} \in V$. Since $g_{\mathbf{1},\perp}(\tau^*)^{-1} \colon V^* \to V$ is a linear isomorphism and $k_0/\tau^*(\mathbf{1}) \neq 0$, $\operatorname{rank}_{\mathrm{ess}}(j^3_{\mathbf{0}} f_0)$ equals the essential cubic rank of $w \mapsto \tau^*(w^3)$ on $V$.

In $A_{\mathrm{line}}$, the $18$-dimensional radical $\operatorname{rad}(A_{\mathrm{line}})$ splits as $U_{\mathrm{act}} \oplus U_{\mathrm{null}}$, where $U_{\mathrm{act}} \coloneqq \operatorname{span}\{a, b, ab, ba, aba, bab, (ab)^2, (ba)^2\}$ ($\dim U_{\mathrm{act}} = 8$) and \[ U_{\mathrm{null}} \coloneqq \operatorname{span}\{\omega_1, \omega_2, \omega_3, \omega_4, \, (ab)^2a, (ba)^2b, \, c, d, f, g\}\] ($\dim U_{\mathrm{null}} = 10$). Every triple product in $\operatorname{rad}(A_{\mathrm{line}})$ containing a factor from $U_{\mathrm{null}}$ vanishes identically, \begin{itemize}
    \item $\omega_1, \dots, \omega_4 \in \operatorname{soc}(A_{\mathrm{line}})$ annihilate $\operatorname{rad}(A_{\mathrm{line}})$ on both sides;
    \item $(ab)^2a$ and $(ba)^2b$ annihilate $\{c, d, f, g\}$ via $ac = db = 0$, while multiplying either by two $\{a, b\}$-paths yields a path of length $\ge 7$, which vanishes since $(ab)^3a = (ba)^3b = 0$;
    \item $\{c, d\}$ and $\{f, g\}$ annihilate $\{a, b\}$ and each other across vertices $2$ and $3$ via $ac = db = cf = gd = 0$, while all alternating length-$3$ paths in $\{c, d\}$ or $\{f, g\}$ vanish via the socle identifications:
    \begin{equation*}
    cdc = (ba)^2b(ac) = 0, \quad dcd = f(gd) = 0, \quad fgf = d(cf) = 0, \quad gfg = (gd)c = 0. 
    \end{equation*}
\end{itemize}
Writing any $w \in \operatorname{rad}(A_{\mathrm{line}})$ as $w = w_{\mathrm{act}} + w_{\mathrm{null}}$ therefore gives $w^3 = w_{\mathrm{act}}^3$, so $U_{\mathrm{null}} \cap V \subset \ker_{\mathrm{lin}}\big(\tau^*(w^3)\big)$ and hence $\operatorname{rank}_{\mathrm{ess}}(j^3_{\mathbf{0}} f_0) \le \dim U_{\mathrm{act}} = 8$ for every admissible pair $(\tau^*, V)$ in $A_{\mathrm{line}}$.

Conversely, in $A_{\mathrm{star}}$, define $\tau^* \in \mathcal{T}_1(A_{\mathrm{star}})$ by $\tau^*(e_i) = \frac{1}{4}$, $\tau^*(\omega_i) = 1$, $\tau^*(x^2) = -1$, and $\tau^*(x) = 0$, vanishing on the commutator paths $a_i$, $m_i \coloneqq a_i a_{i+1}$, and $b_i \coloneqq a_{i+1} a_{i+2} a_{i+3}$ ($i = 1, 2, 3, 4 \bmod 4$); direct inspection confirms that $\hat{\tau}^*$ is block-non-degenerate across $\{e_1, \omega_1, x, x^2\}$, $\{e_j, \omega_j\}_{j=2}^4$, and the paired cycle paths. On the $13$-dimensional subspace
\begin{equation*}
V \coloneqq \operatorname{span}\big\{ x, \; a_1, \, a_2, \, a_3, \, a_4, \; m_1, \, m_2, \, m_3, \, m_4, \; b_1, \, b_2, \, b_3, \, b_4 \big\} \subset \operatorname{rad}(A_{\mathrm{star}}) \cap \ker\tau^*, 
\end{equation*}
$g_{\mathbf{1}}(u, v; \tau^*) = -\tau^*(uv)$ decomposes into orthogonal invertible blocks with $g_{\mathbf{1}}(x, x; \tau^*) = 1$, $g_{\mathbf{1}}(m_1, m_3; \tau^*) = g_{\mathbf{1}}(m_2, m_4; \tau^*) = -1$, and $g_{\mathbf{1}}(a_i, b_j; \tau^*) = -\delta_{ij}$ (all other basis pairings vanish), giving $\det g_{\mathbf{1},\perp}(\tau^*) = 1 \neq 0$. Writing $w = \xi x + \sum_{i=1}^4 (\alpha_i a_i + \mu_i m_i + \beta_i b_i) \in V$, expanding $w^3$ via $\tau^*(x^3) = \tau^*(\omega_i) = 1$ gives
\begin{equation*}
\tau^*(w^3) = \xi^3 + 3\big(\alpha_1 \alpha_2 \mu_3 + \alpha_2 \alpha_3 \mu_4 + \alpha_3 \alpha_4 \mu_1 + \alpha_4 \alpha_1 \mu_2\big). 
\end{equation*}
Vanishing of \[ D_{\delta w}\big(\tau^*(w^3)\big) = 3\delta\xi\,\xi^2 + 3\sum_{i=1}^4 \big[\delta\mu_{i+2}\,\alpha_i\alpha_{i+1} + \mu_{i+2}(\delta\alpha_i\,\alpha_{i+1} + \alpha_i\,\delta\alpha_{i+1})\big] \equiv 0,\] forces $\delta\xi = \delta\mu_{i+2} = \delta\alpha_i = 0$ for all $i = 1, \dots, 4$. Hence $\ker_{\mathrm{lin}}\big(\tau^*(w^3)\big) = \operatorname{span}\{b_1, b_2, b_3, b_4\}$ has dimension $4$, yielding $\operatorname{rank}_{\mathrm{ess}}(j^3_{\mathbf{0}} f_0) = 13 - 4 = 9 > 8$ and separating $A_{\mathrm{star}}$ from $A_{\mathrm{line}}$ by Corollary~\ref{cor:iso_invariance}.
\end{proof}

\section{Navigating algebraic deformation space with ART-inspired metrics}
Let $\mathcal{M}$ be the ambient space of structure constants of real, not-necessarily associative $n$-dimensional algebras. There are uncountably many ways to endow the manifold $\mathcal{M}$ with a metric structure, which in general are not canonical with respect to the moduli space of associative algebras. However, in Section~\ref{sec:def_metric} we show that, starting from a point $\mu_0 \in \mathcal{M}$ representing a non-rigid symmetric Frobenius algebra, a constrained variational principle endows trace-fibered bundles over trace-extendable unobstructed deformation sheets with canonical $GL_{\mathbf{1}}(V)$-equivariant pseudo-Riemannian structures, simultaneously launching families of Maurer--Cartan paths as geodesic sprays alongside their co-evolving trace calibrations. Whereas classical formal Maurer--Cartan integration must make non-canonical coboundary ($1$-cochain gauge) choices at every power-series order $t^k$ while co-evolving cyclic traces through Getzler's Gauss--Manin connection~\cite{getzler:1993} on the infinite-dimensional $(b, B)$ periodic cyclic bar bicomplex, our bundle metric resolves coboundary gauge ambiguity at $t = 0$ by orthogonal projection onto a canonical harmonic cocycle slice and fixes all higher-order jets uniquely via an autonomous finite-dimensional geodesic ODE. Though such sprays cannot originate at a non-rigid associative algebra that fails to be symmetric Frobenius (which otherwise lies on the metric degeneracy locus), one can at least use the algebra's symmetric Frobenius trivial extension for that purpose.

Section~\ref{sec:2d_deformations} illustrates this program on the two-dimensional commutative moduli space.  Section~\ref{sec:nakayama_geodesics} demonstrates the full noncommutative apparatus, with active Gerstenhaber multipliers, co-rotating traces, and complete avoidance of the non-terminating periodic cyclic bar bicomplex, on the six-dimensional symmetric Frobenius Nakayama algebra $N_6$ of infinite global dimension.   Section~\ref{sec:A3_spray} exhibits a two-dimensional geodesic spray launched from the three-dimensional algebra $A_3 = \mathbb{R}[x]/(x^3)$.

\subsection{The Pseudo-Riemannian metric on the associative locus}
\label{sec:def_metric}

\subsubsection{Construction of the action integral}
Let $V$ be an $n$-dimensional real vector space with a distinguished vector $\mathbf{1} \in V$, and let $\mathcal{M} = V \otimes V^* \otimes V^*$ be the space of bilinear multiplications (with structure constants $c_{ih}^k$), where $\mu \in \mathcal{M}$ defines a not-necessarily associative algebra $B_\mu$. For any covector $\tau_\mu \in V^*$ with $\tau_\mu(\mathbf{1}) \neq 0$, generalizing (\ref{eq:g_one}) yields the canonical symmetrized bilinear~form
\begin{equation}
\label{eq:g_sym_mu}
g_{\mathbf{1}}^{B_\mu}\big(u, v; \tau_\mu\big) \coloneq \frac{2}{\tau_\mu(\mathbf{1})^2}\tau_\mu(u)\tau_\mu(v) - \frac{1}{2\tau_\mu(\mathbf{1})}\left(\tau_\mu(uv)+\tau_\mu(vu)\right),
\end{equation}
which is homogeneous of degree $0$ in $\tau_\mu$ and thus depends only on the projective class $\bar{\tau}_\mu = \tau_\mu/\tau_\mu(\mathbf{1})$.

Because no field of linear functionals $\tau_\mu$ is canonically distinguished across $\mathcal{M}$, we work on the extended configuration bundle $\pi \colon \mathcal{E} \to \mathcal{M}$ whose fiber over $\mu \in \mathcal{M}$ is the open domain of \textit{admissible functionals},
\begin{equation}\label{eq:E_mu_fiber}
\mathcal{E}_\mu = \left\{ \tau \in V^* \;\big|\; \tau(\mathbf{1}) \neq 0 \text{ and } \Delta_\mu(\tau) \coloneqq \det g_{\mathbf{1}}^{B_\mu}(\tau) \neq 0 \right\}, 
\end{equation}
with $\Delta_\mu(\tau)$ evaluated in a fixed basis of $V$. The degeneracy hypersurface $\Delta_\mu(\tau) = 0$ is basis-independent and will be called the \emph{trace horizon}. At any $(\mu, \tau_\mu) \in \mathcal{E}$, write $g_\mu \coloneqq g_{\mathbf{1}}^{B_\mu}(\tau_\mu)$ on $V$ and $g_\mu^{-1}$ for its inverse bilinear form on $V^*$, represented in dual bases by $[g_{ij}(\mu; \tau_\mu)]$ and $[g^{ij}(\mu; \tau_\mu)]$, respectively. 

Mirroring the $(1,2)$-tensor structure of $c_{ih}^k$ lifts $g_\mu$ to a horizontal pseudo-Riemannian metric $G_\mu(\tau_\mu)$ on $\mathcal{M}$, given on rank-one tangent variations by
\begin{equation}\label{eq:G_horizontal}
G_\mu\big(u \otimes \alpha_1 \otimes \beta_1, \, v \otimes \alpha_2 \otimes \beta_2; \tau_\mu\big) \coloneqq g_\mu(u, v) \, g_\mu^{-1}(\alpha_1, \alpha_2) \, g_\mu^{-1}(\beta_1, \beta_2).
\end{equation}
Normalizing $g^{ij}$ by $\tau_\mu(\mathbf{1})^{-2}$ defines the fiber-scale-invariant vertical metric
\begin{equation}\label{eq:h_vertical}
h^{ij}(\mu; \tau_\mu) \coloneqq \frac{1}{\tau_\mu(\mathbf{1})^2} g^{ij}(\mu; \tau_\mu). 
\end{equation}
Because $\mathcal{E}$ is an open submanifold of $\mathcal{M} \times V^*$, its tangent space splits canonically as $T_{(\mu, \tau_\mu)}\mathcal{E} = T_\mu \mathcal{M} \oplus T_{\tau_\mu}\mathcal{E}_\mu$, endowing $\mathcal{E}$ with the $GL_{\mathbf{1}}(V)$-invariant pseudo-Riemannian metric (invariant under constant non-zero rescalings $\tau_\mu \mapsto c\tau_\mu$)
\begin{equation}\label{eq:bundle_metric_G}
\mathcal{G}_{(\mu, \tau_\mu)}\big((\dot{\mu}, \dot{\tau}), (\dot{\mu}, \dot{\tau})\big) \coloneqq G_{IJ}(\mu; \tau_\mu)\dot{\mu}^I \dot{\mu}^J + h^{ij}(\mu; \tau_\mu)\dot{\tau}_i \dot{\tau}_j, 
\end{equation}
where $G_{IJ}(\mu; \tau_\mu)$ denotes $G_\mu(\tau_\mu)$ in coordinates on $\mathcal{M}$, and $GL_{\mathbf{1}}(V) \subset GL(V)$ is the subgroup fixing $\mathbf{1}$.

Inside $\mathcal{M}$, the unital associative algebras form the Maurer--Cartan locus
\begin{equation}\label{11/19/26.1}
\mathcal{Z} = \left\{ \mu \in \mathcal{M}: \, [\mu, \mu]_{\mathrm{G}} = 0 \text{ and } \mathbf{1} \cdot_\mu u = u \cdot_\mu \mathbf{1} = u \text{ for all } u \in V \right\}, 
\end{equation}
where $[\cdot, \cdot]_{\mathrm{G}}$ is the Gerstenhaber bracket, and the admissible trace bundle over $\mathcal{Z}$ is cut out inside $\mathcal{E}$ by the commutator-vanishing constraint \[ \mathcal{R}_{\mathrm{tr}}(\mu, \tau_\mu)(u, v) \coloneqq \tau_\mu(u \cdot_\mu v - v \cdot_\mu u) = 0.\] On this locus, $\Delta_\mu(\tau) \neq 0$ is equivalent to non-degeneracy of $\hat{\tau}(x, y) = \tau(xy)$. This makes $B_\mu$ symmetric Frobenius with $\tau \in \mathcal{T}_1(B_\mu)$. Indeed, for any non-zero $u \in \ker \hat{\tau}$ we have $\tau(u) = \hat{\tau}(u, \mathbf{1}) = 0$, while $\mathcal{R}_{\mathrm{tr}}(\mu, \tau) = 0$ yields $\tau(vu) = \tau(uv) = 0$ for all $v \in V$, so both terms on the right-hand side of (\ref{eq:g_sym_mu}) vanish. Hence $g_\mu(u, \cdot) = 0$, while Theorem~\ref{thm:pseudo_finsler} gives the converse, thereby confining the constrained flow to symmetric Frobenius algebras. Smooth co-evolving Maurer--Cartan paths $(\mu(t), \tau(t))$ are then governed by the constrained geodesic action
\begin{multline}
\label{eq:constrained_action}
S[\mu, \tau, \lambda, \Lambda] = \int_{t_0}^{t_1} \bigg( \frac{1}{2} G_{IJ}\big(\mu(t); \tau(t)\big) \dot{\mu}^I \dot{\mu}^J + \frac{1}{2} h^{ij}\big(\mu(t); \tau(t)\big) \dot{\tau}_i \dot{\tau}_j \\
+ \lambda^a(t) \Phi_a\big(\mu(t)\big) + \Lambda_b(t) \mathcal{R}_{\mathrm{tr}}^b\big(\mu(t), \tau(t)\big) \bigg) dt, 
\end{multline}
where $\Phi_a(\mu) = 0$ collects the Gerstenhaber associativity and unitality equations and $\lambda^a(t), \Lambda_b(t)$ are Lagrange multipliers. 

\subsubsection{Derivation of the Constrained Bundle Geodesic Equations}
To express the Euler--Lagrange equations of (\ref{eq:constrained_action}) compactly on $\mathcal{E} \subset \mathcal{M} \times V^*$, concatenate the $n^3$ base structure constants $\mu = (\mu^I)$ and the $n$ dual-basis fiber coordinates $\tau = (\tau_i)$ into the combined bundle coordinate tuple
\begin{equation*}
q = (q^1, \dots, q^{n^3+n}) \coloneqq (\mu^1, \dots, \mu^{n^3}, \, \tau_1, \dots, \tau_n),
\end{equation*}
where capital indices $A, B, C, D \in \{1, \dots, n^3 + n\}$ run over both blocks ($q^I = \mu^I$ for $1 \le I \le n^3$ and $q^{n^3+i} = \tau_i$ for $1 \le i \le n$). Because $\Delta_\mu(\tau) \neq 0$ on $\mathcal{E}$, the ambient bundle metric (\ref{eq:bundle_metric_G}) and its inverse take the non-singular block-diagonal forms
\begin{equation*}
\big[\mathcal{G}_{AB}(q)\big] = \begin{pmatrix} G_{IJ}(\mu; \tau) & 0 \\ 0 & h^{ij}(\mu; \tau) \end{pmatrix}, \qquad
\big[\mathcal{G}^{AB}(q)\big] = \begin{pmatrix} G^{IJ}(\mu; \tau) & 0 \\ 0 & h_{ij}(\mu; \tau) \end{pmatrix},
\end{equation*}
with $h_{ij}(\mu; \tau) = \tau(\mathbf{1})^2 g_{ij}(\mu; \tau)$, where both blocks depend simultaneously on $\mu$ and $\tau$. Collecting the constraints into $\tilde{\Phi}(q) = (\tilde{\Phi}_r(q)) \coloneqq \big(\Phi_a(\mu), \, \mathcal{R}_{\mathrm{tr}}^b(\mu, \tau)\big)$ and the Lagrange multipliers into $\Xi(t) = (\Xi^r(t)) \coloneqq \big(\lambda^a(t), \, \Lambda_b(t)\big)$ condenses the integrand of (\ref{eq:constrained_action}) to
\begin{equation*}
L(q, \dot{q}, \Xi) = \frac{1}{2}\mathcal{G}_{BC}(q)\dot{q}^B \dot{q}^C + \Xi^r(t)\tilde{\Phi}_r(q).
\end{equation*}

Stationarity of (\ref{eq:constrained_action}) with respect to $\Xi^r$ enforces $\tilde{\Phi}(q(t)) = 0$. Writing $\partial_A \coloneqq \partial / \partial q^A$ and using the symmetry of $\mathcal{G}_{BC}$ and $\dot{q}^B \dot{q}^C$, we compute for each $D \in \{1, \dots, n^3+n\}$,
\begin{align*}
\frac{d}{dt}\left(\frac{\partial L}{\partial \dot{q}^D}\right) &= \frac{d}{dt}\big(\mathcal{G}_{DB}(q)\dot{q}^B\big) = \mathcal{G}_{DB}(q)\ddot{q}^B + \frac{1}{2}\big(\partial_C \mathcal{G}_{DB} + \partial_B \mathcal{G}_{DC}\big)\dot{q}^B \dot{q}^C, \\
\frac{\partial L}{\partial q^D} &= \frac{1}{2}(\partial_D \mathcal{G}_{BC})\dot{q}^B \dot{q}^C + \Xi^r(t)\,\partial_D \tilde{\Phi}_r.
\end{align*}
Equating these two expressions and contracting with $\mathcal{G}^{AD}(q)$ yields the constrained bundle geodesic equations
\begin{equation}
\label{eq:EL_constrained}
\ddot{q}^A + \Gamma^A{}_{BC}(q) \dot{q}^B \dot{q}^C = \mathcal{G}^{AD}(q) \Xi^r(t) \frac{\partial \tilde{\Phi}_r}{\partial q^D},
\end{equation}
where $\Gamma^A{}_{BC}(q) \coloneqq \frac{1}{2}\mathcal{G}^{AD}(q)\big(\partial_C \mathcal{G}_{DB} + \partial_B \mathcal{G}_{DC} - \partial_D \mathcal{G}_{BC}\big)$ are the Levi-Civita Christoffel symbols of $(\mathcal{E}, \mathcal{G})$, whose non-vanishing mixed base-fiber components encode the coupling between $\mu(t)$ and $\tau(t)$.

\subsubsection{Unobstructed Sheets, Trace Compatibility, and Orthogonal Splitting}
Because $\mathcal{Z}$ and $\dim[B_\mu, B_\mu]$ can exhibit algebraic singularities or rank jumps at non-rigid algebras, as when a Gerstenhaber--Kuranishi obstruction in $HH^3(B_{\mu_0}, B_{\mu_0})$ cuts $\mathcal{Z}$ into a singular cone at $\mu_0$~\cite{Gerstenhaber:1964,Manetti:2004}, we restrict to a real-analytic local constrained submanifold $\mathcal{E}_{\mathrm{reg}} \subset \mathcal{E} \cap \tilde{\Phi}^{-1}(0)$ through $(\mu_0, \tau_0)$ fibering over a real-analytic sheet $\mathcal{Z}_0 \subset \mathcal{Z}$ (taken either on the full admissible trace space when $\dim[B_\mu, B_\mu]$ is constant on $\mathcal{Z}_0$ or along a real-analytic subbundle of extendable traces across a commutator-rank jump, over a $GL_{\mathbf{1}}(V)$-saturated open unobstructed sheet or a real-analytic slice transverse to the coboundary orbit as in Sections~\ref{sec:nakayama_geodesics}--\ref{sec:A3_spray}). Whenever $\bar{\mathcal{G}} \coloneqq \mathcal{G}|_{T\mathcal{E}_{\mathrm{reg}}}$ is non-degenerate, we define a \emph{constrained geodesic} on $\mathcal{E}_{\mathrm{reg}}$ as a stationary curve of the action (\ref{eq:constrained_action}) with respect to variations tangent to $\mathcal{E}_{\mathrm{reg}} \subset \mathcal{E} \cap \tilde{\Phi}^{-1}(0)$, i.e., an intrinsic Levi-Civita geodesic of $(\mathcal{E}_{\mathrm{reg}}, \bar{\mathcal{G}})$. When $\mathcal{E}_{\mathrm{reg}}$ is an open subset of $\mathcal{E} \cap \tilde{\Phi}^{-1}(0)$ cut out cleanly with $\ker d\tilde{\Phi} = T\mathcal{E}_{\mathrm{reg}}$, $(T\mathcal{E}_{\mathrm{reg}})^{\perp_{\mathcal{G}}}$ is spanned by the constraint gradients $\mathcal{G}^{AD}\partial_D\tilde{\Phi}_r$ and intersects $T\mathcal{E}_{\mathrm{reg}}$ trivially, so every constrained geodesic on $\mathcal{E}_{\mathrm{reg}}$ equivalently solves the ambient multiplier system (\ref{eq:EL_constrained}) for unique real-analytic multipliers $\Xi^r(t)$ (up to linear dependencies among $\partial\tilde{\Phi}_r$).

Along any smooth curve $(\mu(t), \tau(t))$ in $\mathcal{E}_{\mathrm{reg}}$ with $(\mu(0), \tau(0)) = (\mu_0, \tau_0)$ and initial velocity $(v_0, w_0) \coloneqq (\dot{\mu}(0), \dot{\tau}(0))$, differentiating unitality $\mu(t)(\mathbf{1}, u) = \mu(t)(u, \mathbf{1}) = u$ at $t = 0$ forces $v_0(\mathbf{1}, \cdot) = v_0(\cdot, \mathbf{1}) = 0$, placing $v_0$ in the normalized Hochschild cochain complex $C^\bullet(B_{\mu_0}, B_{\mu_0})$~\cite{Gerstenhaber:1964}. By standard deformation theory, differentiating associativity $[\mu(t), \mu(t)]_{\mathrm{G}} = 0$ places $v_0$ in $T_{\mu_0}\mathcal{Z}_0 \subset Z^2(B_{\mu_0}, B_{\mu_0}) = \ker d_{\mu_0}^{(2)}$, while differentiating the unital basis-change action $g(t)\cdot\mu_0$ for $g(t) \in GL_{\mathbf{1}}(V)$ identifies the tangent space to the isomorphism orbit $\mathcal{O}_{\mu_0} = GL_{\mathbf{1}}(V)\cdot\mu_0$ with the $2$-coboundary space $B^2(B_{\mu_0}, B_{\mu_0}) = \operatorname{Im}(d_{\mu_0}^{(1)})$ (where $C^1(B_{\mu_0}, B_{\mu_0}) = \mathfrak{gl}_{\mathbf{1}}(V)$). Simultaneously, differentiating $\tau(t)\big(\mu(t)(x, y) - \mu(t)(y, x)\big) = 0$ at $t = 0$ via the bilinear product rule subjects $(v_0, w_0)$ to the linearized trace-compatibility condition
\begin{equation*}
w_0([x, y]_{\mu_0}) + \tau_0\big(v_0(x, y) - v_0(y, x)\big) = 0 \qquad \text{for all } x, y \in V.
\end{equation*}

Finally, note that although $g_{\mu_0}(\tau_0)$ is typically indefinite when $\operatorname{rad}(B_{\mu_0}) \neq \{0\}$, dualizing the $g_{\mu_0}$-orthogonal splitting $V = \mathbb{R}\mathbf{1} \oplus \ker\tau_0$ shows that $g_{\mu_0}^{-1}$ restricts non-degenerately to $(\mathbb{R}\mathbf{1})^0$, rendering $G_{\mu_0}(\tau_0)$ globally non-degenerate on $C^2(B_{\mu_0}, B_{\mu_0}) = V \otimes (\mathbb{R}\mathbf{1})^0 \otimes (\mathbb{R}\mathbf{1})^0$ for every $\tau_0 \in \mathcal{T}_1(B_{\mu_0})$. Recall that on a finite-dimensional space with a symmetric bilinear form, a subspace $W$ satisfies $W \cap W^\perp = \{0\}$ if and only if the form restricts non-degenerately to $W$. Setting $B^2_{\mathcal{Z}_0} \coloneqq \operatorname{Im}(d_{\mu_0}^{(1)}) \cap T_{\mu_0}\mathcal{Z}_0$, whenever $\bar{\mathcal{G}}_{(\mu_0, \tau_0)}$ and the restriction of $\bar{G}_{\mu_0}(\tau_0) \coloneqq G_{\mu_0}(\tau_0)|_{T_{\mu_0}\mathcal{Z}_0}$ to $B^2_{\mathcal{Z}_0}$ are non-degenerate---rational determinant conditions on $\tau_0$ that hold on a Zariski-open dense subset of $\mathcal{T}_1(B_{\mu_0})$ whenever satisfied at a single trace, as verified explicitly in Sections~\ref{sec:2d_deformations}--\ref{sec:A3_spray}---we obtain the $\bar{G}_{\mu_0}(\tau_0)$-orthogonal Hodge-type splitting
\begin{align}\label{eq:hodge_splitting}
T_{\mu_0}\mathcal{Z}_0 &= B^2_{\mathcal{Z}_0} \oplus \mathcal{H}^2_{\mu_0}(\tau_0), \\
\mathcal{H}^2_{\mu_0}(\tau_0) &\coloneqq \big(B^2_{\mathcal{Z}_0}\big)^{\perp_{\bar{G}_{\mu_0}(\tau_0)}} \cap T_{\mu_0}\mathcal{Z}_0 \hookrightarrow HH^2(B_{\mu_0}, B_{\mu_0}), \nonumber
\end{align}
where $\mathcal{H}^2_{\mu_0}(\tau_0)$ intersects $\operatorname{Im}(d_{\mu_0}^{(1)})$ trivially and embeds into $HH^2(B_{\mu_0}, B_{\mu_0})$ as the subspace of $\mathcal{Z}_0$-unobstructed cohomology classes (canonically whenever $\mathcal{Z}_0$ is $GL_{\mathbf{1}}(V)$-saturated), with $\mathcal{H}^2_{\mu_0}(\tau_0) \cong HH^2(B_{\mu_0}, B_{\mu_0})$ whenever $Z^2(B_{\mu_0}, B_{\mu_0}) = \operatorname{Im}(d_{\mu_0}^{(1)}) + T_{\mu_0}\mathcal{Z}_0$.

\subsubsection{The ART-inspired computation of Maurer-Cartan paths}

\begin{theorem}
\label{thm:deformation_geodesics}
Suppose a real-analytic constrained submanifold $\mathcal{E}_{\mathrm{reg}} \subset \mathcal{E} \cap \tilde{\Phi}^{-1}(0)$ fibering over a real-analytic sheet $\mathcal{Z}_0 \subset \mathcal{Z}$ contains a smooth path $(\tilde{\mu}(s), \tilde{\tau}(s))$ through $(\tilde{\mu}(0), \tilde{\tau}(0)) = (\mu_0, \tau_0)$ with $\tilde{\mu}'(0) \notin \operatorname{Im}(d_{\mu_0}^{(1)})$, i.e., an infinitesimally non-trivial trace-compatible Maurer--Cartan path. Assume further that $\bar{\mathcal{G}} = \mathcal{G}|_{T\mathcal{E}_{\mathrm{reg}}}$ is non-degenerate at $(\mu_0, \tau_0)$ and that $\bar{G}_{\mu_0}(\tau_0) = G_{\mu_0}(\tau_0)|_{T_{\mu_0}\mathcal{Z}_0}$ restricts non-degenerately to $B^2_{\mathcal{Z}_0} = \operatorname{Im}(d_{\mu_0}^{(1)}) \cap T_{\mu_0}\mathcal{Z}_0$. Then $\mathcal{H}^2_{\mu_0}(\tau_0) \neq \{0\}$ in~\textup{(\ref{eq:hodge_splitting})}, and for every initial velocity $(v_0, w_0) \in T_{(\mu_0, \tau_0)}\mathcal{E}_{\mathrm{reg}}$ whose base projection $v_0$ has a non-zero component in $\mathcal{H}^2_{\mu_0}(\tau_0)$, we have
\begin{enumerate}
    \item[\textup{(i)}] On the open locus $\mathcal{E}_{\mathrm{reg}}^{\mathrm{nd}} \coloneqq \{ q \in \mathcal{E}_{\mathrm{reg}} : \det\bar{\mathcal{G}}_q \neq 0 \}$, there exists a unique maximal real-analytic constrained geodesic $(\mu(t), \tau(t))$ initialized by $\big((\mu_0, \tau_0), (v_0, w_0)\big)$ whose base projection $\mu(t) \in \mathcal{Z}$ is a Maurer--Cartan path of unital associative algebras co-evolving with $\tau(t) \in \mathcal{T}_1(B_{\mu(t)})$, defined either for all $t$ or until approaching $\partial\mathcal{E}_{\mathrm{reg}}$, norm infinity, a trace horizon $\Delta_{\mu(t)}(\tau(t)) = 0$, $\tau(t)(\mathbf{1}) = 0$, or pullback-metric degeneracy $\det\bar{\mathcal{G}} = 0$.
    \item[\textup{(ii)}] There is an open interval $(-\epsilon, \epsilon)$ around $t = 0$ on which $\mu(t)$ represents a non-trivial deformation of $B_{\mu_0}$, remaining outside the $GL_{\mathbf{1}}(V)$-isomorphism orbit of $\mu_0$ for all $t \in (-\epsilon, \epsilon) \setminus \{0\}$.
\end{enumerate}
\end{theorem}

\begin{proof}
Because $\tilde{\mu}'(0) \in T_{\mu_0}\mathcal{Z}_0 \setminus \operatorname{Im}(d_{\mu_0}^{(1)})$, its projection onto $\mathcal{H}^2_{\mu_0}(\tau_0)$ across~\textup{(\ref{eq:hodge_splitting})} is non-zero, confirming $\mathcal{H}^2_{\mu_0}(\tau_0) \neq \{0\}$.

For (i), $\mathcal{E}_{\mathrm{reg}}^{\mathrm{nd}}$ is a non-empty open neighborhood of $(\mu_0, \tau_0)$ in $\mathcal{E}_{\mathrm{reg}}$ on which $(\mathcal{E}_{\mathrm{reg}}^{\mathrm{nd}}, \bar{\mathcal{G}})$ is a real-analytic pseudo-Riemannian manifold whose intrinsic geodesic equation is a real-analytic autonomous second-order ODE. By the Cauchy--Kovalevskaya and Picard--Lindel\"of theorems together with standard maximal extension of ODE solutions, the initial data $\big((\mu_0, \tau_0), (v_0, w_0)\big) \in T\mathcal{E}_{\mathrm{reg}}^{\mathrm{nd}}$ determines a unique maximal real-analytic geodesic $(\mu(t), \tau(t))$ in $\mathcal{E}_{\mathrm{reg}}^{\mathrm{nd}}$ (leaving every compact subset of $\mathcal{E}_{\mathrm{reg}}^{\mathrm{nd}}$ whenever an endpoint is finite) along which $\Phi(\mu(t)) = 0$ and $\mathcal{R}_{\mathrm{tr}}(\mu(t), \tau(t)) = 0$ hold identically.  Thus, $\mu(t) \in \mathcal{Z}$ satisfies unitality and $[\mu(t), \mu(t)]_{\mathrm{G}} = 0$ (equivalently, $\gamma(t) \coloneqq \mu(t) - \mu_0$ solves $d_{\mu_0}^{(2)}\gamma(t) + \frac{1}{2}[\gamma(t), \gamma(t)]_{\mathrm{G}} = 0$) while co-evolving with $\tau(t) \in \mathcal{T}_1(B_{\mu(t)})$.

For (ii), the algebraic group action of $GL_{\mathbf{1}}(V)$ on $\mathcal{M}$ renders $\mathcal{O}_{\mu_0} = GL_{\mathbf{1}}(V) \cdot \mu_0$ a locally closed immersed submanifold, so for a sufficiently small open neighborhood $U \subset \mathcal{M}$ of $\mu_0$, $\mathcal{O}_{\mu_0} \cap U$ is an embedded submanifold of $U$ with $T_{\mu_0}(\mathcal{O}_{\mu_0} \cap U) = \operatorname{Im}(d_{\mu_0}^{(1)})$. By (\ref{eq:hodge_splitting}), the non-zero projection of $v_0 = \dot{\mu}(0)$ onto $\mathcal{H}^2_{\mu_0}(\tau_0)$ implies $v_0 \notin \operatorname{Im}(d_{\mu_0}^{(1)})$ since $\mathcal{H}^2_{\mu_0}(\tau_0) \cap \operatorname{Im}(d_{\mu_0}^{(1)}) = \{0\}$, so the Inverse Function Theorem guarantees an interval $(-\epsilon, \epsilon)$ on which $\mu(t) \notin \mathcal{O}_{\mu_0}$ for all $t \in (-\epsilon, \epsilon) \setminus \{0\}$.
\end{proof}

Note that the $GL_{\mathbf{1}}(V)$-equivariance of $\mathcal{G}$ ensures that any unital isomorphism $\varphi \colon B_{\mu_1} \xrightarrow{\sim} B_{\mu_2}$ intertwines $\mathcal{H}^2_{\mu_1}(\varphi^*\tau_2)$ with $\mathcal{H}^2_{\mu_2}(\tau_2)$ and carries the resulting constrained geodesics, trace horizons, and focal caustics isometrically onto one another.

We now consider algorithmic aspects of how the hypotheses of Theorem~\ref{thm:deformation_geodesics} can be verified to be satisfied or not from $(\mu_0, \tau_0)$ alone, without knowing any Maurer--Cartan path in advance. Note first that mere existence of an associative path $\tilde{\mu}(s) \in \mathcal{Z}$ does not guarantee an accompanying trace path $\tilde{\tau}(s)$, since higher-order terms can break the symmetric Frobenius condition (for instance, twisting a commutative factor $\mathbb{R}[x,y]/(x^2, y^2)$ into $\mathbb{R}\langle x,y\rangle/(x^2, y^2, yx - (1+s^2)xy)$ at $O(s^2)$ forces $[x,y]_s = -s^2 xy$ into the socle and empties $\mathcal{T}_1$ for $s \neq 0$). Assuming rational (or real-algebraic) entries, one therefore tests the joint quadratic polynomial system $\tilde{\Phi}(\mu, \tau) = \big(\Phi_a(\mu), \mathcal{R}_{\mathrm{tr}}^b(\mu, \tau)\big) = 0$ directly at $(\mu_0, \tau_0)$ in two stages:
\begin{enumerate}
    \item[\textup{(1)}] \emph{Real-analytic constrained submanifold existence (finite-time real algebraic geometry).}  Form the linearized constraint kernel
    \begin{equation*}
    K(\mu_0, \tau_0) \coloneqq \ker d\tilde{\Phi}|_{(\mu_0, \tau_0)} = \Big\{ (v, w) \in Z^2(B_{\mu_0}, B_{\mu_0}) \oplus V^* \;\Big|\; w([x, y]_{\mu_0}) + \tau_0\big(v(x, y) - v(y, x)\big) = 0 \Big\}.
    \end{equation*}
    A necessary linear screen is $\pi_1\big(K(\mu_0, \tau_0)\big) \not\subset \operatorname{Im}(d_{\mu_0}^{(1)})$. Whenever $\tilde{\Phi}^{-1}(0)$ is cleanly cut out near $(\mu_0, \tau_0)$ of local dimension $\dim K(\mu_0, \tau_0)$ with constant-rank base projection $\pi \colon \tilde{\Phi}^{-1}(0) \to \mathcal{M}$ (certifiable via local Gr\"obner bases and Jacobian minors, and holding in Sections~\ref{sec:2d_deformations}--\ref{sec:A3_spray}), $\mathcal{E}_{\mathrm{reg}} = \tilde{\Phi}^{-1}(0) \cap \mathcal{E}$ is automatically a real-analytic constrained submanifold fibering over a real-analytic sheet $\mathcal{Z}_0 \subset \mathcal{Z}$ with $T_{(\mu_0, \tau_0)}\mathcal{E}_{\mathrm{reg}} = K(\mu_0, \tau_0)$ and $T_{\mu_0}\mathcal{Z}_0 = \pi_1\big(K(\mu_0, \tau_0)\big)$, guaranteeing $(\tilde{\mu}(s), \tilde{\tau}(s))$ without computing the path. More generally, when $\tilde{\Phi}^{-1}(0)$ is singular at $(\mu_0, \tau_0)$, Artin's Nash approximation theorem~\cite{Artin1969}, constructive local resolution of singularities, and Tarski--Seidenberg quantifier elimination / Cylindrical Algebraic Decomposition decide in finite time (albeit with at least exponential worst-case complexity in $n^3$) whether a smooth irreducible component or one-dimensional real-analytic (Nash) branch $\mathcal{Z}_0 \subset \mathcal{Z}$ through $\mu_0$ with $T_{\mu_0}\mathcal{Z}_0 \not\subset \operatorname{Im}(d_{\mu_0}^{(1)})$ supporting a real-analytic trace subbundle $\mathcal{E}_{\mathrm{reg}} \subset \mathcal{E} \cap \tilde{\Phi}^{-1}(0)$ through $(\mu_0, \tau_0)$ exists, and extract their tangent spaces at $(\mu_0, \tau_0)$.
    \item[\textup{(2)}] \emph{Metric non-degeneracy (polynomial-time linear algebra):} Once $T_{\mu_0}\mathcal{Z}_0$ and $T_{(\mu_0, \tau_0)}\mathcal{E}_{\mathrm{reg}}$ are in hand, choose any basis $\{e_1 = \mathbf{1}, e_2, \dots, e_n\}$ of $V$, invert $[g_{ij}(\mu_0; \tau_0)]$ in~\textup{(\ref{eq:g_sym_mu})}, extract a basis $\{b_1, \dots, b_r\}$ of $B^2_{\mathcal{Z}_0} = \operatorname{Im}(d_{\mu_0}^{(1)}) \cap T_{\mu_0}\mathcal{Z}_0$, and test $\det M_{\mathrm{cob}}(\tau_0) \neq 0$, where $M_{\mathrm{cob}}(\tau_0)_{pq} \coloneqq G_{\mu_0}(b_p, b_q; \tau_0)$. Similarly, extract a basis $\{(v_1, w_1), \dots, (v_m, w_m)\}$ of $T_{(\mu_0, \tau_0)}\mathcal{E}_{\mathrm{reg}}$ via Gaussian elimination and test $\det M_{\mathrm{bun}}(\tau_0) \neq 0$, where $M_{\mathrm{bun}}(\tau_0)_{\alpha\beta} \coloneqq G_{\mu_0}(v_\alpha, v_\beta; \tau_0) + h_{\mu_0}(w_\alpha, w_\beta; \tau_0)$. (Across a linear trace fiber $(\mathcal{E}_{\mathrm{reg}})_{\mu_0}$ with coordinates $u$, both determinants are rational in $u$ where the linearized trace constraint has constant rank, so non-vanishing at a single rational functional $\tau_0$ certifies both metric hypotheses on a Zariski-open dense subset of that fiber, while symbolic evaluation in $u$ decides identical vanishing.)
\end{enumerate}
When both stages succeed, integrating the intrinsic geodesic equation of $(\mathcal{E}_{\mathrm{reg}}, \bar{\mathcal{G}})$ constructs the real-analytic Maurer--Cartan path $\mu(t)$ and co-evolving trace $\tau(t)$ from initial data at $t = 0$.

\begin{remark}
\label{rem:MC_comparison}
It is instructive to note the complementary scopes of Theorem~\ref{thm:deformation_geodesics} and classical homological deformation theory~\cite{Gerstenhaber:1964,getzler:1993,Manetti:2004,Dolgushev:2011}. Homological perturbation and Kuranishi methods apply to arbitrary associative deformations in $\mathcal{Z}$ without metric or trace restrictions. Thus, they explicitly compute the $HH^3(B_{\mu_0}, B_{\mu_0})$ obstruction equations cutting out $\mathcal{Z}$ near $\mu_0$, and they construct formal paths even along symmetry-breaking branches on which $B_{\mu(t)}$ ceases to be symmetric Frobenius ($\mathcal{T}_1(B_{\mu(t)}) = \emptyset$). By contrast, Theorem~\ref{thm:deformation_geodesics} is subject to two structural hypotheses:
\begin{enumerate}
    \item[\textup{(a)}] \emph{Trace extendability:} It requires a real-analytic sheet $\mathcal{Z}_0 \subset \mathcal{Z}$ supporting a real-analytic submanifold of extendable admissible traces $\mathcal{E}_{\mathrm{reg}} \subset \mathcal{E} \cap \tilde{\Phi}^{-1}(0)$ (which holds across the full admissible trace fiber when $\dim[B_\mu, B_\mu]$ is constant on $\mathcal{Z}_0$, or along a real-analytic subbundle of surviving traces across an upward commutator-rank jump), confining $\mu(t)$ to symmetric Frobenius deformations (or deformations of the trivial extension $B \ltimes B^*$) up to the trace horizon $\Delta_\mu(\tau) = 0$.
    \item[\textup{(b)}] \emph{Indefinite pullback non-degeneracy:} Because $g_{\mu_0}(\tau_0)$ is generally indefinite when $\operatorname{rad}(B_{\mu_0}) \neq \{0\}$, non-degeneracy of the submanifold pullbacks $\bar{\mathcal{G}}_{(\mu_0, \tau_0)}$ and $\bar{G}_{\mu_0}(\tau_0)|_{B^2_{\mathcal{Z}_0}}$ is not automatic.  While
Zariski-open across a linear trace fiber whenever satisfied at a single trace, either determinant could vanish identically on $\mathcal{T}_1(B_{\mu_0})$ (precluding the orthogonal splitting~\textup{(\ref{eq:hodge_splitting})} or ODE reduction at launch), and launched trajectories can terminate in finite time at internal metric-degeneracy walls $\det\bar{\mathcal{G}} = 0$ (as at $3c_1 = 4u_2$ in Section~\ref{sec:A3_spray}) where $\mathcal{Z}_0$ itself remains smooth.
\end{enumerate}

Conversely, whenever $(\mu_0, \tau_0) \in \mathcal{E}_{\mathrm{reg}}^{\mathrm{nd}}$ satisfies $\det\big(\bar{G}_{\mu_0}(\tau_0)|_{B^2_{\mathcal{Z}_0}}\big) \neq 0$, the variational bundle flow upgrades recursive formal Maurer--Cartan integration and Gauss--Manin cyclic transport in four respects:
\begin{enumerate}
    \item[\textup{(1)}] \emph{Canonical gauge resolution:} Rather than making non-canonical $1$-cochain coboundary choices at every formal order $t^k$ (or choosing an arbitrary contracting homotopy), $\bar{G}_{\mu_0}(\tau_0)$-orthogonal projection onto $\mathcal{H}^2_{\mu_0}(\tau_0)$ fixes the gauge $GL_{\mathbf{1}}(V)$-equivariantly at $t = 0$, after which the autonomous geodesic ODE determines all higher-order jets uniquely along $\mathcal{Z}_0$.
    \item[\textup{(2)}] \emph{Real-analyticity:} Because $\mathcal{G}$ and $\tilde{\Phi}$ are rational and polynomial on $\mathcal{E}_{\mathrm{reg}}^{\mathrm{nd}}$, the Cauchy--Kovalevskaya/Picard--Lindel\"of theorem yields genuine convergent real-analytic Maurer--Cartan paths $\mu(t)$ rather than formal $\R[[t]]$-series.
    \item[\textup{(3)}] \emph{Finite-dimensional trace co-evolution:} The constraint $\mathcal{R}_{\mathrm{tr}}(\mu, \tau) = 0$ and fiber metric $h^{ij}$ co-evolve concrete non-degenerate Frobenius traces $\tau(t)$ in finite dimensions, bypassing Getzler's Gauss--Manin connection on the infinite-dimensional $(b, B)$ periodic cyclic bar bicomplex (even for algebras of infinite global dimension, such as $N_6$ in Section~\ref{sec:nakayama_geodesics}).
    \item[\textup{(4)}] \emph{Canonical geodesic sprays and global deformation invariants:} Classical cohomological methods lack an intrinsic metric or connection and thus cannot define a canonical geodesic spray or integrate formal germs to finite-distance horizons. Via the harmonic exponential map $\exp^\perp_{(\mu_0, \tau_0)}$ (detailed below), our variational metric launches $GL_{\mathbf{1}}(V)$-equivariant multi-parameter sprays across unobstructed cones in $\mathcal{H}^2_{\mu_0}(\tau_0)$, yielding global differential-geometric and topological algebra invariants---including sectional curvatures, Jacobi conjugate caustics, Kossowski metric-wall cusps, Morse/Maslov focal index jumps, stratum wall-crossing times, and horizon/discriminant monodromies.  The latter have no counterpart in formal cochain algebra.
\end{enumerate}
\end{remark}

\subsubsection{Emergence of canonical geodesic sprays and global deformation invariants}
We will now formally develop issue (4) in Remark~\ref{rem:MC_comparison}. The new invariants referred to there emerge from connecting Theorem~\ref{thm:deformation_geodesics} to the catastrophe perspective of Section~3. 

Whenever the fiber metric $h(\mu_0; \tau_0)$ restricts non-degenerately to $\mathcal{T}(B_{\mu_0})$, each compatible harmonic velocity $v_0 \in \mathcal{H}^2_{\mu_0}(\tau_0)$ pairs with a unique $h(\mu_0; \tau_0)$-orthogonal compatible fiber velocity $w_0 \in \mathcal{T}(B_{\mu_0})^{\perp_h}$; more generally, even if $h|_{\mathcal{T}(B_{\mu_0})}$ admits null directions, $w_0 = 0$ is canonically distinguished whenever $\tau_0\big(v_0(x,y) - v_0(y,x)\big) = 0$ (in particular, when $B_{\mu_0}$ is commutative, or along $v_{\mathrm{harm}}$ in Section~\ref{sec:nakayama_geodesics}), while otherwise $w_0$ may be chosen orthogonal to a fixed complement of $\operatorname{rad}\big(h|_{\mathcal{T}(B_{\mu_0})}\big)$. To formalize the resulting family of paths, recall that on the maximally defined star-shaped open neighborhood $\mathcal{D}_{(\mu_0, \tau_0)} \subset T_{(\mu_0, \tau_0)}\mathcal{E}_{\mathrm{reg}}$ of $(0, 0)$ for which the unique $\bar{\mathcal{G}}$-geodesic $t \mapsto \gamma_{(v_0, w_0)}(t)$ of Theorem~\ref{thm:deformation_geodesics}(i) with $\gamma_{(v_0, w_0)}(0) = (\mu_0, \tau_0)$ and $\dot{\gamma}_{(v_0, w_0)}(0) = (v_0, w_0)$ exists up to $t = 1$, the intrinsic pseudo-Riemannian exponential map of $(\mathcal{E}_{\mathrm{reg}}^{\mathrm{nd}}, \bar{\mathcal{G}})$ evaluates at unit affine parameter via $(v_0, w_0) \mapsto \gamma_{(v_0, w_0)}(1)$ (so each ray is traced out by $t \mapsto \gamma_{(v_0, w_0)}(t)$ via radial homogeneity $\gamma_{(v_0, w_0)}(t) = \gamma_{(t v_0, t w_0)}(1)$). Denoting by $v_0 \mapsto \big(v_0, w_0(v_0)\big) \in T_{(\mu_0, \tau_0)}\mathcal{E}_{\mathrm{reg}}$ the linear lift that pairs each unobstructed harmonic cocycle $v_0 \in \mathcal{H}^2_{\mu_0}(\tau_0)$ with its selected compatible fiber velocity $w_0(v_0)$, restricting the intrinsic exponential map to this coboundary-orthogonal harmonic subspace defines the \emph{harmonic constrained exponential map}
\begin{equation*}
\exp^\perp_{(\mu_0, \tau_0)} \colon \mathcal{D}^\perp_{(\mu_0, \tau_0)} \longrightarrow \mathcal{E}_{\mathrm{reg}}, \qquad \exp^\perp_{(\mu_0, \tau_0)}(v_0) \coloneqq \gamma_{(v_0, w_0(v_0))}(1),
\end{equation*}
where $\mathcal{D}^\perp_{(\mu_0, \tau_0)} \coloneqq \big\{ v_0 \in \mathcal{H}^2_{\mu_0}(\tau_0) : (v_0, w_0(v_0)) \in \mathcal{D}_{(\mu_0, \tau_0)} \big\}$ (intersected with the unobstructed Maurer--Cartan cone when $HH^3(B_{\mu_0}, B_{\mu_0}) \neq 0$). Launching the constrained geodesic spray $t \mapsto \exp^\perp_{(\mu_0, \tau_0)}(t v_0)$ across unobstructed directions in $\mathcal{H}^2_{\mu_0}(\tau_0) \hookrightarrow HH^2(B_{\mu_0}, B_{\mu_0})$ sweeps out a canonical ``deformation fan'' in $\mathcal{E}_{\mathrm{reg}}$, exhibiting a two-tier catastrophe hierarchy.

First, for fixed $\tau_0$, suppose the unobstructed Maurer--Cartan cone in $\mathcal{H}^2_{\mu_0}(\tau_0)$ has dimension $d \ge 2$, and consider a smooth one-parameter family of constrained geodesics $s \mapsto \gamma_s(t) = (\mu_s(t), \tau_s(t))$ in $(\mathcal{E}_{\mathrm{reg}}^{\mathrm{nd}}, \bar{\mathcal{G}})$ all launched from $\gamma_s(0) = (\mu_0, \tau_0)$ with smoothly varying initial velocities $\dot{\gamma}_s(0)$. Measuring the first-order separation between neighboring rays in the spray via the transverse variation vector field $J(t) \coloneqq \left.\partial_s \gamma_s(t)\right|_{s=0} \in T_{\gamma(t)}\mathcal{E}_{\mathrm{reg}}$ along the central ray $\gamma(t) = \gamma_0(t)$, differentiating the intrinsic geodesic equation $\bar{\nabla}_{\dot{\gamma}_s}\dot{\gamma}_s = 0$ with respect to $s$ at $s = 0$ yields the classical geodesic deviation (Jacobi) equation,
\begin{equation*}
\bar{\nabla}_{\dot{\gamma}} \bar{\nabla}_{\dot{\gamma}} J + \bar{R}(J, \dot{\gamma})\dot{\gamma} = 0, \qquad J(0) = 0, \quad \bar{\nabla}_{\dot{\gamma}}J(0) = (\delta v_0, \delta w_0),
\end{equation*}
where $\bar{\nabla}$ is the Levi-Civita connection of $\bar{\mathcal{G}}$, $\bar{R}(X, Y)Z \coloneqq \bar{\nabla}_X \bar{\nabla}_Y Z - \bar{\nabla}_Y \bar{\nabla}_X Z - \bar{\nabla}_{[X,Y]}Z$ is its Riemann curvature endomorphism, and $(\delta v_0, \delta w_0) = \left.\partial_s \dot{\gamma}_s(0)\right|_{s=0}$ is the angular perturbation of the launch velocity. Taking the $\bar{\mathcal{G}}$-inner product of this deviation equation with $J(t)$ shows that the convergence or divergence of neighboring Maurer--Cartan rays is controlled by the curvature term $\langle \bar{R}(J, \dot{\gamma})\dot{\gamma}, J \rangle_{\bar{\mathcal{G}}}$ (which forms the numerator of the intrinsic sectional curvature of the $2$-plane spanned by $\dot{\gamma}$ and~$J$).

In turn, because $(\mathcal{E}_{\mathrm{reg}}^{\mathrm{nd}}, \bar{\mathcal{G}})$ is an embedded pseudo-Riemannian submanifold of the ambient bundle $(\mathcal{E}, \mathcal{G})$, decomposing ambient covariant derivatives into tangential and normal parts via $\nabla^{\mathcal{G}}_X Y = \bar{\nabla}_X Y + \mathrm{II}(X, Y)$ expresses this intrinsic curvature via the Gauss submanifold curvature equation (the higher-codimension generalization of Gauss's \emph{Theorema Egregium}),
\begin{equation*}
\langle \bar{R}(J, \dot{\gamma})\dot{\gamma}, J \rangle_{\bar{\mathcal{G}}} = \langle R^{\mathcal{G}}(J, \dot{\gamma})\dot{\gamma}, J \rangle_{\mathcal{G}} + \langle \mathrm{II}(\dot{\gamma}, \dot{\gamma}), \mathrm{II}(J, J) \rangle_{\mathcal{G}} - \langle \mathrm{II}(\dot{\gamma}, J), \mathrm{II}(\dot{\gamma}, J) \rangle_{\mathcal{G}},
\end{equation*}
where $R^{\mathcal{G}}$ is the Riemann curvature of the unconstrained bundle metric (\ref{eq:bundle_metric_G}) and $\mathrm{II} \colon T\mathcal{E}_{\mathrm{reg}}^{\mathrm{nd}} \times T\mathcal{E}_{\mathrm{reg}}^{\mathrm{nd}} \to (T\mathcal{E}_{\mathrm{reg}}^{\mathrm{nd}})^{\perp_{\mathcal{G}}}$ is the normal second fundamental form encoding the extrinsic bending of the associativity and trace-cyclicity constraints inside $\mathcal{E}$. On a Riemannian chamber, whenever these combined ambient and extrinsic constraint curvatures render the sectional curvature sufficiently positive along $\gamma(t)$ (or more generally whenever the tidal operator $-\bar{R}(\cdot, \dot{\gamma})\dot{\gamma}$ is focusing), neighboring Maurer--Cartan rays can be pulled back together before reaching the trace horizon $\Delta_{\mu(t)}(\tau(t)) = 0$, producing conjugate points---that is, bundle caustics of the exponential map $\exp^\perp_{(\mu_0, \tau_0)}$ where $J(t_*) = 0$ for some $t_* > 0$ (or moduli caustics where $J(t_*)$ becomes tangent to the $GL_{\mathbf{1}}(V)$-isomorphism orbit).

Second, varying the projective initial trace $[\tau_0] \in \mathbb{P}\big(\mathcal{T}_1(B_{\mu_0})\big)$ (represented by $\bar{\tau}_0 = \tau_0/\tau_0(\mathbf{1})$) at fixed cohomology class $[\eta] \in HH^2(B_{\mu_0}, B_{\mu_0})$ modulates the trajectory through three distinct mechanisms:
\begin{enumerate}
    \item[\textup{(a)}] \emph{Tilting of the Harmonic Slice:} Varying $\tau_0$ changes $G_{\mu_0}(\tau_0)$ (and $\bar{G}_{\mu_0}(\tau_0)$ on a transverse slice), rotating the orthogonal complement $\mathcal{H}^2_{\mu_0}(\tau_0)$ and shifting $v_0(\tau_0) \in [\eta]$ by a $\tau_0$-dependent coboundary.
    \item[\textup{(b)}] \emph{Infinitesimal Trace-Compatibility Selection Rules:} In noncommutative algebras, $w_0([x,y]_{\mu_0}) + \tau_0\big(v_0(x,y) - v_0(y,x)\big) = 0$ couples $(v_0, w_0)$ directly to $\tau_0$, enforcing $\tau_0$-dependent selection rules on pairs commuting in $B_{\mu_0}$.
    \item[\textup{(c)}] \emph{Moduli-Space Curving and Caustic Unfolding:} Because $\bar{\nabla}$ and $\bar{R}$ depend nonlinearly on $\tau_0$, varying $\tau_0$ across $\operatorname{Aut}(B_{\mu_0})$-orbits alters $O(t^2)$ covariant acceleration---deflecting rays launched along the same class $[\eta]$ into non-isomorphic algebras when $d \ge 2$ and unfolding conjugate focal loci and trace horizons.
\end{enumerate}

Viewing the family of maps $(v_0; \bar{\tau}_0) \mapsto \exp^\perp_{(\mu_0, \bar{\tau}_0)}(v_0)$ across normalized initial calibrations $\bar{\tau}_0 = \tau_0/\tau_0(\mathbf{1})$, the unobstructed harmonic coordinates $v_0 \in \mathcal{D}^\perp_{(\mu_0, \bar{\tau}_0)} \subset \mathcal{H}^2_{\mu_0}(\bar{\tau}_0)$ thus act as transverse state variables tracing out Maurer--Cartan rays $t \mapsto \exp^\perp_{(\mu_0, \bar{\tau}_0)}(t v_0)$, while $\bar{\tau}_0 \in \mathbb{P}\big(\mathcal{T}_1(B_{\mu_0})\big)$ serves as an intrinsic unfolding control space, enriching $\dim HH^2(B_{\mu_0}, B_{\mu_0})$ with Lagrangian caustic classes, Morse/Maslov focal indices, and horizon monodromy invariants.

\subsection{The simplest nontrivial example: deformations of real two-dimensional unital associative algebras}
\label{sec:2d_deformations}

To illustrate this framework in an example, consider the moduli space of two-dimensional unital associative algebras over $\mathbb{R}$ with ordered basis $\{e_1, e_2\} = \{\mathbf{1}, e\}$. Any such algebra $B_\mu$ is automatically associative and commutative, and is completely determined by two structure constants $(\alpha, \beta) \in \mathbb{R}^2$ via the multiplication table
\begin{equation*}
    e^2 = \alpha e + \beta \mathbf{1}, 
\end{equation*}
corresponding to $c_{22}^1 = \beta$ and $c_{22}^2 = \alpha$ (with $c_{1i}^k = c_{i1}^k = \delta_i^k$). Because associativity and commutativity hold identically across this two-dimensional family, both the Gerstenhaber associativity constraint $\Phi_a(\mu) = 0$ and the commutator-vanishing trace-compatibility constraint $\mathcal{R}_{\mathrm{tr}}^b(\mu, \tau) = 0$ are satisfied identically for every $(\mu, \tau) \in \mathbb{R}^2 \times V^*$. Consequently, no constraint forces are required to maintain associativity or trace cyclicity, and the Lagrange multiplier terms in the constrained action (\ref{eq:constrained_action}) vanish identically in this two-dimensional unital algebra setting.

Parametrizing the dual space covector $\tau \in V^*$ via $\rho \coloneqq \tau(\mathbf{1}) \neq 0$, $\theta \coloneqq \tau(e)$, and the normalized trace ratio $u \coloneqq \theta / \rho$, the fundamental bilinear form \mbox{$g_{ij} \coloneq g_{\mathbf{1}}^{B_\mu}(e_i, e_j; \tau)$} from (\ref{eq:g_sym_mu}) and its inverse $g^{ij}$ evaluate to
\begin{equation}\label{eq:2d_g_matrices}
    \big[g_{ij}\big] = \begin{pmatrix} 1 & u \\ u & 2u^2 - \alpha u - \beta \end{pmatrix},  \qquad
    \big[g^{ij}\big] = \frac{1}{\Delta_\mu(\tau)} \begin{pmatrix} 2u^2 - \alpha u - \beta & -u \\ -u & 1 \end{pmatrix}, 
\end{equation}
where the trace-horizon determinant is the scalar function on the bundle $\mathcal{E}$ given by
\begin{equation}\label{eq:2d_Delta}
    \Delta_\mu(\tau) \coloneqq \det \big[g_{ij}\big] = u^2 - \alpha u - \beta = \frac{\theta^2 - \alpha\rho\theta - \beta\rho^2}{\rho^2}. 
\end{equation}
Here, the subscript $\mu$ records that the algebra structure constants $\mu = (\alpha, \beta)$ supply the coefficients of the determinant function $\Delta_\mu(\tau)$ evaluated on the covector $\tau = (\rho, \theta) \in V^*$. Whenever the dependence on $(\mu, \tau)$ is clear from context, we write $\Delta \coloneqq \Delta_\mu(\tau)$ as a shorthand. Clearing the conformal denominator $\rho^2$ yields the polynomial trace-horizon locus $\theta^2 - \alpha\rho\theta - \beta\rho^2 = 0$.

Because $e_1 = \mathbf{1}$ is fixed, tangent variations on the $(\alpha, \beta)$-plane have non-vanishing components $\dot{c}_{22}^1 = \dot{\beta}$ and $\dot{c}_{22}^2 = \dot{\alpha}$. Contracting indices via $G_{IJ}\dot{\mu}^I\dot{\mu}^J = (g^{22})^2 g_{kl}\dot{c}_{22}^k \dot{c}_{22}^l$ with $g^{22} = 1/\Delta$, the horizontal base metric tensor $G_{IJ}(\alpha, \beta; u)$ on the $(\alpha, \beta)$-plane (ordered as $(\alpha, \beta)$) and its inverse $G^{IJ}(\alpha, \beta; u)$ are
\begin{equation}\label{eq:2d_G_horizontal}
    \big[G_{IJ}\big] = \frac{1}{\Delta^2} \begin{pmatrix} \Delta + u^2 & u \\ u & 1 \end{pmatrix},  \qquad
    \big[G^{IJ}\big] = \Delta \begin{pmatrix} 1 & -u \\ -u & \Delta + u^2 \end{pmatrix}. 
\end{equation}
Since $\det[G_{IJ}] = \Delta^{-3}$, the horizontal metric on the $(\alpha,\beta)$-plane is Riemannian on the positive chamber $\Delta > 0$, pseudo-Riemannian (Lorentzian) on the negative chamber $\Delta < 0$, and singular along the trace horizon $\Delta = 0$.

\subsubsection{Cohomology Rigidity, Orthogonal Splitting, and Launching from the Dual Numbers}
Up to isomorphism, real two-dimensional unital associative algebras are classified by the algebraic discriminant $\Delta_{\mathrm{alg}} = \alpha^2 + 4\beta$ of the characteristic polynomial $x^2 - \alpha x - \beta = 0$, \begin{enumerate}
    \item \emph{Split Semisimple ($\mathbb{R} \oplus \mathbb{R}$, $\Delta_{\mathrm{alg}} > 0$):} Rigid, i.e., $HH^2(\mathbb{R} \oplus \mathbb{R}, \mathbb{R} \oplus \mathbb{R}) = 0$.
    \item \emph{Dual Numbers ($\mathbb{R}[\epsilon]/(\epsilon^2)$, $\Delta_{\mathrm{alg}} = 0$):} Non-rigid, with \[ \dim HH^2\big(\mathbb{R}[\epsilon]/(\epsilon^2), \mathbb{R}[\epsilon]/(\epsilon^2)\big) = 1.  \]
    \item \emph{Complex Numbers ($\mathbb{C}$, $\Delta_{\mathrm{alg}} < 0$):} Rigid.
\end{enumerate}
Because $\mathbb{R} \oplus \mathbb{R}$ and $\mathbb{C}$ are cohomologically rigid, no nontrivial deformation ray can be initiated from them. The dual numbers at $\mu_0 = (\alpha_0, \beta_0) = (0, 0)$ serve as the unique non-rigid launchpad at the origin of the $(\alpha, \beta)$-plane.

Under a unital basis change $e' = a\mathbf{1} + b e$ ($b \neq 0$) around the dual numbers $e^2 = 0$, an infinitesimal shift $a$ with $b = 1$ yields $(a\mathbf{1} + e)^2 = 2a(a\mathbf{1} + e) - a^2\mathbf{1}$, producing the first-order variation $(\dot{\alpha}, \dot{\beta}) = (2a, 0)$. Thus, at $\mu_0 = (0, 0)$ equipped with an initial covector $\tau_0 = (\rho_0, \theta_0)$ of trace ratio $u_0 = \theta_0/\rho_0 \neq 0$ (so that $\Delta_{\mu_0}(\tau_0) = u_0^2 > 0$ lies in the positive regular chamber), \begin{itemize}
    \item The $GL_{\mathbf{1}}(V)$ coboundary orbit tangent space $\operatorname{Im}(d_{\mu_0}^{(1)})$ is spanned by the $\alpha$-axis vector $v_{\mathrm{cob}} = (1, 0)$, on which $G_{\mu_0}(v_{\mathrm{cob}}, v_{\mathrm{cob}}; u_0) = 2/u_0^2 \neq 0$.
    \item By (\ref{eq:2d_G_horizontal}), a tangent vector $v = (\dot{\alpha}, \dot{\beta})$ satisfies the $G_{\mu_0}(\tau_0)$-orthogonality condition $G_{\mu_0}(v_{\mathrm{cob}}, v; u_0) = u_0^{-4}(2u_0^2 \dot{\alpha} + u_0 \dot{\beta}) = 0$ if and only if $\dot{\beta} = -2u_0 \dot{\alpha}$. Thus, the $G_{\mu_0}(\tau_0)$-orthogonal harmonic complement $\mathcal{H}^2_{\mu_0}(\tau_0) = \big(\operatorname{Im}(d_{\mu_0}^{(1)})\big)^{\perp_{G_{\mu_0}(\tau_0)}} \cap Z^2(B_{\mu_0}, B_{\mu_0})$ in (\ref{eq:hodge_splitting}) is spanned by
    \begin{equation*}
        v_{\mathrm{harm}} = (1, -2u_0). 
    \end{equation*}
    \item More generally, the coordinate $\beta$-axis direction $v_\beta = (0, 1) = \frac{1}{2u_0}(1, 0) - \frac{1}{2u_0}(1, -2u_0)$ has a strictly non-zero orthogonal projection onto $\mathcal{H}^2_{\mu_0}(\tau_0)$ whenever $u_0 \neq 0$, and therefore satisfies the transversality hypothesis of Theorem~\ref{thm:deformation_geodesics}.
\end{itemize}

\subsubsection{Instantaneous Horizontal Geometry and the Hyperbolic Plane at the Onset of the Maurer--Cartan Path}
In the general framework of Section~\ref{sec:def_metric}, a canonical deformation ray launched from $(\mu_0, \tau_0) \in \mathcal{E}_{\mathrm{reg}}$ along a cohomology class $[v_0] \in HH^2(B_{\mu_0}, B_{\mu_0})$ is initialized by pairing the harmonic representative $v_0 = (\dot{\alpha}(0), \dot{\beta}(0))$ with the $h(\mu_0; \tau_0)$-orthogonal compatible fiber velocity $w_0 = \dot{\tau}(0) \in V^*$ satisfying the linearized trace-compatibility constraint $w_0([x,y]_{\mu_0}) + \tau_0(v_0(x,y) - v_0(y,x)) = 0$. Unlike non-commutative deformations in dimensions $n \ge 3$---where non-vanishing commutators actively force $w_0 \neq 0$ and drive fiber rotation from the very first instant (as demonstrated explicitly in Section~\ref{sec:nakayama_geodesics} below)---every algebra in our two-dimensional unital family is commutative, so the linearized compatibility constraint is satisfied identically by the zero fiber velocity, \begin{equation*}
    w_0 = \dot{\tau}(0) = 0 \qquad \Longleftrightarrow \qquad \dot{\rho}(0) = 0 \quad \text{and} \quad \dot{u}(0) = 0. 
\end{equation*}
Consequently, either canonical Maurer--Cartan deformation ray in dimension two is launched \emph{purely horizontally} along the instantaneous constant-trace slice $u = u_0 \coloneqq \theta_0/\rho_0$.

Because the horizontal metric $G_{IJ}(\alpha,\beta; u)$ in (\ref{eq:2d_G_horizontal}) depends explicitly on $u$ through $\Delta = u^2 - \alpha u - \beta$, the horizontal slice $u = u_0$ is not a totally geodesic submanifold of the four-dimensional bundle $\mathcal{E}_{\mathrm{reg}}$; the non-vanishing metric gradient $\partial_u G_{IJ} \neq 0$ acts as a variational force in the Euler--Lagrange equations of (\ref{eq:constrained_action}) that accelerates $u(t)$ away from $u_0$ at second order, $u(t) = u_0 + \frac{1}{2}\ddot{u}(0)t^2 + O(t^3)$. However, because $\dot{u}(0) = 0$ and $\dot{\rho}(0) = 0$, all mixed fiber-base Christoffel cross-terms $2\Gamma^K{}_{Iu}\dot{\mu}^I\dot{u}$ in the horizontal geodesic equations vanish identically at $t = 0$. Thus, at the onset of the Maurer--Cartan path, the base projection $(\alpha(t), \beta(t))$ of the true co-evolving bundle geodesic agrees through second order (up to $O(t^3)$ error) with the horizontal slice geodesic $(\bar{\alpha}(t), \bar{\beta}(t))$ of the instantaneous fixed-$u_0$ metric, \begin{equation*}
    \alpha(t) = \bar{\alpha}(t) + O(t^3),  \qquad \beta(t) = \bar{\beta}(t) + O(t^3), 
\end{equation*}
sharing the same initial position, initial velocity, and initial covariant acceleration (geodesic curvature) at $t = 0$.

To determine this second-order horizontal geometry at the launch of the deformation path, we evaluate the horizontal slice metric on the $(\alpha, \beta)$-plane at the initial trace ratio (writing $u = u_0$ to streamline notation). From (\ref{eq:constrained_action}) and (\ref{eq:2d_G_horizontal}), the instantaneous horizontal motion is governed by the reduced two-dimensional Lagrangian $L_{2\mathrm{D}} = \frac{1}{2}G_{IJ}\dot{\mu}^I\dot{\mu}^J$,
\begin{equation}\label{eq:2d_L_horizontal}
    L_{2\mathrm{D}}(\alpha, \beta, \dot{\alpha}, \dot{\beta}; u) = \frac{1}{2\Delta^2}\Big( (\Delta + u^2)\dot{\alpha}^2 + 2u\dot{\alpha}\dot{\beta} + \dot{\beta}^2 \Big), 
\end{equation}
where $\Delta = u^2 - \alpha u - \beta$. Using $\partial_\alpha \Delta = -u$ and $\partial_\beta \Delta = -1$ at fixed $u$, direct computation of the horizontal Levi-Civita connection $\Gamma^K{}_{IJ} = \frac{1}{2}G^{KL}(\partial_I G_{JL} + \partial_J G_{IL} - \partial_L G_{IJ})$ from (\ref{eq:2d_G_horizontal}) yields the six horizontal Christoffel symbols in $(\alpha, \beta)$ coordinates,
\begin{equation}\label{eq:2d_christoffel}
\begin{aligned}
    \Gamma^\alpha{}_{\alpha\alpha} &= \frac{u}{\Delta},  & \Gamma^\alpha{}_{\alpha\beta} &= \frac{1}{2\Delta},  & \Gamma^\alpha{}_{\beta\beta} &= 0,  \\[4pt]
    \Gamma^\beta{}_{\alpha\alpha} &= -\frac{1}{2},  & \Gamma^\beta{}_{\alpha\beta} &= \frac{u}{2\Delta},  & \Gamma^\beta{}_{\beta\beta} &= \frac{1}{\Delta}. 
\end{aligned}
\end{equation}
Notice the clean pole structure of (\ref{eq:2d_christoffel}): despite the $\Delta^{-2}$ factor in the metric tensor (\ref{eq:2d_G_horizontal}) and the $\Delta^{-3}$ derivatives of its individual entries, off-diagonal cancellations against $G^{KL}$ collapse all higher poles, leaving four components with simple $\Delta^{-1}$ poles along the trace horizon $\Delta = 0$, one constant component $\Gamma^\beta{}_{\alpha\alpha} = -1/2$, and $\Gamma^\alpha{}_{\beta\beta} = 0$.

In our setting where all algebras are commutative and associative, $\Xi$ vanishes from (\ref{eq:EL_constrained}). Thus, the latter equation becomes $\ddot{q}^K + \Gamma^K{}_{IJ}\dot{q}^I\dot{q}^J = 0$. Substituting (\ref{eq:2d_christoffel}) into that equation gives the coupled second-order horizontal slice ODEs on the $(\alpha, \beta)$-plane,
\begin{subequations}\label{eq:2d_horizontal_odes}
\begin{align}
    \ddot{\alpha} + \frac{1}{\Delta}\Big( u\dot{\alpha}^2 + \dot{\alpha}\dot{\beta} \Big) &= 0,  \\
    \ddot{\beta} - \frac{1}{2}\dot{\alpha}^2 + \frac{1}{\Delta}\Big( u\dot{\alpha}\dot{\beta} + \dot{\beta}^2 \Big) &= 0. 
\end{align}
\end{subequations}

Notably, the horizontal slice system (\ref{eq:2d_horizontal_odes}) is completely integrable via the affine coordinate change $(\alpha, \beta) \mapsto (\alpha, \Delta)$ with $\dot{\beta} = -u\dot{\alpha} - \dot{\Delta}$, which diagonalizes the reduced Lagrangian (\ref{eq:2d_L_horizontal}) and horizontal line element to
\begin{equation}\label{eq:2d_poincare_metric}
    L_{2\mathrm{D}}(\alpha, \Delta, \dot{\alpha}, \dot{\Delta}) = \frac{\dot{\alpha}^2}{2\Delta} + \frac{\dot{\Delta}^2}{2\Delta^2},  \qquad ds^2_{2\mathrm{D}} = \frac{d\alpha^2}{\Delta} + \frac{d\Delta^2}{\Delta^2}. 
\end{equation}
On the positive Riemannian chamber $\Delta > 0$, introducing $x = \alpha/2$ and $y = \sqrt{\Delta} > 0$ transforms (\ref{eq:2d_poincare_metric}) into $ds^2_{2\mathrm{D}} = 4(dx^2 + dy^2)/y^2$, demonstrating that the instantaneous horizontal slice on the positive chamber is isometric to the Poincar\'e upper half-plane model of the hyperbolic plane $\mathbb{H}^2$ of constant negative Gaussian curvature $K = -1/4$. On the negative Lorentzian chamber $\Delta = -y^2 < 0$, $ds^2_{2\mathrm{D}} = 4(-dx^2 + dy^2)/y^2$ is isometric to the Lorentzian hyperbolic plane $\mathbb{H}^{1,1}$. Here, $\mathbb{H}^2$ and $\mathbb{H}^{1,1}$ denote real hyperbolic spaces, not to be confused with the quaternion algebra $\mathbb{H}$. Because $\alpha$ is an ignorable coordinate of $L_{2\mathrm{D}}$ (i.e., $\partial_\alpha L_{2\mathrm{D}} = 0$), the instantaneous horizontal equations (\ref{eq:2d_horizontal_odes}) decouple directly into
\begin{equation}\label{eq:2d_decoupled_odes}
    \ddot{\alpha} - \frac{\dot{\alpha}\dot{\Delta}}{\Delta} = 0 \quad \Longrightarrow \quad \dot{\alpha}(t) = P_\alpha\,\Delta(t),  \qquad \text{and} \qquad \ddot{\Delta} - \frac{\dot{\Delta}^2}{\Delta} + \frac{1}{2}P_\alpha^2\Delta^2 = 0, 
\end{equation}
where $P_\alpha = \dot{\alpha}(0)/\Delta(0)$ is the horizontal Noether momentum and $E = \frac{1}{2}P_\alpha^2\Delta + \frac{\dot{\Delta}^2}{2\Delta^2}$ is the horizontal kinetic energy. For $P_\alpha = 0$ ($\dot{\alpha}(0) = 0$), (\ref{eq:2d_decoupled_odes}) integrates to $\bar{\alpha}(t) \equiv \alpha_0$ and $\bar{\Delta}(t) = \Delta_0 e^{\pm \sqrt{2E}\,t}$ (consistent with $\Gamma^\alpha{}_{\beta\beta} = 0$ in (\ref{eq:2d_christoffel})), while for $P_\alpha \neq 0$ it integrates in closed form to the hyperbolic secant profile
\begin{equation*}
    \bar{\Delta}(t) = \frac{2E}{P_\alpha^2 \cosh^2\!\big(\sqrt{E/2}\,(t - t_*)\big)},  \qquad \bar{\alpha}(t) = \alpha_* + \frac{2\sqrt{2E}}{P_\alpha}\tanh\!\big(\sqrt{E/2}\,(t - t_*)\big). 
\end{equation*}
Within this instantaneous horizontal metric, the $\Delta^{-2}$ blow-up of the line element (\ref{eq:2d_poincare_metric}) and $\Delta^{-1}$ divergence of the Christoffel symbols (\ref{eq:2d_christoffel}) erect an infinite metric barrier along the trace horizon $\Delta = 0$, rendering the positive horizontal slice geodesically complete and placing the horizon line $u^2 - \alpha u - \beta = 0$ at the conformal boundary at infinity ($t \to \pm\infty$).

\subsubsection{Full Co-Evolving 4D Bundle Geodesics}
While the fixed-$u_0$ Poincar\'e half-plane geometry governs the second-order horizontal trajectory at the onset of either canonical Maurer--Cartan deformation ray, continuing the trajectory to finite times $t > 0$ requires accounting for the $O(t^2)$ variational backreaction of the evolving algebra on the trace fiber, which pulls the trajectory off the initial $u_0$-slice into the full four-dimensional bundle $\mathcal{E}_{\mathrm{reg}}$.

When both the algebra structure constants $(\alpha(t), \beta(t))$ and the trace parameters $(\rho(t), \theta(t))$ co-evolve dynamically on $\mathcal{E}_{\mathrm{reg}}$, the scale-invariant fiber kinetic energy from (\ref{eq:h_vertical}), (\ref{eq:constrained_action}), and (\ref{eq:2d_g_matrices}) is
\begin{equation*}
    T_{\mathrm{fiber}} = \frac{1}{2}h^{ij}\dot{\tau}_i\dot{\tau}_j = \frac{1}{2\rho^2\Delta}\Big( (\Delta + u^2)\dot{\rho}^2 - 2u\dot{\rho}\dot{\theta} + \dot{\theta}^2 \Big). 
\end{equation*}
Introducing the logarithmic radial scaling coordinate $\varsigma \coloneqq \log|\rho|$ (so that $d\varsigma^2 = d\rho^2/\rho^2$ splits off orthogonally) alongside $\theta = \rho u$ ($\dot{\theta} = u\dot{\rho} + \rho\dot{u}$) diagonalizes the vertical fiber metric block to $d\varsigma^2 + du^2/\Delta$. In the original product coordinates $(\alpha, \beta, \varsigma, u)$, the four-dimensional bundle metric $\mathcal{G}$ on $\mathcal{E}_{\mathrm{reg}}$ is the orthogonal direct sum of the horizontal base metric $G_{IJ}$ from (\ref{eq:2d_G_horizontal}) and the vertical fiber metric $h^{ij}$,
\begin{equation*}
    ds^2_{\mathcal{E}} = \underbrace{\frac{d\alpha^2}{\Delta} + \frac{(u\,d\alpha + d\beta)^2}{\Delta^2}}_{\text{Horizontal base metric } G_{IJ}} + \underbrace{d\varsigma^2 + \frac{du^2}{\Delta}}_{\text{Vertical fiber metric } h^{ij}}. 
\end{equation*}

To analyze the dynamics relative to the trace horizon, we pass to the adapted bundle coordinates $(\alpha, \Delta, \varsigma, u)$ by eliminating $\beta$ via $\beta = u^2 - \alpha u - \Delta$. Because the horizon determinant $\Delta$ depends simultaneously on the base coordinates $(\alpha, \beta)$ and the fiber coordinate $u$, its exterior derivative $d\Delta = (2u - \alpha)du - (u\,d\alpha + d\beta)$ mixes horizontal and vertical variations: even if the base algebra $(\alpha, \beta)$ is held fixed ($d\alpha = d\beta = 0$), vertical motion along the trace fiber ($du \neq 0$) induces a coordinate drift $(2u - \alpha)du$ in $\Delta$. Substituting $-(u\,d\alpha + d\beta) = d\Delta - (2u - \alpha)du$ into the line element expresses $\mathcal{G}$ in a Kaluza--Klein-type warped bundle form (that is, a bundle metric coupling base and fiber coordinates via the square of a connection-shifted $1$-form),
\begin{equation}\label{eq:4d_kk_metric}
    ds^2_{\mathcal{E}} = \frac{d\alpha^2}{\Delta} + \frac{\big(d\Delta - (2u - \alpha)du\big)^2}{\Delta^2} + d\varsigma^2 + \frac{du^2}{\Delta}, 
\end{equation}
with total four-dimensional Lagrangian
\begin{equation}\label{eq:4d_kk_lagrangian}
    \mathcal{L}_{4\mathrm{D}}(\alpha, \Delta, \varsigma, u, \dot{\alpha}, \dot{\Delta}, \dot{\varsigma}, \dot{u}) = \frac{\dot{\alpha}^2}{2\Delta} + \frac{\big(\dot{\Delta} - (2u - \alpha)\dot{u}\big)^2}{2\Delta^2} + \frac{1}{2}\dot{\varsigma}^2 + \frac{\dot{u}^2}{2\Delta}. 
\end{equation}
Geometrically, (\ref{eq:4d_kk_metric}) exhibits both ingredients of this ansatz: the vertical trace-ratio fiber $du^2$ is conformally warped by the factor $\Delta^{-1}$ (which stretches fiber lengths to infinity as $\Delta \to 0$), while the $1$-form $(2u - \alpha)du$ acts as a local Ehresmann connection (or Kaluza--Klein gauge potential) twisting the horizon coordinate $\Delta$ over the $u$-fiber. Subtracting the vertical fiber drift $(2u - \alpha)\dot{u}$ from $\dot{\Delta}$ defines the horizontally projected connection velocity
\begin{equation*}
    \omega \coloneqq \dot{\Delta} - (2u - \alpha)\dot{u} = -(u\dot{\alpha} + \dot{\beta}), 
\end{equation*}
which isolates the purely horizontal motion perpendicular to the trace fiber. In terms of $\omega$, the logarithmic radial scale $\varsigma(t) = \varsigma_0 + v_\varsigma t$ decouples completely ($\ddot{\varsigma} = 0$), and the Euler--Lagrange equations derived from (\ref{eq:4d_kk_lagrangian}) take the closed, projectively invariant form
\begin{subequations}\label{eq:4d_kk_odes}
\begin{align}
    \frac{d}{dt}\left(\frac{\dot{\alpha}}{\Delta}\right) &= \frac{\omega\dot{u}}{\Delta^2},  \\
    \frac{d}{dt}\left(\frac{\omega}{\Delta^2}\right) &= -\frac{\dot{\alpha}^2 + \dot{u}^2}{2\Delta^2} - \frac{\omega^2}{\Delta^3},  \\
    \ddot{\varsigma} &= 0,  \\
    \frac{d}{dt}\left(\frac{\dot{u}}{\Delta} - \frac{(2u - \alpha)\omega}{\Delta^2}\right) &= -\frac{2\omega\dot{u}}{\Delta^2}. 
\end{align}
\end{subequations}
Notice that setting the initial fiber velocity to $\dot{u}(0) = 0$ in (\ref{eq:4d_kk_odes}a)--(\ref{eq:4d_kk_odes}b) causes the right-hand connection-coupling term $\omega\dot{u}/\Delta^2$ and the fiber kinetic contribution $\dot{u}^2/(2\Delta^2)$ to vanish identically at $t = 0$, immediately recovering the instantaneous horizontal Poincar\'e half-plane system (\ref{eq:2d_decoupled_odes}) to second order while (\ref{eq:4d_kk_odes}d) drives the $O(t^2)$ onset of trace evolution via $\ddot{u}(0) = (2u_0 - \alpha_0)\Delta(0)\frac{d}{dt}\big(\omega/\Delta^2\big)\big|_{t=0} - \frac{\dot{\alpha}(0)\omega(0)}{\Delta(0)}$.

By Theorem~\ref{thm:deformation_geodesics}, launching a geodesic ray from the dual numbers $(\alpha, \beta) = (0, 0)$ on $\mathcal{E}_{\mathrm{reg}}$ with initial horizontal velocity transverse to the coboundary orbit $\operatorname{Im}(d_{\mu_0}^{(1)}) = \operatorname{span}\{(1,0)\}$ produces a co-evolving trajectory $(\mu(t), \tau(t))$ whose base projection transits across the moduli strata: for $\Delta_{\mathrm{alg}}(t) = \alpha(t)^2 + 4\beta(t) > 0$ it enters the split semisimple stratum ($\mathbb{R} \oplus \mathbb{R}$), at $t = 0$ it passes through the dual numbers ($\Delta_{\mathrm{alg}} = 0$), and for $\Delta_{\mathrm{alg}}(t) < 0$ it enters the complex elliptic stratum ($\mathbb{C}$). Moreover, introducing the dimensionless horizon ratio
\begin{equation*}
    \eta \coloneqq \frac{\Delta_\mu(\tau)}{u^2} = 1 - \frac{\alpha}{u} - \frac{\beta}{u^2}
\end{equation*}
on the open domain $u \neq 0$ collapses the two-variable dependence of the horizontal Christoffel symbols (\ref{eq:2d_christoffel}) onto single-variable universal profiles: rescaling each component by its natural homogeneity weight in $u$ yields $u\,\Gamma^\alpha{}_{\alpha\alpha} = u^2\,\Gamma^\beta{}_{\beta\beta} = 1/\eta$, $u^2\,\Gamma^\alpha{}_{\alpha\beta} = u\,\Gamma^\beta{}_{\alpha\beta} = 1/(2\eta)$, $\Gamma^\beta{}_{\alpha\alpha} = -1/2$, and $u^3\,\Gamma^\alpha{}_{\beta\beta} = 0$. Consequently, even as $\alpha(t)$, $\beta(t)$, and $u(t)$ co-evolve dynamically in four-dimensional bundle space, the rescaled horizontal connection components depend solely on $\eta(t)$, with the four horizon-singular components undergoing universal $\eta^{-1}$ divergence as $\eta \to 0^+$, as illustrated in Figure~\ref{fig:christoffel_horizon}.

\begin{figure}[htbp]
    \centering
    \begin{tikzpicture}
        \begin{axis}[
            width=0.85\textwidth,
            height=0.52\textwidth,
            xmin=-0.03, xmax=1.0,
            x dir=reverse,
            ymin=-2, ymax=16,
            samples=200,
            domain=0.065:1.0,
            axis lines=box,
            grid=both,
            grid style={dashed, gray!30},
            xlabel={Dimensionless horizon ratio $\eta = \Delta_\mu(\tau)/u^2$ (from $\eta = 1$ at the dual numbers to $\eta = 0$)},
            ylabel={Rescaled Christoffel components},
            legend pos=north west,
            legend cell align={left},
            legend style={font=\small, fill=white, fill opacity=0.85, draw=gray!50},
            tick label style={font=\small},
            label style={font=\small}
        ]
            \addplot[thick, blue, smooth] {1 / x};
            \addlegendentry{$u\,\Gamma^\alpha{}_{\alpha\alpha} = u^2\,\Gamma^\beta{}_{\beta\beta} = \frac{1}{\eta}$}

            \addplot[thick, red, densely dashed, smooth] {1 / (2 * x)};
            \addlegendentry{$u^2\,\Gamma^\alpha{}_{\alpha\beta} = u\,\Gamma^\beta{}_{\alpha\beta} = \frac{1}{2\eta}$}

            \addplot[thick, teal!80!black, dashdotted, domain=0:1.0] {0};
            \addlegendentry{$u^3\,\Gamma^\alpha{}_{\beta\beta} = 0$}

            \addplot[thick, purple, densely dotted, domain=0:1.0] {-0.5};
            \addlegendentry{$\Gamma^\beta{}_{\alpha\alpha} = -\frac{1}{2}$}

            \draw[thick, black!70, dashed] (axis cs:0,-2) -- (axis cs:0,16);
            \node[anchor=south east, font=\footnotesize, rotate=90] at (axis cs:0.015, 4) {Trace Horizon $\eta = 0$ ($\Delta_\mu(\tau) = 0$)};
        \end{axis}
    \end{tikzpicture}
    \caption{Behavior of the homogeneity-rescaled horizontal Christoffel symbol components (\ref{eq:2d_christoffel}) as a function of the dimensionless horizon ratio $\eta = \Delta_\mu(\tau)/u^2 = 1 - \alpha/u - \beta/u^2 \in (0, 1]$, decreasing from $\eta = 1$ at the dual numbers ($\alpha = 0, \beta = 0, u \neq 0$) toward the $\eta^{-1}$ trace-horizon barrier at $\eta = 0$.}
    \label{fig:christoffel_horizon}
\end{figure}

\subsection{An example with a noncommutative algebra of infinite global dimension}
\label{sec:nakayama_geodesics}

In this section we present an example of what the methodology of Theorem~\ref{thm:deformation_geodesics} can accomplish that prior homological deformation and cyclic-transport methods cannot. 

Whereas in our two-dimensional commutative illustration (Section~\ref{sec:2d_deformations}) both $\Phi_a(\mu) = 0$ and $\mathcal{R}_{\mathrm{tr}}^b(\mu, \tau) = 0$ held identically on the ambient plane, so that neither normal constraint forces nor initial fiber velocities ($w_0 = 0$) were activated, the six-dimensional cyclic Nakayama bound quiver algebra $N_6 = \mathbb{R}Q_2 / \mathcal{I}_3$ on
\begin{equation*}
Q_2: \quad \begin{tikzcd}
1 \arrow[r,shift left=1ex,"a"] & 2 \arrow[l,shift left=1ex,"b"]
\end{tikzcd},
\end{equation*}
modulo the ideal $\mathcal{I}_3$ of paths of length $\ge 3$, activates the full noncommutative apparatus of Section~\ref{sec:def_metric}, including quadratic Gerstenhaber multipliers, trace-asymmetry tilting of $\mathcal{H}^2_{\mu_0}(\tau_0)$, first-order trace co-rotation $w_0 \neq 0$, and smooth transit from infinite global dimension into a noncommutative semisimple stratum. With left-to-right concatenation and basis $\{e_1, e_2, a, b, \omega_1, \omega_2\}$ ($\mathbf{1} = e_1 + e_2$, $\omega_1 \coloneqq ab$, $\omega_2 \coloneqq ba$), $N_6$ is non-local, noncommutative, and self-injective of Loewy length $3$. Because its simple modules have $4$-periodic projective resolutions $\cdots \to P_1 \to P_2 \to P_2 \to P_1 \to S_1 \to 0$, $N_6$ has infinite global dimension and $HH^k(N_6, N_6) \neq 0$ in infinitely many degrees, rendering its periodic cyclic bar bicomplex infinite-dimensional.

At $\mu_0 = N_6$, $[N_6, N_6] = \operatorname{span}\{a, b, \omega_1 - \omega_2\}$ has dimension $3$. Any $\tau_0 \in \mathcal{T}(N_6)$ vanishes on $[N_6, N_6]$ and is determined by $\rho_i \coloneqq \tau_0(e_i)$ ($\rho \coloneqq \rho_1 + \rho_2 = \tau_0(\mathbf{1}) \neq 0$) and $s_0 \coloneqq \tau_0(\omega_1) = \tau_0(\omega_2)$. Setting
\begin{equation*}
\delta \coloneqq \frac{\rho_1 - \rho_2}{\rho} = \frac{\tau_0(e_1 - e_2)}{\tau_0(\mathbf{1})} \in \mathbb{R}, \qquad v \coloneqq \frac{s_0}{\rho} = \frac{\tau_0(\omega_1)}{\tau_0(\mathbf{1})} \neq 0,
\end{equation*}
and passing to the adapted basis $\{\mathbf{1}, \, \tilde{h}, \, \tilde{\omega}_+, \, \omega_-, \, a, \, b\}$ with $h \coloneqq e_1 - e_2$, $\omega_\pm \coloneqq \omega_1 \pm \omega_2$, $\tilde{h} \coloneqq h - \delta\mathbf{1} \in \ker\tau_0$, and $\tilde{\omega}_+ \coloneqq \omega_+ - 2v\mathbf{1} \in \ker\tau_0$, the bilinear form $g_{\mathbf{1}}^{N_6}(\tau_0)$ from (\ref{eq:g_sym_mu}) and its inverse decouple into $4 \times 4$ and $2 \times 2$ blocks,
\begin{equation}\label{eq:N6_g}
\big[g_{ij}\big] = \begin{pmatrix}
1 & 0 & 0 & 0 \\
0 & \delta^2 - 1 & 2\delta v & -2v \\
0 & 2\delta v & 4v^2 & 0 \\
0 & -2v & 0 & 0
\end{pmatrix} \oplus \begin{pmatrix} 0 & -v \\ -v & 0 \end{pmatrix},
\end{equation}
\begin{equation}\label{eq:N6_ginv}
\big[g^{ij}\big] = \begin{pmatrix}
1 & 0 & 0 & 0 \\
0 & 0 & 0 & -\frac{1}{2v} \\[2pt]
0 & 0 & \frac{1}{4v^2} & \frac{\delta}{4v^2} \\[4pt]
0 & -\frac{1}{2v} & \frac{\delta}{4v^2} & \frac{1}{4v^2}
\end{pmatrix} \oplus \begin{pmatrix} 0 & -\frac{1}{v} \\[2pt] -\frac{1}{v} & 0 \end{pmatrix}.
\end{equation}
All $\delta^2$ terms cancel upon inversion, and $\det[g_{ij}] = (-16v^4)(-v^2) = 16v^6$ (giving $\Delta_{N_6}(\tau_0) = v^6 \neq 0$ in the quiver basis, since the basis transition matrix has determinant $4$), confirming $\tau_0 \in \mathcal{T}_1(N_6)$ for all $v \neq 0$ and $\delta \in \mathbb{R}$.

\subsubsection{Peirce Deformation Sheet, Harmonic Slice Tilting, and First-Order Trace Rotation}
Because $S = \mathbb{R}e_1 \oplus \mathbb{R}e_2 \subset N_6$ is separable and rigid, every unital associative deformation of $N_6$ is $GL_{\mathbf{1}}(V)$-equivalent to an $S$-bimodule multiplication preserving the Peirce decomposition $V = \bigoplus_{i,j=1}^2 e_i V e_j$, with $\dim HH^2(N_6, N_6) = 1$.   Restricting to the $GL_{\mathbf{1}}(V)$-transverse Peirce-preserving slice $\mathcal{M}_0 \subset \mathcal{M}$ of six-parameter multiplications $(p_1, p_2, q_1, q_2, s_1, s_2)$ given by
\begin{align*}
ab &= p_1 e_1 + q_1 \omega_1, & \omega_1^2 &= s_1 \omega_1, & \omega_1 a &= a\omega_2 = b\omega_1 = \omega_2 b = 0, \\
ba &= p_2 e_2 + q_2 \omega_2, & \omega_2^2 &= s_2 \omega_2, & &
\end{align*}
centered at $N_6$ ($q_i = 1$, $p_i = s_i = 0$), associativity on $q_1, q_2 \neq 0$ via $(ab)a = a(ba)$, $(ab)\omega_1 = a(b\omega_1)$, and $(ba)\omega_2 = b(a\omega_2)$ cuts out the curved three-dimensional Maurer--Cartan sheet $\mathcal{Z}_0 \subset \mathcal{M}_0$,
\begin{equation}\label{eq:N6_MC}
p_1 + q_1 s_1 = 0, \qquad p_2 + q_2 s_2 = 0, \qquad q_1 s_1 - q_2 s_2 = 0.
\end{equation}
Setting $\varepsilon \coloneqq -q_1 s_1 = -q_2 s_2$ (so $p_i = \varepsilon$ and $s_i = -\varepsilon/q_i$), $\varepsilon = 0$ parametrizes the two-dimensional coboundary orbit $(q_1, q_2)$ of $N_6$ via $\omega_i \mapsto q_i\omega_i$, whereas any $\varepsilon \neq 0$ yields orthogonal central idempotents $z_i = s_i^{-1}\omega_i$ splitting $B_\mu \cong \mathbb{R}z_1 \oplus \mathbb{R}z_2 \oplus M_2(\mathbb{R})$. Crucially, $\dim[B_\mu, B_\mu] \equiv 3$ remains constant across $\varepsilon = 0$, since $[\mathbb{R}^2 \oplus M_2(\mathbb{R}), \, \mathbb{R}^2 \oplus M_2(\mathbb{R})] = \mathfrak{sl}_2(\mathbb{R})$.

For any $B_\mu \in \mathcal{Z}_0$, $[e_1, a]_\mu = a$ and $[b, e_1]_\mu = b$ force $\tau(a) = \tau(b) = 0$, while $[a, b]_\mu = \varepsilon(e_1 - e_2) + q_1\omega_1 - q_2\omega_2$ imposes the trace-cyclicity constraint
\begin{equation}\label{eq:N6_trace_constraint}
\mathcal{R}_{\mathrm{tr}}(\mu, \tau) \coloneqq \varepsilon\big(\tau(e_1) - \tau(e_2)\big) + q_1\tau(\omega_1) - q_2\tau(\omega_2) = 0,
\end{equation}
defining the six-dimensional constrained Peirce subbundle $\mathcal{E}_{\mathrm{reg}} \coloneqq \big\{ (\mu, \tau) \in \mathcal{Z}_0 \times V^* : \tau \in \mathcal{T}_1(B_\mu) \big\}$. On $T_{N_6}\mathcal{Z}_0$ (parametrized by $\dot{q}_\pm \coloneqq \frac{1}{2}(\dot{q}_1 \pm \dot{q}_2)$ and $\dot{\varepsilon}$, with $\dot{p}_i = \dot{\varepsilon}$ and $\dot{s}_i = -\dot{\varepsilon}$), contracting the $\{a, b\}$-variations $\dot{\mu}(a, b), \dot{\mu}(b, a)$ via $(g^{ab})^2 = v^{-2}$ and the $\{\tilde{\omega}_+, \omega_-\}$-variations $\dot{\mu}(\omega_i, \omega_i) = -\dot{\varepsilon}\omega_i$ via (\ref{eq:N6_g})--(\ref{eq:N6_ginv}) cancels all higher powers of $\delta$, yielding
\begin{equation}\label{eq:N6_G0_induced}
\bar{G}_{N_6}\big((\dot{q}_+, \dot{q}_-, \dot{\varepsilon}); \delta, v\big) = \frac{1}{v^2}\big(\dot{\varepsilon} + 2v\dot{q}_+\big)^2 - \frac{1}{v^2}\big(\delta\dot{\varepsilon} + 2v\dot{q}_-\big)^2 + \frac{1 + \delta^2}{v^2}\dot{\varepsilon}^2.
\end{equation}
On the coboundary plane $\dot{\varepsilon} = 0$, (\ref{eq:N6_G0_induced}) restricts to $4(\dot{q}_+^2 - \dot{q}_-^2) = 4\dot{q}_1\dot{q}_2$ ($\det = -16 \neq 0$), and its $\bar{G}_{N_6}$-orthogonal complement is spanned by the harmonic $2$-cocycle ($\dot{\varepsilon}(0) = 1$),
\begin{equation}\label{eq:N6_vharm}
v_{\mathrm{harm}}(\delta, v) = \big(\dot{q}_1(0), \, \dot{q}_2(0), \, \dot{\varepsilon}(0)\big) = \left( -\frac{1 + \delta}{2v}, \; -\frac{1 - \delta}{2v}, \; 1 \right),
\end{equation}
demonstrating Mechanism~(a): varying the trace asymmetry $\delta$ tilts $v_{\mathrm{harm}}(\delta, v)$ linearly across the coboundary fibers.

Differentiating (\ref{eq:N6_trace_constraint}) at $t = 0$ gives the linearized compatibility law
\begin{equation*}
w_0(\omega_1 - \omega_2) = -\rho\Big[ \delta\,\dot{\varepsilon}(0) + v\big(\dot{q}_1(0) - \dot{q}_2(0)\big) \Big].
\end{equation*}
Although $\tilde{h}^*$ is null for $h(N_6; \tau_0)\big|_{\mathcal{T}(N_6)}$ ($g^{\tilde{h}\tilde{h}} = g^{\tilde{h}\tilde{\omega}_+} = 0$ in (\ref{eq:N6_ginv})), its off-diagonal pairing $g^{\tilde{h}\omega_-} = -\frac{1}{2v}$ against $w_0(\omega_-) = w_0(\omega_1 - \omega_2)$ renders $\bar{\mathcal{G}} = \mathcal{G}|_{T\mathcal{E}_{\mathrm{reg}}}$ non-degenerate across all six tangent directions at $(N_6, \tau_0)$, satisfying Theorem~\ref{thm:deformation_geodesics}. Moreover, when $\delta \neq 0$, any un-tilted cocycle with $\dot{q}_1(0) = \dot{q}_2(0)$---such as $(0, 0, 1)$---forces $w_0(\omega_1 - \omega_2) = -\rho\delta \neq 0$ at $t = 0$ (illustrating Mechanism~(b) in contrast to $w_0 = 0$ in Section~\ref{sec:2d_deformations}), whereas along $v_{\mathrm{harm}}(\delta, v)$ in (\ref{eq:N6_vharm}) the tilt $\dot{q}_1(0) - \dot{q}_2(0) = -\delta/v$ cancels $\delta\dot{\varepsilon}(0)$ so that $(v_{\mathrm{harm}}, 0)$ is both $\bar{G}_{N_6}$- and $\bar{\mathcal{G}}$-orthogonal to the coboundary fibers, with $\ddot{\tau}(0) \neq 0$ entering at second order.

\subsubsection{Exact Constrained Geodesic Equations and Active Gerstenhaber Multipliers}
At $\delta = 0$ ($\tau_0(e_i) = \frac{1}{2}\rho$), the involution $\iota \in GL_{\mathbf{1}}(V)$ swapping $e_1 \leftrightarrow e_2$, $a \leftrightarrow b$, and $\omega_1 \leftrightarrow \omega_2$ is an isometry of $\mathcal{G}$. By Palais's principle of symmetric criticality, the $\iota$-fixed locus ($p_i = p$, $q_i = q$, $s_i = s$, $\tau(e_i) = \frac{1}{2}\rho$, $\tau(\omega_i) = \rho v$) is totally geodesic in $\mathcal{E}_{\mathrm{reg}}$ and satisfies (\ref{eq:N6_trace_constraint}) identically, leaving the single quadratic Gerstenhaber constraint
\begin{equation}\label{eq:N6_sym_constraint}
\bar{\Phi}(p, q, s) \coloneqq p + qs = 0.
\end{equation}
On this $\iota$-symmetric subspace (prior to imposing (\ref{eq:N6_sym_constraint})), evaluating (\ref{eq:g_sym_mu}) in $\{\mathbf{1}, \, h, \, \tilde{\omega}_+, \, \omega_-, \, a, \, b\}$ yields
\begin{equation}\label{eq:N6_g_deformed}
g_{\mathbf{1}}^{B_\mu}(\tau) = \begin{pmatrix}
1 & 0 & 0 & 0 \\
0 & -1 & 0 & -2v \\
0 & 0 & 2\Omega & 0 \\
0 & -2v & 0 & -2sv
\end{pmatrix} \oplus \begin{pmatrix} 0 & -\Pi \\ -\Pi & 0 \end{pmatrix},
\end{equation}
where $\Omega \coloneqq v(2v - s)$ and $\Pi \coloneqq \frac{1}{2}p + qv$ ($\det g_{\mathbf{1}}^{B_\mu}(\tau) = 4\Omega^2\Pi^2$). Using $g^{\tilde{\omega}_+\tilde{\omega}_+} = g^{\omega_-\omega_-} = \frac{1}{2\Omega}$ and $g^{ab} = -\frac{1}{\Pi}$ to contract the $\{a, b\}$-variations $\dot{\mu}(a, b), \dot{\mu}(b, a)$, the $\{\tilde{\omega}_+, \omega_-\}$-variations $\dot{\mu}(\omega_i, \omega_i) = \dot{s}\omega_i$, and $\dot{\tau} = \dot{\rho}\mathbf{1}^* + 2\rho\dot{v}\tilde{\omega}_+^*$ via (\ref{eq:G_horizontal})--(\ref{eq:h_vertical}) cancels all intermediate $sv$ terms, giving the constrained action integrand
\begin{equation}\label{eq:N6_L_constrained}
\mathcal{L}_{N_6}(p, q, s, \varsigma, v; \lambda_{\mathrm{G}}) = \frac{(\dot{p} + 2v\dot{q})^2}{2\Pi^2} + \frac{\dot{s}^2}{\Omega} + \frac{1}{2}\dot{\varsigma}^2 + \frac{\dot{v}^2}{\Omega} + \lambda_{\mathrm{G}}(t)(p + qs),
\end{equation}
which is Riemannian on $\Omega > 0, \Pi \neq 0$. Varying $p$ determines the Gerstenhaber constraint force along the constrained flow (for $P_\xi$ below) as
\begin{equation*}
\lambda_{\mathrm{G}}(t) = \frac{d}{dt}\left(\frac{\dot{p} + 2v\dot{q}}{\Pi^2}\right) + \frac{(\dot{p} + 2v\dot{q})^2}{2\Pi^3} = -\frac{\dot{v} P_\xi}{\Pi(2v - s)},
\end{equation*}
which is actively non-zero whenever $P_\xi \neq 0$ and $\dot{v} \neq 0$, and vanishes along the harmonic ray $P_\xi = 0$.

Eliminating $p = -qs$ (so $\Pi = \frac{1}{2}q(2v - s)$ and $\frac{\dot{p} + 2v\dot{q}}{\Pi} = \frac{2\dot{q}}{q} - \frac{2\dot{s}}{2v - s}$) with $\xi \coloneqq \log|q|$ reduces (\ref{eq:N6_L_constrained}) to the unconstrained four-dimensional geodesic Lagrangian
\begin{equation}\label{eq:N6_L_reduced}
\bar{\mathcal{L}}_{N_6}(\xi, s, \varsigma, v, \dot{\xi}, \dot{s}, \dot{\varsigma}, \dot{v}) = 2\left(\dot{\xi} - \frac{\dot{s}}{2v - s}\right)^2 + \frac{\dot{s}^2}{v(2v - s)} + \frac{1}{2}\dot{\varsigma}^2 + \frac{\dot{v}^2}{v(2v - s)},
\end{equation}
with cyclic conservation laws $\ddot{\varsigma} = 0$ and $P_\xi \coloneqq \partial_{\dot{\xi}}\bar{\mathcal{L}}_{N_6} = 4\big(\dot{\xi} - \frac{\dot{s}}{2v - s}\big) = \text{const}$. Launching from $(N_6, \tau_0)$ ($\xi(0) = s(0) = \dot{v}(0) = 0$, $v(0) = v_* > 0$) along $v_{\mathrm{harm}}(0, v_*)$ in (\ref{eq:N6_vharm}) ($\dot{s}(0) = -1$, $\dot{\xi}(0) = -\frac{1}{2v_*}$) sets $P_\xi = 0$, locking $\dot{\xi} = \frac{\dot{s}}{2v - s}$ for all $t$ and reducing the $(s, v)$ co-evolution to the conformally flat Riemannian Lagrangian
\begin{equation}\label{eq:N6_2D_Routhian}
\mathcal{L}_{\mathrm{eff}}(s, v, \dot{s}, \dot{v}) = \frac{\dot{s}^2 + \dot{v}^2}{v(2v - s)}
\end{equation}
on $\Omega = v(2v - s) > 0$ (reconstructing $q(t) = e^{\xi(t)}$ and $\varepsilon(t) = -q(t)s(t)$ by quadrature), smoothly carrying the Maurer--Cartan path from $N_6$ ($s = 0$, infinite global dimension) into $\mathbb{R} \oplus \mathbb{R} \oplus M_2(\mathbb{R})$ ($s \neq 0$) up to the trace horizons $v = 0$ and $s = 2v$.

This explicit reduction brings into sharp relief what classical homological methods cannot accomplish on $N_6$. First, while one can write down \emph{ad hoc} Maurer--Cartan paths on $\mathcal{Z}_0$ by inspection---most naively fixing $q_1(t) = q_2(t) \equiv 1$ with $p_i(t) = t$ and $s_i(t) = -t$ in~(\ref{eq:N6_MC})---any smooth positive profile $q_i(t)$ with $q_i(0) = 1$ and $s_i(t) = -t/q_i(t)$ represents the same class in $HH^2(N_6, N_6)$, and classical deformation theory provides no canonical criterion to select among these coboundary gauges; by contrast, (\ref{eq:N6_vharm}) shows that the naive path $\dot{q}_i(0) = 0$ is never coboundary-orthogonal, and $P_\xi = 0$ in~(\ref{eq:N6_L_reduced}) uniquely locks the nonlinear gauge evolution $\dot{q}/q = \dot{s}/(2v - s)$ for all $t$. Second, because $N_6$ has infinite global dimension, co-evolving cyclic traces classically via Getzler's Gauss--Manin connection requires working on the non-terminating, infinite-dimensional $(b, B)$ periodic cyclic bar bicomplex, whereas~(\ref{eq:N6_2D_Routhian}) co-evolves the concrete trace ratio $v(t)$ via a two-degree-of-freedom Riemannian ODE. Third, homological trace transport tracks only commutator vanishing and is completely blind to Frobenius non-degeneracy, whereas the conformal denominator $\Omega = v(2v - s)$ in~(\ref{eq:N6_2D_Routhian}) explicitly detects the exact trace horizons $v = 0$ and $s = 2v$ at which $\det g_{\mathbf{1}}^{B_\mu}(\tau)$ degenerates.

\subsection{An example featuring nontrivial geodesic spray and global singularity invariants}
\label{sec:A3_spray}

In both Section~\ref{sec:2d_deformations} and Section~\ref{sec:nakayama_geodesics}, the non-rigid launch algebra satisfied $\dim HH^2 = 1$, confining the harmonic deformation fan to a single ray line. To demonstrate a genuine geodesic spray ($d = 2$) that steers forward trajectories across distinct moduli strata and exhibits both local and global catastrophe phenomena, consider the three-dimensional indecomposable  symmetric Frobenius algebra $A_3 = \mathbb{R}[x]/(x^3)$ with ordered linear basis $\{e_1, e_2, e_3\} = \{\mathbf{1}, x, x^2\}$ ($e_2^2 = e_3$, $e_2 e_3 = e_3^2 = 0$), for which $\dim HH^2(A_3, A_3) = 2$.

Every linear functional $\tau = (\rho, \theta_1, \theta_2) \in V^*$ is a trace on commutative deformations of $A_3$, parametrized on $\rho = \tau(\mathbf{1}) \neq 0$ by $\varsigma \coloneqq \log|\rho|$ and the normalized trace ratios $u_1 \coloneqq \theta_1/\rho$, $u_2 \coloneqq \theta_2/\rho$. At the basepoint $\mu_0 = A_3$, passing to the trace-centered basis $\{\mathbf{1}, \, \tilde{e}_2, \, \tilde{e}_3\}$ with $\tilde{e}_2 \coloneqq e_2 - u_1\mathbf{1}$ and $\tilde{e}_3 \coloneqq e_3 - u_2\mathbf{1}$ evaluates $g_{\mathbf{1}}^{A_3}(\tau_0)$ in (\ref{eq:g_sym_mu}) and its inverse to
\begin{equation}\label{eq:A3_g0}
\big[g_{ij}(\mu_0; \tau_0)\big] = \begin{pmatrix}
1 & 0 & 0 \\
0 & u_1^2 - u_2 & u_1 u_2 \\
0 & u_1 u_2 & u_2^2
\end{pmatrix},  \qquad
\big[g^{ij}(\mu_0; \tau_0)\big] = \begin{pmatrix}
1 & 0 & 0 \\
0 & -\frac{1}{u_2} & \frac{u_1}{u_2^2} \\[6pt]
0 & \frac{u_1}{u_2^2} & \frac{u_2 - u_1^2}{u_2^3}
\end{pmatrix}, 
\end{equation}
with horizon determinant $\Delta_{A_3}(\tau_0) = -u_2^3 \neq 0$ for all $u_2 \neq 0$ and $u_1 \in \mathbb{R}$.

Because every unital associative deformation near $A_3$ remains generated by $e_2$ (and is hence automatically commutative), fixing the natural $GL_{\mathbf{1}}(V)$-transverse cyclic-generator slice $e_2 \cdot_\mu e_2 = e_3$ cuts out a smooth, unobstructed three-dimensional Maurer--Cartan slice $\mathcal{Z}_0 \subset \mathcal{Z}$ of monic cubic quotients $\mathbb{R}[x]/(x^3 - c_2 x^2 - c_1 x - c_0)$ defined for $(c_0, c_1, c_2) \in \mathbb{R}^3$ by
\begin{equation}\label{eq:A3_Z0}
e_2^2 = e_3,  \qquad e_2 e_3 = e_3 e_2 = c_0\mathbf{1} + c_1 e_2 + c_2 e_3,  \qquad e_3^2 = c_0 c_2\mathbf{1} + (c_0 + c_1 c_2)e_2 + (c_1 + c_2^2)e_3, 
\end{equation}
over which $\mathcal{E}_{\mathrm{reg}} \coloneqq \{(\mu, \tau) \in \mathcal{Z}_0 \times V^* : \rho \neq 0, \, \Delta_\mu(\tau) \neq 0\}$ is a real-analytic six-dimensional constrained subbundle endowed with the induced pullback metric $\bar{\mathcal{G}} = \mathcal{G}|_{T\mathcal{E}_{\mathrm{reg}}}$. 
Within $\mathcal{Z}_0$, residual unital generator shifts $\bar{e}_2 = e_2 + \epsilon\mathbf{1}$ accompanied by $\bar{e}_3 = \bar{e}_2^2 = e_3 + 2\epsilon e_2 + \epsilon^2\mathbf{1}$ preserve the slice condition $\bar{e}_2^2 = \bar{e}_3$ and send $(c_0, c_1, c_2) \mapsto (\epsilon^3, -3\epsilon^2, 3\epsilon)$ at $\mu_0$ (since $\bar{e}_2^3 = (e_2 + \epsilon\mathbf{1})^3 = \epsilon^3\mathbf{1} - 3\epsilon^2\bar{e}_2 + 3\epsilon\bar{e}_2^2$ when $e_2^3 = 0$), spanning the one-dimensional coboundary tangent line $\operatorname{Im}(d_{\mu_0}^{(1)}) \cap T_{\mu_0}\mathcal{Z}_0 = \mathbb{R}\partial_{c_2}$. At the reference calibration $u_1 = 0, u_2 \neq 0$, we have $[g_{ij}] = \begin{psmallmatrix} 1 & 0 & u_2 \\ 0 & -u_2 & 0 \\ u_2 & 0 & 2u_2^2 \end{psmallmatrix}$ and $[g^{ij}] = \begin{psmallmatrix} 2 & 0 & -u_2^{-1} \\ 0 & -u_2^{-1} & 0 \\ -u_2^{-1} & 0 & u_2^{-2} \end{psmallmatrix}$ in $\{e_1, e_2, e_3\}$. Because $g^{12} = g^{23} = 0$ and $\dot{\mu}(e_1, \cdot) = \dot{\mu}(e_2, e_2) = 0$ along $\mathcal{Z}_0$, the contraction (\ref{eq:G_horizontal}) decouples cleanly into $2 g^{22}g^{33}\big\|\dot{\mu}(e_2, e_3)\big\|_g^2 + (g^{33})^2\big\|\dot{\mu}(e_3, e_3)\big\|_g^2$ with $\dot{\mu}(e_2, e_3) = \dot{c}_0\mathbf{1} + \dot{c}_1 e_2 + \dot{c}_2 e_3$ and $\dot{\mu}(e_3, e_3) = \dot{c}_0 e_2 + \dot{c}_1 e_3$, yielding the exact induced horizontal quadratic form on $T_{\mu_0}\mathcal{Z}_0$,
\begin{equation}\label{eq:A3_G0_induced}
\bar{G}_{A_3}\big((\dot{c}_0, \dot{c}_1, \dot{c}_2); \, 0, u_2\big) = -\frac{2}{u_2^3}\dot{c}_0^2 + \frac{4}{u_2^2}\dot{c}_1^2 - \frac{1}{u_2}\left(\frac{\dot{c}_0}{u_2} + 2\dot{c}_2\right)^2.
\end{equation}
On the coboundary line $v_{\mathrm{cob}} = (0, 0, 1)$, $\bar{G}_{A_3}(v_{\mathrm{cob}}, v_{\mathrm{cob}}) = -\frac{4}{u_2} \neq 0$ is non-degenerate, and $\bar{G}_{A_3}(v_{\mathrm{cob}}, v) = -\frac{2}{u_2}\big(\frac{\dot{c}_0}{u_2} + 2\dot{c}_2\big) = 0$ cuts out the two-dimensional harmonic slice $\mathcal{H}^2_{\mu_0}(\tau_0) \cong HH^2(A_3, A_3)$ within $T_{\mu_0}\mathcal{Z}_0$ defined by $\dot{c}_2 = -\frac{1}{2u_2}\dot{c}_0$, with orthogonal harmonic basis
\begin{equation}\label{eq:A3_vharm}
v_{\mathrm{harm}}^{(0)} = \left(1, \, 0, \, -\frac{1}{2u_2}\right),  \qquad v_{\mathrm{harm}}^{(1)} = (0, \, 1, \, 0), 
\end{equation}
and restricted metric $\bar{G}_{A_3}\big|_{\mathcal{H}^2} = -\frac{2}{u_2^3}\dot{c}_0^2 + \frac{4}{u_2^2}\dot{c}_1^2$ (strictly positive-definite Riemannian on the chamber $u_2 < 0$, and Lorentzian on $u_2 > 0$). Perturbing $u_1 \neq 0$ mixes $\dot{c}_1$ into (\ref{eq:A3_G0_induced}) via $g^{23} \neq 0$ in (\ref{eq:A3_g0}), tilting $\mathcal{H}^2_{\mu_0}(\tau_0)$ across the coboundary fibers.

Because every algebra in $\mathcal{Z}_0$ is commutative, $\mathcal{R}_{\mathrm{tr}} \equiv 0$, so launching a directional harmonic spray of amplitude $R > 0$ across angles $\vartheta \in [0, 2\pi)$ from $(\mu_0, \tau_0)$ ($u_1(0) = 0$, $u_2(0) = u_{2,0} \neq 0$) initializes at zero fiber velocity,
\begin{equation}\label{eq:A3_spray_init}
v_0(\vartheta) = R\cos\vartheta\,v_{\mathrm{harm}}^{(0)} + R\sin\vartheta\,v_{\mathrm{harm}}^{(1)} = \left(R\cos\vartheta, \; R\sin\vartheta, \; -\frac{R\cos\vartheta}{2u_{2,0}}\right), \quad w_0 = \dot{\tau}(0) = 0,
\end{equation}
sweeping out the $d = 2$ deformation fan $(t, \vartheta) \mapsto \exp^\perp_{(A_3, \tau_0)}\big(t v_0(\vartheta)\big)$ of Section~\ref{sec:def_metric} (tracing an ellipse in $\mathcal{H}^2_{\mu_0}(\tau_0)$ on the Riemannian chamber $u_{2,0} < 0$, and crossing the Lorentzian null lines on $u_{2,0} > 0$). The $\mathbb{Z}_2$-grading involution $\sigma \in GL_{\mathbf{1}}(V)$ given by $\sigma(\mathbf{1}) = \mathbf{1}$, $\sigma(e_2) = -e_2$, $\sigma(e_3) = e_3$ acts on $\mathcal{E}_{\mathrm{reg}}$ as $(c_0, c_1, c_2, \varsigma, u_1, u_2) \mapsto (-c_0, c_1, -c_2, \varsigma, -u_1, u_2)$ and is an isometry of $\bar{\mathcal{G}}$, rendering its fixed locus $\{c_0 = c_2 = u_1 = 0\}$ totally geodesic by Palais's principle of symmetric criticality. Along this locus ($B_{c_1} = \mathbb{R}[x]/(x^3 - c_1 x)$), $g_{\mathbf{1}}^{B_{c_1}}(\tau)$ block-diagonalizes into $g_{22} = -u_2$ and $\begin{psmallmatrix} 1 & u_2 \\ u_2 & 2u_2^2 - c_1 u_2 \end{psmallmatrix}$ with horizon determinant $\Delta_{c_1}(\tau) = -u_2^2(u_2 - c_1)$ and inverse entries $g^{22} = -u_2^{-1}$, $g^{33} = u_2^{-1}(u_2 - c_1)^{-1}$. Summing $2 g^{22}g^{33}g_{22}\dot{c}_1^2 + (g^{33})^2 g_{33}\dot{c}_1^2$ for the non-zero velocities $\dot{\mu}(e_2, e_3) = \dot{\mu}(e_3, e_2) = \dot{c}_1 e_2$ and $\dot{\mu}(e_3, e_3) = \dot{c}_1 e_3$ alongside (\ref{eq:h_vertical}) (where $\tau = \rho(\mathbf{1}^* + u_2 e_3^*)$ cancels all $(\dot{\varsigma}, \dot{u}_2)$ cross-terms) yields the exact global reduced Lagrangian governing the principal rays $\vartheta = \pm\frac{\pi}{2}$,
\begin{equation}\label{eq:A3_c1_lagrangian}
\mathcal{L}_{(1)}(c_1, \varsigma, u_2, \dot{c}_1, \dot{\varsigma}, \dot{u}_2) = \frac{4u_2 - 3c_1}{2u_2(u_2 - c_1)^2}\dot{c}_1^2 + \frac{1}{2}\dot{\varsigma}^2 + \frac{\dot{u}_2^2}{2u_2(u_2 - c_1)}.
\end{equation}
Thus, for $t > 0$ up to the trace horizons $u_2 = 0$ or $c_1 = u_2$ (and internal metric wall $3c_1 = 4u_2$), the forward ray $\vartheta = \frac{\pi}{2}$ ($\dot{c}_1(0) = R > 0$) remains strictly inside the completely split semisimple stratum $\mathbb{R} \oplus \mathbb{R} \oplus \mathbb{R}$ ($c_1(t) > 0$), whereas the antipodal ray $\vartheta = -\frac{\pi}{2}$ ($\dot{c}_1(0) = -R < 0$) remains strictly inside $\mathbb{R} \oplus \mathbb{C}$ ($c_1(t) < 0$).

Across general launch angles $\vartheta$ in (\ref{eq:A3_spray_init}), the cubic discriminant $\Delta_{\mathrm{alg}}(c_0, c_1, c_2) = c_1^2 c_2^2 + 4c_1^3 - 4c_0 c_2^3 - 27c_0^2 - 18c_0 c_1 c_2$ along $(\mu(t), \tau(t))$ expands near $t = 0$ as
\begin{multline*}
\Delta_{\mathrm{alg}}\big(\mu(t)\big) = -27 R^2 \cos^2\vartheta \, t^2 \\ + \left(4 R^3 \sin^3\vartheta + \frac{9 R^3}{u_{2,0}}\cos^2\vartheta\sin\vartheta - 27 R\cos\vartheta\,\ddot{c}_0(0)\right)t^3 + O(t^4). 
\end{multline*}
Because $\sigma$-symmetry forces $\ddot{c}_0(0)$ to vanish whenever $\cos\vartheta = 0$, every non-vertical ray ($\cos\vartheta \neq 0$) initially enters $\mathbb{R} \oplus \mathbb{C}$ ($\Delta_{\mathrm{alg}} < 0$) while activating $\sigma$-odd trace evolution $\ddot{u}_1(0) \neq 0$. For shallow angles $0 < |\frac{\pi}{2} - \vartheta| \ll 1$, the cubic term $4R^3\sin^3\vartheta\,t^3$ overtakes $-27R^2\cos^2\vartheta\,t^2$ at a finite transit time $t_{\mathrm{trans}}(\vartheta) = \frac{27}{4R}\cos^2\vartheta + O(\cos^4\vartheta) > 0$, driving the forward geodesic across the discriminant wall $\Delta_{\mathrm{alg}} = 0$ (passing instantaneously through $\mathbb{R} \oplus \mathbb{R}[\epsilon]/(\epsilon^2)$) into $\mathbb{R} \oplus \mathbb{R} \oplus \mathbb{R}$.

\subsubsection{Local and Global Catastrophes, Maslov Index, and Monodromy}
This $d = 2$ deformation spray exhibits an instructive interplay between internal (Section~3) and bundle/moduli (Section~4) singularity invariants:

\begin{enumerate}
    \item \emph{Evolution of Local ART Catastrophes:} At $t = 0$ ($\mu_0 = A_3$) with $u_1 = 0$ and $u_2 = -w < 0$, evaluating the Section~\ref{sec:catastrophes} ambient Legendre transform on $T^*A_3$ along the $g_{\mathbf{1}}(\tau_0)$-orthonormalized radical channel $V = \mathbb{R}\mathbf{e}_1 \subset \operatorname{rad}(A_3) \cap \ker\tau_0$ ($\mathbf{e}_1 = w^{-1/2}e_2$, $W(\tau_0) = \mathbb{R}(e_3 + w\mathbf{1})$) for $\mathbf{p} = (1 - h)\mathbf{1}^*_{\tau_0} + \sigma\mathbf{e}_1^*$ imposes $y_0 = (1 - h)^{-1}$ and yields the exact dispersion equation $\log(1 - h) + \frac{\sigma^2}{2(1 - h)^2} = 0$. Solving for $h = H_\perp(\sigma; \tau_0)$ in powers of $\sigma$ gives $H_\perp(\sigma; \tau_0) = \frac{1}{2}\sigma^2 + \frac{3}{8}\sigma^4 + O(\sigma^6)$, so under $S_0(\sigma) = \frac{1}{2}k_0\sigma^2$ ($k_0 \neq 0$) the focal phase germ at $z^* = k_0$ is $f_0(\sigma) = -\frac{3}{8}k_0\sigma^4 + O(\sigma^6)$. Thus the algebra $A_3 = \mathbb{R}[x]/(x^3)$ realizes the isolated Arnol'd cusp catastrophe $A_3^\pm$ ($\mu_{\mathrm{Milnor}} = 3$; standard cusp $A_3^+$ for $k_0 < 0$, dual cusp $A_3^-$ for $k_0 > 0$), across which the transverse Hessian $\partial_\sigma^2\Phi|_{\sigma=0} = k_0 - z$ flips sign and induces an axial Maslov (Morse) focal index jump of $1$. Along forward Maurer--Cartan rays ($t > 0$), entering the semisimple strata $\Delta_{\mathrm{alg}} \neq 0$ annihilates $\operatorname{rad}(B_{\mu(t)})$ (so $P_{\boldsymbol{\tau}} \equiv 0$ by Proposition~\ref{prop:WM_factorization}), whereas $\sigma$-breaking rays ($\cos\vartheta \neq 0$) generate $u_1(t) \neq 0$ and $c_0(t) \neq 0$ that unfold the $A_3^\pm$ cusp on $\ker\tau(t)$ into Arnol'd fold ($A_2$, $\mu_{\mathrm{Milnor}} = 2$) caustics; at $t = t_{\mathrm{trans}}(\vartheta)$, crossing $\Delta_{\mathrm{alg}} = 0$ ($B_\mu \cong \mathbb{R} \oplus \mathbb{R}[\epsilon]/(\epsilon^2)$) momentarily restores a one-dimensional radical whose transverse slices in $\ker\tau$ likewise exhibit generic $A_2$ folds.
    
    \item \emph{Global Bundle Curvature, Metric-Wall Cusps, and Horizon/Moduli Monodromy:} Under the homothetic $GL_{\mathbf{1}}(V)$-scaling $e_2 \mapsto \lambda^{1/2}e_2$, $e_3 \mapsto \lambda e_3$ ($\lambda > 0$), introducing $\chi \coloneqq \log|u_2|$ and the dimensionless ratio $q \coloneqq c_1 / u_2$ renders $\chi$ cyclic in (\ref{eq:A3_c1_lagrangian}), reducing the $(c_1, u_2)$ surface metric and its exact Gaussian curvature to
    \begin{align}\nonumber
    ds^2_{(1)} &= \frac{4 - 3q}{(1 - q)^2}(dq + q\,d\chi)^2 + \frac{1}{1 - q}d\chi^2,  \\  K_{(1)}(q) &= -\frac{18q^4 - 81q^3 + 147q^2 - 149q + 71}{4(4 - 3q)^2(1 - q)}.  \label{eq:A3_K}
    \end{align}
    On the entire Riemannian chamber $u_2 < 0, \, q < 1$ (containing both rays $\vartheta = \pm\frac{\pi}{2}$ launched from $q = 0$, where $K_{(1)}(0) = -\frac{71}{64}$), the quartic numerator of (\ref{eq:A3_K}) is strictly decreasing and bounded below by $6$, so $K_{(1)}(q) < 0$ everywhere and diverges as $K_{(1)}(q) \sim -\frac{3}{2(1 - q)} \to -\infty$ at $q \to 1^-$: $(c_1, u_2)$-variations are strictly defocusing with zero Morse/Maslov conjugate index all the way to the infinite-distance trace horizon $q = 1$ ($c_1 = u_2$). Conversely, inside the indefinite regime $1 < q < \frac{4}{3}$, the locus $q = \frac{4}{3}$ ($3c_1 = 4u_2$) forms an internal metric-degeneracy wall where $\Delta_{c_1}(\tau) = \frac{1}{3}u_2^3 \neq 0$ remains strictly inside $\mathcal{T}_1(B_{c_1})$ while $K_{(1)}(q) \sim -\frac{13}{12(4 - 3q)^2} \to -\infty$. Geodesics of the reduced surface Lagrangian $L_{\mathrm{surf}} \coloneqq \frac{1}{2}ds^2_{(1)}(\dot{\gamma}, \dot{\gamma})$ reaching $q = \frac{4}{3}$ with conserved momentum $P_\chi \coloneqq \partial_{\dot{\chi}} L_{\mathrm{surf}} \neq 0$ and conserved energy $E \coloneqq L_{\mathrm{surf}}$ satisfying $2E + \frac{1}{3}P_\chi^2 > 0$ form finite-time Kossowski semicubical cusps ($A_2$ singularities of transverse signature-changing metrics~\cite{Kossowski1987}), $(q - \frac{4}{3})^3 \propto -(\chi - \chi_*)^2$, concretely realizing the metric-degeneracy termination in Theorem~\ref{thm:deformation_geodesics}(i). Simultaneously, crossing the trace horizon $q = 1$ flips the sign of $1 - q$ in $ds^2_{(1)}$ (inducing a $\mathbb{Z}_2$ metric/Maslov index jump of $1$), while encircling $q = 1$ in the complexified $q$-plane induces a $(-1)$ monodromy on the volume half-density $\sqrt{\det g_{(1)}} \propto (4 - 3q)^{1/2}(1 - q)^{-3/2}$. Finally, projecting $\mathcal{Z}_0$ modulo shifts $\bar{e}_2 = e_2 + \epsilon\mathbf{1}$ onto the depressed-cubic (Tschirnhaus) moduli plane $(\tilde{c}_0, \tilde{c}_1) = \big(c_0 + \frac{1}{3}c_1 c_2 + \frac{2}{27}c_2^3, \, c_1 + \frac{1}{3}c_2^2\big)$ maps the $d = 2$ spray directly onto Arnol'd's semicubical cusp discriminant $\Delta_{\mathrm{alg}} = 4\tilde{c}_1^3 - 27\tilde{c}_0^2 = 0$, endowing complexified loops in $\mathbb{C}^2 \setminus \{\Delta_{\mathrm{alg}} = 0\}$ (and around the complexified discriminant cone in $\mathcal{H}^2_{\mu_0}(\tau_0) \otimes \mathbb{C}$) with non-abelian $S_3$ idempotent monodromy (Figure~\ref{fig:A3_spray_catastrophes}).
\end{enumerate}

\begin{figure}[htbp]
    \centering
    \begin{tikzpicture}
        \begin{axis}[
            name=panelA,
            width=0.94\textwidth,
            height=0.43\textwidth,
            xmin=-1.35, xmax=1.35,
            ymin=-0.68, ymax=1.92,
            clip mode=individual,
            axis lines=box,
            axis background/.style={fill=blue!3},
            xlabel={Tschirnhaus coordinate $\tilde{c}_0 = c_0 + \frac{1}{3}c_1 c_2 + \frac{2}{27}c_2^3$},
            ylabel={$\tilde{c}_1 = c_1 + \frac{1}{3}c_2^2$},
            title={\textbf{(a)} $d=2$ harmonic geodesic spray and wall-crossing on the Tschirnhaus moduli plane},
            title style={font=\small},
            tick label style={font=\footnotesize},
            label style={font=\small}
        ]
            \addplot[domain=-1.32:1.32, samples=120, smooth, fill=teal!12, draw=none]
                ({0.38490018 * x^3}, {x^2}) \closedcycle;

            \draw[gray!50, dotted] (axis cs:-1.35,0) -- (axis cs:1.35,0);
            \draw[gray!50, dotted] (axis cs:0,-0.68) -- (axis cs:0,1.92);

            \addplot[domain=-1.32:1.32, samples=160, smooth, line width=1.3pt, black!85]
                ({0.38490018 * x^3}, {x^2});

            \draw[thick, teal!75!black, dashdotted, -{Stealth[length=2.2mm]}]
                (axis cs:0,0) -- (axis cs:0,1.62);
            \draw[thick, blue!75!black, dashdotted, -{Stealth[length=2.2mm]}]
                (axis cs:0,0) -- (axis cs:0,-0.52);

            \addplot[domain=0:1.48, samples=80, smooth, semithick, red!80!black, -{Stealth[length=2.1mm]}]
                ({0.38*x - 0.04*x^2}, {x + 0.08*x^2});
            \addplot[domain=0:1.48, samples=80, smooth, semithick, red!80!black, -{Stealth[length=2.1mm]}]
                ({-0.38*x + 0.04*x^2}, {x + 0.08*x^2});

            \addplot[domain=0:1.10, samples=60, smooth, semithick, blue!65!black, densely dashed, -{Stealth[length=2mm]}]
                ({0.88*x + 0.05*x^2}, {0.42*x - 0.04*x^2});
            \addplot[domain=0:1.10, samples=60, smooth, semithick, blue!65!black, densely dashed, -{Stealth[length=2mm]}]
                ({-0.88*x - 0.05*x^2}, {0.42*x - 0.04*x^2});

            \filldraw[red!80!black] (axis cs:0.253,0.761) circle (1.6pt);
            \filldraw[red!80!black] (axis cs:-0.253,0.761) circle (1.6pt);
            \draw[thin, red!75!black, <-] (axis cs:0.28,0.75) -- (axis cs:0.46,0.65);
            \node[anchor=west, font=\scriptsize, text=red!75!black, fill=blue!3, inner sep=1pt] at (axis cs:0.47,0.64)
                {transit $t_{\mathrm{trans}}(\vartheta)$};

            \draw[semithick, purple!80!black, densely dotted, -{Stealth[length=1.7mm]}]
                (axis cs:0.55,1.27) ellipse [x radius=0.13, y radius=0.11];
            \node[anchor=west, font=\scriptsize, text=purple!80!black, fill=blue!3, inner sep=1pt] at (axis cs:0.70,1.27)
                {$S_3$ monodromy};

            \filldraw[fill=white, draw=black, thick] (axis cs:0,0) circle (2.2pt);
            \draw[thin, black!70, <-] (axis cs:0.04,-0.04) -- (axis cs:0.20,-0.24);
            \node[anchor=west, font=\scriptsize, fill=blue!3, inner sep=1pt] at (axis cs:0.21,-0.26)
                {$\mu_0 = A_3$ ($A_3^\pm$ cusp)};

            \node[font=\scriptsize, text=teal!55!black, align=center, fill=teal!12, inner sep=1pt] at (axis cs:0,1.18)
                {$\mathbb{R}\oplus\mathbb{R}\oplus\mathbb{R}$\\($\Delta_{\mathrm{alg}}>0$)};
            \node[anchor=south, font=\scriptsize, text=teal!75!black] at (axis cs:0,1.62)
                {$\vartheta = +\frac{\pi}{2}$};
            \node[anchor=east, font=\scriptsize, text=blue!75!black] at (axis cs:-0.04,-0.42)
                {$\vartheta = -\frac{\pi}{2}$};
            \node[font=\scriptsize, text=blue!60!black, align=center, fill=blue!3, inner sep=1pt] at (axis cs:-0.88,0.92)
                {$\mathbb{R}\oplus\mathbb{C}$\\($\Delta_{\mathrm{alg}}<0$, $A_2$ folds)};
            \draw[thin, black!70, <-] (axis cs:-0.48,1.16) -- (axis cs:-0.65,1.42);
            \node[anchor=south, font=\scriptsize, text=black!85, fill=blue!3, inner sep=1pt] at (axis cs:-0.72,1.42)
                {$\Delta_{\mathrm{alg}}=0$: $\mathbb{R}\oplus\mathbb{R}[\epsilon]/(\epsilon^2)$};
        \end{axis}

        \begin{axis}[
            name=panelB,
            at={(panelA.below south west)},
            anchor=north west,
            yshift=-0.55cm,
            width=0.57\textwidth,
            height=0.42\textwidth,
            xmin=-0.55, xmax=1.82,
            ymin=-12.5, ymax=13.5,
            clip=true,
            clip mode=individual,
            axis lines=box,
            xlabel={Dimensionless ratio $q = c_1 / u_2$},
            ylabel={Curvature $K_{(1)}(q)$},
            title={\textbf{(b)} Bundle curvature $K_{(1)}(q)$ \& horizons},
            title style={font=\small},
            tick label style={font=\footnotesize},
            label style={font=\small},
            restrict y to domain=-18:18
        ]
            \fill[teal!7]    (axis cs:-0.55,-12.5) rectangle (axis cs:1.0,13.5);
            \fill[orange!10] (axis cs:1.0,-12.5)   rectangle (axis cs:1.3333,13.5);
            \fill[gray!6]    (axis cs:1.3333,-12.5) rectangle (axis cs:1.82,13.5);

            \draw[gray!60, dotted] (axis cs:-0.55,0) -- (axis cs:1.82,0);

            \addplot[domain=-0.55:0.865, samples=160, smooth, very thick, blue!75!black]
                {-(18*x^4 - 81*x^3 + 147*x^2 - 149*x + 71) / (4 * (4 - 3*x)^2 * (1 - x))};
            \addplot[domain=1.088:1.298, samples=100, smooth, very thick, blue!75!black]
                {-(18*x^4 - 81*x^3 + 147*x^2 - 149*x + 71) / (4 * (4 - 3*x)^2 * (1 - x))};
            \addplot[domain=1.378:1.82, samples=100, smooth, very thick, blue!75!black]
                {-(18*x^4 - 81*x^3 + 147*x^2 - 149*x + 71) / (4 * (4 - 3*x)^2 * (1 - x))};

            \draw[thick, black!75, dashed] (axis cs:1.0,-12.5) -- (axis cs:1.0,13.5);

            \draw[thick, red!80!black, dashdotted] (axis cs:1.3333,-12.5) -- (axis cs:1.3333,13.5);

            \draw[thin, black!75, <-] (axis cs:0.98, 5.2) -- (axis cs:0.68, 5.2);
            \node[anchor=east, font=\scriptsize, align=right, fill=teal!7, inner sep=1.5pt] at (axis cs:0.67, 5.2)
                {Trace horizon $q=1$\\($\mathbb{Z}_2$ Maslov jump)};

            \draw[thin, red!80!black, <-] (axis cs:1.31, 10.2) -- (axis cs:0.86, 10.2);
            \node[anchor=east, font=\scriptsize, text=red!80!black, align=right, fill=teal!7, inner sep=1.5pt] at (axis cs:0.85, 10.2)
                {Metric wall $q=\frac{4}{3}$\\($K_{(1)}\to -\infty$)};

            \filldraw[fill=white, draw=blue!75!black, thick] (axis cs:0,-1.1094) circle (2pt);
            \node[anchor=north, font=\scriptsize, text=blue!75!black, fill=teal!7, inner sep=1pt] at (axis cs:-0.02,-2.1)
                {$K_{(1)}(0) = -\frac{71}{64}$};

            \node[anchor=south west, font=\scriptsize, text=teal!55!black] at (axis cs:-0.50,-11.8)
                {Defocusing ($q<1$): $K_{(1)}<0$};
        \end{axis}

        \begin{axis}[
            name=panelC,
            at={(panelB.south east)},
            anchor=south west,
            xshift=1.05cm,
            width=0.35\textwidth,
            height=0.42\textwidth,
            xmin=1.02, xmax=1.46,
            ymin=-0.85, ymax=0.85,
            clip mode=individual,
            axis lines=box,
            xtick={1.0, 1.2, 1.3333},
            xticklabels={$1$, $1.2$, $\frac{4}{3}$},
            ytick={-0.5, 0, 0.5},
            xlabel={Ratio $q = c_1 / u_2$},
            ylabel={$\chi = \log|u_2|$},
            ylabel style={yshift=-4pt},
            title={\textbf{(c)} Kossowski $A_2$ cusps},
            title style={font=\small},
            tick label style={font=\footnotesize},
            label style={font=\small}
        ]
            \fill[orange!10] (axis cs:1.02,-0.85) rectangle (axis cs:1.3333,0.85);
            \fill[gray!6]    (axis cs:1.3333,-0.85) rectangle (axis cs:1.46,0.85);

            \draw[thick, red!80!black, dashdotted] (axis cs:1.3333,-0.85) -- (axis cs:1.3333,0.85);

            \addplot[domain=-0.52:0.52, samples=80, smooth, semithick, orange!85!black, -{Stealth[length=1.8mm]}]
                ({1.3333 - x^2}, {0.26 + 3.6*x^3});
            \addplot[domain=-0.52:0.52, samples=80, smooth, semithick, purple!80!black, densely dashed, -{Stealth[length=1.8mm]}]
                ({1.3333 - x^2}, {-0.26 + 3.6*x^3});

            \filldraw[red!80!black] (axis cs:1.3333, 0.26) circle (1.6pt);
            \filldraw[red!80!black] (axis cs:1.3333,-0.26) circle (1.6pt);

            \node[anchor=south west, font=\scriptsize, align=left, fill=orange!10, inner sep=1.5pt] at (axis cs:1.04,-0.78)
                {$P_\chi \neq 0$ rays\\cusp at $q=\frac{4}{3}$};
        \end{axis}
    \end{tikzpicture}
    \caption{Schematic summary of moduli-plane and global bundle singularities for the $d = 2$ geodesic spray launched from $A_3 = \mathbb{R}[x]/(x^3)$. \textbf{(a)} On the Tschirnhaus moduli plane $(\tilde{c}_0, \tilde{c}_1)$, Arnol'd's semicubical cusp discriminant $\Delta_{\mathrm{alg}} = 4\tilde{c}_1^3 - 27\tilde{c}_0^2 = 0$ (heavy solid curve) separates $\mathbb{R}\oplus\mathbb{R}\oplus\mathbb{R}$ (shaded wedge) from $\mathbb{R}\oplus\mathbb{C}$; shallow-angle rays (solid curves) initially enter $\mathbb{R}\oplus\mathbb{C}$ before crossing the discriminant wall at $t_{\mathrm{trans}}(\vartheta)$ into $\mathbb{R}\oplus\mathbb{R}\oplus\mathbb{R}$, whereas wider oblique rays (dashed) remain in $\mathbb{R}\oplus\mathbb{C}$, vertical $\sigma$-invariant rays $\vartheta = \pm\frac{\pi}{2}$ (dash-dotted) stay within a single stratum, and complexified loops around $\Delta_{\mathrm{alg}} = 0$ (dotted) carry $S_3$ monodromy. \textbf{(b)} Along the $\sigma$-fixed subbundle ($q = c_1/u_2$), the Gaussian curvature $K_{(1)}(q)$ in~(\ref{eq:A3_K}) is strictly negative across the Riemannian defocusing chamber $q < 1$, diverges across the infinite-distance trace horizon $q = 1$ (dashed vertical line), and diverges to $-\infty$ at the internal metric-degeneracy wall $q = \frac{4}{3}$ (dash-dotted vertical line). \textbf{(c)} On the $(q, \chi)$ surface, trajectories ($P_\chi \neq 0$, $2E + \frac{1}{3}P_\chi^2 > 0$) reaching $q = \frac{4}{3}$ in $1 < q < \frac{4}{3}$ form finite-time Kossowski semicubical cusps $(q - \frac{4}{3})^3 \propto -(\chi - \chi_*)^2$.}
    \label{fig:A3_spray_catastrophes}
\end{figure}

\section{Future Directions}
The geometric and variational frameworks developed here for finite-dimensional unital associative algebras over $\R$ point to several natural directions for further investigation even if we remain in the real algebra finite-dimensional algebra class. 

For example, in the deformation-bundle setting of Section~\ref{sec:def_metric}, both low-dimensional commutative slices examined in Sections~\ref{sec:2d_deformations} and~\ref{sec:A3_spray} exhibit strictly negative Gaussian curvature ($K < 0$) across their Riemannian chambers, rendering radial geodesic sprays defocusing up to the trace horizon. Characterizing which symmetric Frobenius algebras admit positively curved ($K > 0$) focusing chambers that generate genuine Jacobi conjugate caustics prior to horizon termination---alongside lifting deformations of non-Frobenius algebras via their trivial extensions $B \ltimes B^*$ and classifying bundle caustic bifurcations, Kossowski metric-degeneracy cusps, and trace-horizon obstructions across broader noncommutative classes---remains open.

Also, the particular separation results of Section~3 suggest sharper invariance questions. In Section~\ref{sec:homological_separation}, $A_{\mathrm{hyp}}$ and $A_{\mathrm{ell}}$ share identical radical-filtration and metric-signature data yet differ in their Hochschild Gerstenhaber structures, whereas in Section~\ref{sec:derived_separation}, $A_{\mathrm{star}}$ and $A_{\mathrm{line}}$ are derived equivalent (sharing all Hochschild Gerstenhaber and cyclic invariants) yet differ in radical nilpotency index. Can ART invariants separate non-isomorphic algebras when radical-filtration data and derived equivalence agree simultaneously?

A closely related homological question concerns the one-sided Yoneda algebra $\operatorname{Ext}_A^\bullet(S, S)$ of $S = A/\operatorname{rad}(A)$, which is not a derived invariant and distinguishes split bound quiver (including Brauer tree) algebras via their quiver presentations. For split bound quiver algebras whose minimal relations all lie in $J^3$ or higher powers of the arrow ideal $J$, however, the binary Yoneda product on $\operatorname{Ext}_A^1(S, S)$ vanishes, relegating those relations to higher $A_\infty$ Massey products $(m_3, m_4, \dots)$. Because the ART Legendre transform instead captures higher-order relation and socle data directly inside higher-jet phase germs via finite-dimensional calculus, two natural problems emerge: clarifying the dictionary between higher-jet ART invariants and $A_\infty$ Massey operations, and determining whether ART caustics can separate non-isomorphic symmetric Frobenius algebras that share isomorphic Hochschild Gerstenhaber algebras, cyclic (co)homologies, and associative Yoneda rings $\operatorname{Ext}_A^\bullet(S, S)$ simultaneously. In parallel, lifting $A$ across the canonical Cuntz--Quillen tower of infinitesimal thickenings $A_{(k)} = R/I^{k+1}$ of Remark~\ref{rem:CQ} (passing to $A_{(k)} \ltimes A_{(k)}^*$ as needed) restores trace sensitivity to noncommutative relations hidden inside $[A, A]$ at baseline ($k = 0$), and may likewise resolve isomorphism ambiguities left open by passing from a non-Frobenius algebra $B$ to its single trivial extension $B \ltimes B^*$.

But looking beyond real finite-dimensional algebras opens several broader frontiers. The most immediate extension is to finite-dimensional complex algebras. Viewing an $n$-dimensional complex algebra $A$ as a $2n$-dimensional real algebra $A_{\R}$ allows any non-degenerate $\C$-linear trace $\tau_{\C}$ with $\operatorname{Re}\tau_{\C}(\mathbf{1}) \neq 0$ (achievable after complex scalar rescaling whenever $\tau_{\C}(\mathbf{1}) \neq 0$) to induce an admissible real trace $\tau_{\R} \coloneqq \operatorname{Re}\tau_{\C} \in \mathcal{T}_1(A_{\R})$, so that our real pseudo-Finsler machinery applies directly (passing to the trivial extension if needed). Alternatively, formulating the ART potential and resulting trace-parametrized structures holomorphically on $T^*A \cong \C^{2n}$ links complex associative algebras directly to complex Lagrangian singularity theory.  Comparing the invariants arising from these real and holomorphic viewpoints warrants investigation.

Extending the framework to the infinite-dimensional topological and differential operator algebras central to noncommutative geometry, cyclic cohomology, and deformation quantization~\cite{kontsevich:2003} is an intriguing but substantially harder prospect. Beyond analytic domain and Fr\'echet-manifold issues, one must replace finite-dimensional linear traces with regularized or Dixmier-type traces~\cite{LSZ2013} while preserving a workable Legendre transform, though successful extensions along the lines of Sections~3 and~4 could carry significant payoff.

Finally, one may ask whether an ART-like program extends to algebras over fields of positive characteristic. Over a finite field $\mathbb{F}_q$, the discrete topology precludes continuous path integration and real-analytic wave propagation, and arbitrary trace-weighted determinant powers on a semisimple complement $S$ (cf.\ Proposition~\ref{prop:WM_factorization}) have no global analogue. Nevertheless, algebraic analogues of differential and integral calculus allow the local jet core of the ART framework to be formulated over $\mathbb{F}_q$~\cite{Greuel2016}.

Specifically, on the algebraic group scheme of units $\underline{A}^\times$ over $\mathbb{F}_q$, the Maurer--Cartan $1$-form $t^{-1}dt$ is well-defined within the algebraic de~Rham complex of K\"ahler differentials $\Omega^\bullet_{\underline{A}^\times / \mathbb{F}_q}$. Writing $t = \mathbf{1} + u$ in the formal $k$-jet neighborhood $\widehat{\mathcal{O}}_{\underline{A}^\times, \mathbf{1}} / \mathfrak{m}_{\mathbf{1}}^{k+1}$ of $\mathbf{1}$ (encompassing both semisimple and radical perturbations), trace cyclicity yields $d\big(\tau(u^m)\big) = m\,\tau(u^{m-1}du)$ for every $\tau \in \mathcal{T}(A)$ and $m \ge 1$. Consequently, whenever $\mathcal{T}_1(A) \neq \emptyset$ and $p = \operatorname{char}(\mathbb{F}_q) > k \ge 2$, every denominator and factorial up to $k$ is invertible in $\mathbb{F}_q$, so $\frac{1}{\tau(\mathbf{1})}\tau(t^{-1}dt)$ admits the formal $k$-jet antiderivative $\Psi_k(\mathbf{1} + u; \tau) = \frac{1}{\tau(\mathbf{1})}\sum_{m=1}^k \frac{(-1)^{m-1}}{m}\tau(u^m)$ and truncated ART $k$-jet $F_k(\mathbf{1} + u; \tau) = \sum_{m=0}^k \frac{1}{m!}\Psi_k(\mathbf{1} + u; \tau)^m \pmod{\mathfrak{m}_{\mathbf{1}}^{k+1}}$ (with $\Psi_k$ stabilizing on $\mathbf{1} + \operatorname{rad}(A)$ to the exact polynomial $P_\tau$ once $k \ge N - 1$ for radical nilpotency index $N$).

Furthermore, since $p > k \ge 2$ implies $p \neq 2$, both the fundamental tensor $g_{\mathbf{1}}(u, v; \tau)$ in~(\ref{eq:g_one}) and the bundle metric $\mathcal{G}_{(\mu, \tau_\mu)}$ in~(\ref{eq:bundle_metric_G}) are well-defined rational expressions wherever $\Delta_\mu(\tau) \neq 0$. Formal Legendre transforms, $k$-jet phase reversions (including cases with $k \ge N$, such as $A_3^\pm$ on $\mathbb{F}_q[x]/(x^3)$ for $p \ge 5$), and local Milnor algebras $Q_f = \mathbb{F}_q[[\boldsymbol{\sigma}]] / \langle \nabla f \rangle$ can therefore be analyzed directly in this characteristic range~\cite{Greuel2016}. It remains to be seen whether this will lead to novel insights or practical utility.

\end{document}